\documentclass[12pt,amscd]{amsart}
\usepackage{amssymb}
\numberwithin{equation}{section}
\newtheorem{thm}[equation]{Theorem}

\newtheorem{cor}[equation]{Corollary}

\newtheorem{lem}[equation]{Lemma}

\newtheorem{prop}[equation]{Proposition}

\newenvironment{pf}{\proof[\proofname]}{\endproof}
\newenvironment{pf*}[1]{\proof[#1]}{\endproof}

\theoremstyle{definition}

\newtheorem{defn}[equation]{Definition}

\theoremstyle{remark}

\newtheorem*{rmk}{Remark}
\newtheorem*{rmks}{Remarks}

\newtheorem*{ack}{Acknowledgement}

\newcommand{\comment}[1]{}

\makeatother

\begin{document}
\baselineskip=18truept


\def\C {{\mathbb C}}
\def\Cn {{\mathbb C}^n}
\def\R {{\mathbb R}}
\def\Rn {{\mathbb R}^n}
\def\Z {{\mathbb Z}}
\def\N {{\mathbb N}}
\def\cal#1{{\mathcal #1}}
\def\bb#1{{\mathbb #1}}

\def\dbar {\bar \partial }
\def\dir {{\mathcal D}}
\def\lev#1{{\mathcal L}(#1)}
\def\lap {\Delta }
\def\ol {{\mathcal O}}
\def\E {{\mathcal E}}
\def\J {{\mathcal J}}
\def\U {{\mathcal U}}
\def\V {{\mathcal V}}
\def\z {\zeta }
\def\Harm {\text {Harm}\, }
\def\grad {\nabla }
\def\dexh {\{ M_k \} _{k=0}^{\infty } }
\def\sing#1{#1_{\text {sing}}}
\def\reg#1{#1_{\text {reg}}}

\def\lquotient{\big{\backslash}}

\def\setof#1#2{\{ \, #1 \mid #2 \, \} }

\def\image#1{\text{\rm im}\,\bigl[#1\bigr]}
\def\kernel#1{\text{\rm ker}\,\bigl[#1\bigr]}

\def\holecl {M\setminus \overline M_0}
\def\hole {M\setminus M_0}

\def\nd{\frac {\partial }{\partial\nu } }
\def\ndof#1{\frac {\partial#1}{\partial\nu } }

\def\pdof#1#2{\frac {\partial#1}{\partial#2}}

\def\cinf{\ensuremath{C^{\infty}} }
\def\cinfns{\ensuremath{C^{\infty}}}

\def\diam{\text {\rm diam} \, }

\def\real{\text {\rm Re}\, }

\def\imag{\text {\rm Im}\, }

\def\supp{\text {\rm supp}\, }

\def\Vol{\text {\rm vol} \, }

\def\restrict#1{\upharpoonright_{#1}}

\def\sm{\setminus}



\def\anal{analytic }
\def\analns{analytic}

\def\bdd{bounded }
\def\bddns{bounded}

\def\cpt{compact }
\def\cptns{compact}

\def\cpx{complex }
\def\cpxns{complex}

\def\cont{continuous }
\def\contns{continuous}

\def\dime{dimension }
\def\dimens{dimension }

\def\exh{exhaustion }
\def\exhns{exhaustion}

\def\fn{function }
\def\fnns{function}

\def\fns{functions }
\def\fnsns{functions}

\def\holo{holomorphic }
\def\holons{holomorphic}

\def\mero{meromorphic }
\def\merons{meromorphic}

\def\holoconvex{holomorphically convex }
\def\holoconvexns{holomorphically convex}

\def\ircomp{irreducible component }
\def\concomp{connected component }
\def\ircompns{irreducible component}
\def\concompns{connected component}
\def\ircomps{irreducible components }
\def\concomps{connected components }
\def\ircompsns{irreducible components}
\def\concompsns{connected components}

\def\irred{irreducible }
\def\irredns{irreducible}

\def\con{connected }
\def\conns{connected}

\def\comp{component }
\def\compns{component}
\def\comps{components }
\def\compsns{components}

\def\mfld{manifold }
\def\mfldns{manifold}
\def\mflds{manifolds }
\def\mfldsns{manifolds}

\def\nbd{neighborhood }
\def\nbds{neighborhoods }
\def\nbdns{neighborhood}
\def\nbdsns{neighborhoods}

\def\harm{harmonic }
\def\harmns{harmonic}
\def\plh{pluriharmonic }
\def\plhns{pluriharmonic}
\def\plsh{plurisubharmonic }
\def\plshns{plurisubharmonic}

\def\qplsh#1{$#1$-plurisubharmonic}
\def\hplsh{$(n-1)$-plurisubharmonic }
\def\hplshns{$(n-1)$-plurisubharmonic}

\def\para{parabolic }
\def\parans{parabolic}

\def\rel{relatively }
\def\relns{relatively}

\def\str{strictly }
\def\strns{strictly}

\def\strg{strongly }
\def\strgns{strongly}

\def\cvx{convex }
\def\cvxns{convex}

\def\wrt{with respect to }
\def\wrtns{with respect to}

\def\st {such that }
\def\stns {such that}

\def\hm {harmonic measure }
\def\hmns {harmonic measure}

\def\hmib {harmonic measure of the ideal boundary of }
\def\hmibns {harmonic measure of the ideal boundary of}

\def\Vert{\text{{\rm Vert}}\, }
\def\Edge{\text{{\rm Edge}}\, }

\def\til#1{\tilde{#1}}
\def\wtil#1{\widetilde{#1}}

\def\what#1{\widehat{#1}}

\def\seq#1#2{\{#1_{#2}\} }


\def\vphi {\varphi }


\def\inv{   ^{-1}  }

\def\ssp#1{^{(#1)}}

\def\set#1{\{ #1 \}}

\title[Filtered ends and proper holomorphic mappings]
{Filtered ends, proper holomorphic mappings of K\"ahler manifolds
to Riemann surfaces, and K\"ahler groups}
\author[T.~Napier]{Terrence Napier$^{*}$}
\address{Department of Mathematics\\Lehigh University\\Bethlehem, PA 18015}
\email{tjn2@lehigh.edu}
\thanks{$^{*}$Research partially
supported by NSF grant DMS0306441}
\author[M.~Ramachandran]{Mohan Ramachandran}
\address{Department of Mathematics\\SUNY at Buffalo\\Buffalo, NY 14260}
\email{ramac-m@math.buffalo.edu}

\subjclass[2000]{32C17} \keywords{Fundamental groups, potential
theory}

\date{July 6, 2006}

\begin{abstract}
The main result of this paper is that a \con bounded geometry
complete K\"ahler manifold which has at least~$3$ filtered ends
admits a proper \holo mapping onto a Riemann surface.  As an
application, it is also proved that any properly ascending HNN
extension with finitely generated base group, as well as
Thompson's groups V, T, and F, are not K\"ahler. The results and
techniques also yield a different proof of the theorem of Gromov
and Schoen that, for a \con \cpt K\"ahler manifold whose
fundamental group admits a proper amalgamated product
decomposition, some finite unramified cover admits a surjective
\holo mapping onto a curve of genus at least~$2$.

This version of this paper contains details not in the version
submitted for publication.
\end{abstract}

\maketitle

\section*{Introduction} \label{introduction}

\noindent \textit{This is a version of a paper which is similar to
another paper of the same title submitted for publication. This
version contains details not in the version submitted for
publication.}

The main goal of this paper is the following (the required
definitions appear later in this introduction):
\begin{thm}\label{Main Theorem from Introduction}
Let $X$ be a \con complete K\"ahler manifold satisfying at least
one of the following hypotheses:
\begin{enumerate}
\item [(i)] $X$ has \bdd geometry, \item[(ii)] $X$ admits a
positive Green's \fn $G$ that vanishes at infinity, or
\item[(iii)] $X$ is weakly $1$-complete.
\end{enumerate}
If $\tilde e(X)\geq 3$, then $X$ admits a proper \holo
mapping onto a Riemann surface.
\end{thm}
In particular, this paper provides a unified framework for some of
the results of \cite{Gro1}, \cite{L}, \cite{Gro2}, \cite{GroS},
\cite{NR1}--\cite{NR4}, and \cite{DG}. This framework relies on
the notion of filtered ends, first introduced by Kropholler and
Roller \cite{KroR} in the group theoretic context and later given
a topological interpretation by Geoghegan~\cite{Ge} (we work with
Geoghegan's topological notion in this paper). Theorem~\ref{Main
Theorem from Introduction} for $X$ satisfying conditions (i) and
(ii) together was first proved by Delzant and Gromov~\cite{DG}
using harmonic maps into trees (as in \cite{GroS}, \cite{KoS1},
\cite{KoS2}, \cite{Sun}).  They applied their result to the
problem of determining which hyperbolic groups are K\"ahler.  In
this paper, we consider a different approach which yields the more
complete result Theorem~\ref{Main Theorem from Introduction} and
which is more elementary in the sense that it only uses \harm
\fnsns .  Cousin's example \cite{Co} of a $2$-ended weakly
$1$-complete covering of an Abelian variety which has only
constant \holo \fns demonstrates that one cannot weaken the
hypotheses to $\tilde e(X)\geq 2$. On the other hand, a slightly
stronger version (Theorem~\ref{General filtered ends theorem}) is
obtained in the case in which each end separately (in an
appropriate sense) has \bdd geometry, is weakly $1$-complete, or
admits a positive Green's \fn $G$ that vanishes at infinity.

Theorem~\ref{Main Theorem from Introduction} and elementary facts
from geometric group theory together give the following:
\begin{thm}\label{Gromov-Schoen Theorem}
Let $X$ be a \con \cpt K\"ahler manifold with fundamental group
$\Lambda=\pi_1(X)$ satisfying at least one of the following:
\begin{enumerate}
\item[(a)] {\rm (Gromov and Schoen \cite{GroS})} $\Lambda$ admits
a proper amalgamated product decomposition (i.e. \(\Lambda=\Gamma
_1* _\Gamma \Gamma _2\) where the index of $\Gamma $ in $\Gamma
_1$ is at least~$3$ and the index of $\Gamma $ in $\Gamma _2$ is
at least~$2$); or

\item[(b)] {\rm (See \cite{NR3})} $\Lambda$ is a properly
ascending HNN extension.
\end{enumerate}
Then some finite (unramified) covering of $X$ admits a surjective
\holo mapping onto a curve of genus $g\geq 2$.
\end{thm}
\begin{rmks}
1. Conversely, if a \con \cpt manifold $M$  admits a surjective
\cont mapping onto a curve $S$ of genus $g\geq 2$, then $\pi
_1(M)$ admits a proper amalgamated product decomposition. For such
a decomposition exists for $\pi _1(S)$ by Van Kampen's theorem and
one may pull this back to $\pi _1(M)$.

\noindent 2. By considering the action of $\Lambda$ on the
associated tree and applying Theorem~0.1, one gets the theorem in
both cases (a)~and~(b) simultaneously (see the proof of
Theorem~\ref{Action on tree version of Gromov-Schoen theorem}). On
the other hand, according to a theorem of Baumslag and Shalen (see
Theorem~6 of Chapter~4 of \cite{Baumslag Book}), a finitely
presented group which can be expressed as a properly ascending HNN
extension, but not one with finitely generated base group, is
virtually a proper amalgamated product. Thus
Theorem~\ref{Gromov-Schoen Theorem} is actually contained within
the theorem of Gromov and Schoen (i.e. the case~(a)) together with
part~(ii) of Theorem~\ref{HNN and Thompson not Kahler theorem from
intro} below and the theorem of Baumslag and Shalen.
\end{rmks}

Theorem~\ref{Main Theorem from Introduction} and its consequences
also lead to new restrictions on K\"ahler groups (i.e. fundamental
groups of \cpt K\"ahler manifolds).

\begin{thm}\label{HNN and Thompson not Kahler theorem from intro}
The following groups are not K\"ahler:
\begin{enumerate}
\item[(i)] Thompson's groups V, T, and F; and

\item[(ii)] Any properly ascending HNN extension with finitely
generated base group.
\end{enumerate}
\end{thm}
The question as to whether or not $F$ is K\"ahler was first posed
by Geoghegan (see \cite{Brown}) and the first proof that $F$ is
not K\"ahler appeared in \cite{NR4}. Since $F$ is a properly
ascending HNN extension with finitely generated base group, this
may now be viewed as a special case of part~(ii) of
Theorem~\ref{HNN and Thompson not Kahler theorem from intro}.
Daniel Farley has independently obtained the result that $V$ and
$T$ are not K\"ahler. For more on K\"ahler groups, the reader may
refer to \cite{Ar} and \cite{ABCKT}.

Before sketching the proof of Theorem~\ref{Main Theorem from
Introduction}, we make some remarks which put these results in
context and we recall the required definitions.  For $X$ a \con
\cpt K\"ahler manifold, a natural and much studied problem is to
determine when $X$ admits a surjective \holo mapping onto a curve
of genus $g\geq 2$. According to the classical theorem of
Castelnuovo and de~Franchis (see \cite{Be}, \cite{BPV}), this is
the case if and only if there exists a pair of linearly
independent \holo $1$-forms $\omega _1,\omega _2$ \st $\omega
_1\wedge \omega _2\equiv 0$. The main point is that the
meromorphic \fn $f\equiv \omega _1/\omega _2$ actually has no
points of indeterminacy (this is a general fact about closed \holo
$1$-forms; see, for example, \cite{NR2} for an elementary proof).
Stein factorization then gives the required mapping. For some of
the many other results in this context, the reader may refer to
the work of Beauville (see \cite{Cat1}), \cite{CarT}, \cite{Siu},
\cite{Simpson-VHS}, \cite{GroS}, \cite{JY1}, \cite{JY2},
\cite{Ar}, and \cite{ABCKT}.

For a \con non\cpt complete K\"ahler manifold $(X,g)$, the
analogous problem is to determine when $X$ admits a proper \holo
mapping onto a Riemann surface. We will mainly consider the case
in which $X$ has bounded geometry of order~$k\geq 2$ (in the sense
that there exists a constant $C>0$ and, for each point $p\in X$, a
biholomorphism $\Psi $ of the unit ball $B=B(0;1)\subset \C ^n$
onto a \nbd $U$ of $p$ in $X$ \st $\Psi (0)=p$ and, on $B$, $C\inv
g_{\C ^n} \leq \Psi ^*g\leq Cg_{\C ^n}$ and $| D^m\Psi ^*g| \leq C
\text{ for } m=0,1,2,\dots , k$); the case in which $X$ is weakly
$1$-complete (i.e.~$X$ admits a \cont \plsh \exh \fnns ); or the
case in which $X$ admits a positive   Green's \fn which vanishes
at infinity. As shown in \cite{Gro1}, \cite{Gro2}, \cite{L}, and
\cite{NR1}, if such an $X$ has at least $3$ ends, then $X$ admits
a proper \holo mapping onto a Riemann surface. The main step is to
produce \plh \fns $\rho _1,\rho _2$ which have different limits at
infinity along the various ends and whose \holo differentials
$\omega _1=\partial \rho _1$, $\omega _2=\partial \rho _2$ satisfy
$\omega _1\wedge \omega _2\equiv 0$. In particular,  $1$,~$\rho
_1$,~$\rho _2$, and hence $d\rho _1$,~$d\rho _2$, are linearly
independent. If $\omega _1$, $\omega _2$ are linearly independent,
then one Stein factors the \holo map
$$
f=\frac {\omega _1}{\omega _2} :X\to {\bb P} ^1.
$$
Otherwise, one gets a \holo \fn $f=\rho _1+c\rho _2$, for some
constant $c\in \C $, which one may Stein factor. Thus it has been
known for some time that the ends structure is relevant to the
problem of finding a proper \holo mapping onto a Riemann surface.

We now recall the definitions of ends and filtered ends. Depending
on the context, by an {\it end} of a \con manifold $M$, we will
mean either a \comp $E$ of $M\setminus K$ with non\cpt closure,
where $K$ is a given \cpt subset of $M$, or an element of
$$
\lim _{\leftarrow } \pi _0 (M\setminus K)
$$
where the limit is taken as $K$ ranges over the \cpt subsets of
$M$ (or the \cpt subsets of $M$ whose complement $M\setminus K$
has no \rel \cpt \compsns). The number of ends of $M$ will be
denoted by $e(M)$. For a \cpt set $K$ \st $M\setminus K$ has no
\rel \cpt \compsns , we get an {\it ends decomposition}
$$
M\setminus K=E_1\cup \cdots \cup E_m,
$$
where $E_1, \dots , E_m$ are the distinct \comps of $M\setminus
K$.

As in the work of Geoghegan \cite{Ge}, for $\Upsilon :\widetilde M
\to M$ the universal covering of $M$, consider the set
$$
\lim _{\leftarrow } \pi _0 [\Upsilon \inv (M\setminus K)],
$$
where the limit is taken as $K$ ranges over the \cpt subsets of
$M$ (or the \cpt subsets of $M$ whose complement $M\setminus K$
has no \rel \cpt \compsns). Following \cite{Ge}, we will call
elements of the above set {\it filtered ends}.  The number of
filtered ends of $M$ will be denoted by $\tilde e(M)$. Clearly,
$\tilde e(M)\geq e(M)$. In fact, for $k\in \N$, we have $\tilde
e(M)\geq k$ if and only if there exists an ends decomposition
$M\setminus K=E_1\cup \cdots \cup E_m$ for $M$ \stns , for $\Gamma
_j=\image{\pi _1(E_j) \to \pi _1(M)}$ for $j=1,\dots , m$, we have
$$
\sum _{j=1}^m[\pi _1(M):\Gamma _j]\geq k.
$$
Moreover, if $\widehat M\to M$ is a \con covering space, then
$\tilde e(\widehat M)\leq \tilde e(M)$ with equality if the
covering is finite.

To illustrate some of the arguments in the proof of
Theorem~\ref{Main Theorem from Introduction}, let us consider the
case in which $e(X)=2$ and $X$ admits a positive   Green's \fn $G$
that vanishes at infinity. In this case, $X$ admits an ends
decomposition $X\setminus K=E_1\cup E_2$ \st the image $\Gamma $
of $\pi _1(E_1)$ in $\pi _1(X)$ is a proper subgroup. By standard
arguments, there exists a \plh \fn $\rho :X \to (0,1)$ with finite
energy \st
$$
\lim _{x\to \infty } \rho \restrict{\overline E_1}(x)=1 \qquad
\text{and} \qquad \lim _{x\to \infty } \rho \restrict{\overline
E_2}(x)=0.
$$
In particular, $X$ is weakly $1$-complete. Taking $\Upsilon
:\widehat X\to X$ to be a \con covering space (not the universal
covering) with $\Upsilon _* \pi _1(\widehat X)=\Gamma $, we see
that $\Upsilon $ maps some \comp $\Omega _1$ of $\Upsilon \inv
(E_1)$ isomorphically onto $E_1$ and the set $\Omega _2=\Upsilon
\inv (E_1)\setminus \Omega _1\neq \emptyset $. Again, there exists
a \plh \fn $\rho _2:\widehat X \to (0,1)$ with finite energy \st
$$
\lim _{x\to \infty } \rho _2\restrict{\overline \Omega _1}(x)=1
\qquad \text{and} \qquad \liminf _{x\to \infty } \rho
_2\restrict{\widehat X\setminus \Omega _1}(x)=0.
$$
If $\Gamma $ is of finite index, then $e(\widehat X)\geq 3$ and
hence, by \cite{NR1}, $\widehat X$ admits a proper \holo mapping
onto a Riemann surface. Since $\widehat X\to X$ is a finite
covering in this case, $X$ also admits such a mapping.  If $\Gamma
$ is of infinite index, then the lift $\rho _1=\rho \circ \Upsilon
$ does not have finite energy and so $d\rho _1 $ and $d\rho _2$
must be linearly independent. On the other hand, since $\rho _1$
and $\rho _2$ have \cpt levels in $\Omega _1$ over values
near~$1$, we must have $\partial\rho _1\wedge\partial\rho _2\equiv
0$ (see, for example, \cite{NR3}, Lemma~2.1). It follows that some
(nonempty) open subset of $\Omega _1\cong E_1$ admits a proper
\holo mapping onto a Riemann surface. Standard arguments now imply
that this is the case for $X$.

For the general case, we will again pass to the appropriate
covering spaces. We will produce suitable \plh \fns $\rho _1$ and
$\rho _2$ with prescribed values at infinity along {\it filtered}
ends by applying the theory of massive sets as in Grigor'yan
\cite{Gri}. A version of the cup product lemma (see Lemma~\ref{Cup
product for sublevel of plh function lemma} below) will give
$\partial\rho _1\wedge\partial\rho _2\equiv 0$ (Gromov \cite{Gro2}
was the first to notice the cup product lemma in the context of
bounded geometry and subsequent refinements were formed and used
by others). The idea of filtered ends was applied in some special
cases in the context of Lefschetz type theorems in \cite{NR1},
\cite{NR2}, and \cite{NR3}. The arguments in this paper formalize
this approach.

According to a theorem of Simpson~\cite{Simpson-Lefschetz thm
integral leaves}, if $\Upsilon\colon\what X\to X$ is a \con
covering space of a smooth projective variety $X$, $\rho$ is a
nonconstant \plh \fn on $\what X$ \st $\partial\rho$ descends to a
\holo $1$-form $\alpha$ on $X$ and \stns, for any $\zeta\in\C$,
the fiber $F=\rho\inv(\zeta)$ satisfies
\[
\sum_{L\text{ a component of F}}\left[\pi_1(\what
X):\image{\pi_1(L)\to\pi_1(\what X)}\right]>1,
\]
then there exists a surjective \holo mapping $\Phi$ of $X$ onto a
curve $Y$ and a \holo $1$-form $\beta$ on $Y$ with
$\alpha=\Phi^*\beta$ (Simpson also obtains a version for $\alpha$
the push-forward of $d\rho$ for a real-valued \plh \fn $\rho$).
Simpson's result and Theorem~\ref{Main Theorem from Introduction}
appear to be related, but neither is known to imply the other.

Another application of Theorem~\ref{Main Theorem from
Introduction}, which will be applied in the proof of part~(i) of
Theorem~\ref{HNN and Thompson not Kahler theorem from intro}, is
the following (cf. Theorem~\ref{General special ends, two
ends-infinitely generated maps to Riemann surface theorem}):
\begin{thm}\label{Two ends and infinitely generated maps to Riemann surface theorem}
Let $X$ be a \con non\cpt complete K\"ahler manifold which has
\bdd geometry of order~$2$ or which is weakly $1$-complete or
which admits a positive Green's \fn $G$ that vanishes at infinity.
Assume that $\tilde e(X)\geq 2$ and $\pi _1(X)$ is infinitely
generated. Then $X$ admits a proper \holo mapping onto a Riemann
surface.
\end{thm}
This is an immediate consequence of Theorem~\ref{Main Theorem from
Introduction} and the following fact:
\begin{lem}\label{Two ends and infinitely generated lemma}
Let $M$ be a \con $\cinf$ manifold \st $\tilde e(M)\geq 2$ and
$\Lambda =\pi _1(M)$ is infinitely generated. Then $\tilde
e(M)=\infty$.
\end{lem}
\begin{pf}
If $e(M)=1$, then the image $\Theta$ of the fundamental group of
some end $E=M\setminus K$ is a proper subgroup of $\Lambda$. If
$[\Lambda:\Theta]=\infty$, then $\tilde e(M)=\infty$ as claimed.
If $[\Lambda:\Theta]<\infty$, then $\Theta$ is infinitely
generated and the finite covering $\widehat M\to M$ with
$\image{\pi_1(\widehat M)\to\Lambda}=\Theta$ satisfies $e(\widehat
M)\geq 2$ and $\tilde e(\widehat M)=\tilde e(M)$. Thus we may
assume without loss of generality that $e(M)\geq 2$.

We may fix a $\cinf $ \rel \cpt domain $\Omega$ in $M$ \st
$M\setminus\overline\Omega$ has exactly two \comps $E_0'$ and
$E_1'$, each with non\cpt closure. Fix a point $x_0\in\Omega$, let
$\Gamma\equiv\image{\pi _1(\Omega ,x_0)\to\pi _1(M,x_0)}$, and,
for $i=0,1$, let $E_i$ be the end defined by $E_i=\overline
{E_i'}\cup\Omega$ and let $\Gamma _i\equiv\image{\pi
_1(E_i,x_0)\to\pi _1(M,x_0)}$. Then $\Gamma _0\cap\Gamma
_1=\Gamma$ (for example, by Van Kampen's theorem), so
$$
[\Lambda :\Gamma _0]\geq [\Gamma _1:\Gamma _0\cap\Gamma
_1]=[\Gamma _1:\Gamma ].
$$

If $\Gamma _1$ is finitely generated, then $[\Lambda :\Gamma
_1]=\infty$ because $\Lambda$ is infinitely generated. Hence
$\tilde e(M)=\infty$ in this case. If $\Gamma _1$ is infinitely
generated, then, since $\Gamma$ is finitely generated, $[\Gamma
_1:\Gamma ]=\infty$. Hence $[\Lambda :\Gamma _0]=\infty$ by the above
inequality and,
again, $\tilde e(M)=\infty$.
\end{pf}

We close this section with a proof that Thompson's groups $V$ and
$T$ are {\it not} K\"ahler (part~(i) of Theorem~\ref{HNN and
Thompson not Kahler theorem from intro}). The proof that any
properly ascending HNN extension with finitely generated base
group (for example, Thompson's group $F$) is not K\"ahler
(part~(ii) of Theorem~\ref{HNN and Thompson not Kahler theorem
from intro}) will be given in Section~\ref{Amalgamations and
mappings to Riemann surfaces section}. The group $V$ is the group
of right-\cont bijections $\lambda\colon [0,1]\to [0,1]$ \st
$\lambda$ maps dyadic rational numbers to dyadic rational numbers,
$\lambda$ is differentiable except at finitely many dyadic
rational numbers, and $\lambda$ is affine with derivative a power
of~$2$ on each interval on which it is differentiable.  The group
$T$ is the subgroup consisting of all $v\in V$ \st $v$ induces a
homeomorphism of the circle $[0,1]/0\sim 1$, and $F$ is the
subgroup consisting of all homeomorphisms in $V$.  For the proof
that $V$ and $T$ are not K\"ahler, we will apply the following
fact:

\begin{lem}\label{V and T end of pair lemma}
If $M$ is a \con \cpt manifold with fundamental group $V$ or $T$,
then there exists a \con covering $\widehat M\to M$ \st $\tilde
e(\widehat M)=\infty$.
\end{lem}
\begin{pf}
For a subset $S\subset [0,1]$, let $V_S$ denote the subgroup of
$V$ consisting of those elements whose restriction to $S$ is the
identity and let $T_S=T\cap V_S$. Applying the work of Farley (see
Proposition~6.1 of \cite{Farley}) and the work of
Sageev~\cite{Sageev}, one sees that $V_{[0,1/2)}$ and
$T_{[0,1/2)}$ admit finite index subgroups $G$ and $H$,
respectively, \st the group pairs $(V,G)$ and $(T,H)$ are
multi-ended (actually, we have $G=V_{[0,1/2)}$, since
$V_{[0,1/2)}\cong V$ and $V$ is infinite and simple). Now
\[
G\subset V_{[0,1/2)}\subset \Gamma\equiv \bigcup_{n=1}^\infty
V_{[0,4\inv+2^{-n})}
\]
As a proper increasing union of a sequence of groups, $\Gamma$ is
infinitely generated. If $\pi_1(M)=V$, then, forming covering
spaces $\check M\to\widehat M\to M$ with $\image{\pi_1(\check
M)\to\pi_1(M)}=G$ and $\image{\pi_1(\widehat
M)\to\pi_1(M)}=\Gamma$, we get $\tilde e(\widehat M)\geq \tilde
e(\check M)\geq e(\check M)\geq 2$. Therefore, by Lemma~\ref{Two
ends and infinitely generated lemma}, we have $\tilde e(\widehat
M)=\infty$. A similar proof applies for $\pi_1(M)=T$.
\end{pf}
\begin{pf*}{Proof that $V$ and $T$ are not K\"ahler (part~(i)
of Theorem~\ref{HNN and Thompson not Kahler theorem from intro})}
If $X$ is a \con \cpt K\"ahler manifold with $\Lambda=\pi_1(X)$
equal to $V$ or $T$, then, by Lemma~\ref{V and T end of pair
lemma} and Theorem~\ref{Mapping of compact to curve theorem}, some
finite covering $X'$ of $X$ admits a surjective \holo mapping with
\con fibers onto a curve $S$ of genus $g\geq 2$. In particular,
some finite index subgroup $\Lambda'$ of $\Lambda$ admits a
surjective homomorphism onto a co\cpt Fuchsian group
$\Theta=\pi_1(S)$. However, $V$ and $T$ are infinite simple
groups, so it follows that $\Lambda=\Lambda'\cong\Theta$. But any
co\cpt Fuchsian group is not simple, so we have arrived at a
contradiction and, therefore, $V$ and $T$ are not K\"ahler.
\end{pf*}

As mentioned above, Farley has independently obtained the result
that $V$ and $T$ are not K\"ahler. In fact, he has recently shown
that $e(V,V_{[0,1/2)})=e(T,T_{[0,1/2)})=\infty$, so one may obtain
the desired proper \holo mapping onto a Riemann surface in the
proof by applying~\cite{NR1} in place of Theorem~\ref{Main Theorem
from Introduction} (in his first proof, Farley applied
Theorem~\ref{Main Theorem from Introduction} and a theorem of
Klein~\cite{Klein} according to which, if $G$ is a finitely
generated group, $K<H<G$, $[G:H]=[H:K]=\infty$, and $\tilde
e(G,K)>1$, then $\til e(G,H)=\infty$).

Section~\ref{Massive sets section} contains a summary of the
required facts from the theory of massive sets \cite{Gri}. In
Section~\ref{Special ends, cup product section}, we consider the
required versions of the cup product lemma. Section~\ref{Filtered
ends and mappings to Riemann surfaces section} contains the proof
of Theorem~\ref{Main Theorem from Introduction}.
Section~\ref{Mappings of Kahler manifolds to curves section}
contains consequences for \cpt K\"ahler manifolds. Finally,
Section~\ref{Amalgamations and mappings to Riemann surfaces
section} contains the details of the proof of the theorem of
Gromov and Schoen (Theorem~\ref{Gromov-Schoen Theorem}) using
Theorem~\ref{Main Theorem from Introduction}, as well as a proof
that a properly ascending HNN extension with finitely generated
base group is not K\"ahler (part~(ii) of Theorem~\ref{HNN and
Thompson not Kahler theorem from intro}). Section~\ref{Principal
fns Evans-Selberg section} (which does not appear in the version
submitted for publication) provides, for the convenience of the
reader, the proof of Sario's existence theorem of principal
functions \cite{RS} and Nakai's construction of the Evans-Selberg
potential \cite{Nakai1}, \cite{Nakai2}, \cite{SaNo}. These facts
were applied in \cite{NR1} (and, therefore, indirectly here).
However, it is difficult to find proofs for a general oriented
Riemannian manifold in a convenient form in the literature.

\begin{ack}  We would like to thank Misha
Kapovich for bringing the work of Delzant and Gromov to our
attention and for helping us understand the geometric group theory
component. We would also like to thank Matt Brin, Ross Geoghegan,
and John Meier for many helpful conversations on filtered ends and
Thompson's groups. We would like to further thank Ross Geoghegan
for providing a copy of his forthcoming book.
\end{ack}

\section{Massive sets}\label{Massive sets section}

By applying the work of Sario, Nakai, and their collaborators (see
\cite{Nakai1}, \cite{Nakai2}, \cite{SaNa}, \cite{SaNo}, \cite{RS},
and also \cite{LT}), one can produce independent \harm \fns by
prescribing limiting values along ends. This is the main fact from
potential theory applied in \cite{NR1}. In order to produce
independent \harm \fns on a manifold with a small number of ends,
it is natural to consider massive sets (in place of ends) as
studied by Grigor'yan \cite{Gri}. Throughout this section, $(M,g)$
will denote a \con Riemannian manifold.

\begin{defn}\label{Massive set definition}
Let $U$ be an open subset of $M$.
\begin{enumerate}
\item[(a)] A bounded nonnegative \cont sub\harm \fn $\alpha $ on $M$
\st
$$
\alpha \equiv 0 \text{ on }M\setminus U \qquad \text{and}\qquad
\sup _M\alpha =\sup _U\alpha >0
$$
is called an {\it admissible sub\harm \fnns} for $U$.
\item[(b)] If there exists an admissible sub\harm \fn for $U$,
then $U$ is called {\it massive}.
\item[(c)] If there exists an admissible sub\harm \fn for $U$ with
finite energy, then $U$ is called $D$-{\it massive}.
\end{enumerate}
\end{defn}

\begin{rmks}
1. If $M$ is a K\"ahler manifold, then a \plsh \fn is sub\harmns .
If $\alpha $ is a \plsh admissible sub\harm \fn for $U$, then we
will simply call $\alpha $ an {\it admissible \plsh \fnns} for $U$
and we will say that $U$ is {\it plurimassive}. If, in addition,
$\alpha $ has finite energy, then we will say that $U$ is {\it
pluri-$D$-massive}.

\noindent 2. To say that a \cont \fn $\alpha $ has {\it finite
energy} is to say that $\alpha \in W^{1,2}_{\text{loc}}(M,g)$ and
$$
\int _M|d\alpha |^2_g \, dV_g <\infty .
$$

\noindent 3. If $M$ contains a proper massive subset $U$, then $M$
is hyperbolic; i.e.~$M$ admits a positive   Green's \fn $G$.
Moreover, if $\alpha $ is an admissible sub\harm \fn for~$U$, then
$\lim _{j\to\infty }G(x_j,\cdot )=0$ for any sequence $\seq xj$ in
$M$ \st $\alpha (x_j)\to \sup _U \alpha $ as $j\to \infty $.

\noindent 4. An end $E \subset M$ with $E\neq M$ is a hyperbolic
end if and only if $E$ is a massive set. In fact, a hyperbolic end
$E$ is $D$-massive. For if we fix a $\cinf $ \rel \cpt domain
$\Omega $ in $M$ \st $\partial E\subset \Omega $ and $M\setminus
\Omega $ has no \cpt \comps and we let $u\colon
M\setminus\overline\Omega\to[0,1)$ be the \harm measure of the
ideal boundary of $M$ \wrt $M\setminus \overline \Omega $, then,
for $0<\epsilon <1$, the \fn
$$
\alpha \equiv \left\{
\begin{aligned}
\max (u-\epsilon ,0)
& \quad \text{on } E\setminus \Omega \\
0 &\quad \text {on } (M\setminus E)\cup \Omega
\end{aligned}
\right.
$$
is a finite energy admissible sub\harm \fn for $E$. Observe also
that, for a sequence $\seq xj$ in $E$, we have, as $j\to \infty $,
$$
G(x_j,\cdot ) \to 0 \iff u(x_j)\to 1 \iff \alpha (x_j)\to \sup
_M\alpha .
$$
\end{rmks}

\begin{prop}[See \cite{Gri}]\label{Massive sets give harmonic
functions prop} Suppose $U$ is a proper massive subset of $M$ and
$\alpha :M\to [0,1)$ is an admissible sub\harm \fn for $U$ with
$\sup _M\alpha =1$. Then there exists a \harm \fn $\rho :M\to
(0,1]$ with the following properties:
\begin{enumerate}
\item[(i)] $\alpha \leq \rho \leq 1$ on $M$;
\item[(ii)] If $M\setminus \overline U$ is massive and $\beta
:M\to [0,1)$ is an admissible sub\harm \fn for $M\setminus
\overline U$, then $0<\rho \leq 1-\beta $ on $M$.
\item[(iii)] If $\alpha $ has finite energy (hence $U$ is
$D$-massive), then $\rho $ has finite energy. In fact,
$$
\int _M|d\rho |_g^2 \, dV_g\leq \int _M|d\alpha |_g^2 \, dV_g.
$$
\end{enumerate}
\end{prop}
\begin{rmk}
If $M\setminus \overline U$ is massive, then $\rho $ is
nonconstant and the maximum principle gives $0<\rho <1$ on $M$.
\end{rmk}

\begin{pf}
We may choose a sequence of $\cinf $ domains ${\seq \Omega
m}_{m=1}^\infty $ in $M$ \st
$$
\Omega _m\Subset \Omega _{m+1} \text{ for } m=1,2,3,\dots
\qquad\text{and}\qquad \bigcup _{m=1}^\infty \Omega _m=M.
$$
For each $m$, let $\rho _m:\overline \Omega _m \to [0,1)$ be the
\cont \fn satisfying
$$
\lap \rho _m=0 \text{ in } \Omega _m \qquad\text{and}\qquad \rho
_m=\alpha \text{ on } \partial \Omega _m.
$$
Since $\alpha $ is sub\harmns , we have $\alpha \leq \rho _m$ on
$\overline\Omega _m$ and, in particular, on $\partial \Omega
_m\subset \Omega _{m+1}$, we have $\rho _m=\alpha \leq \rho
_{m+1}$. Thus $\rho _m\leq \rho _{m+_1}$ on $\Omega _m$ and hence
$\rho _m\nearrow \rho $ for some \harm \fn $\rho :M\to (0,1]$ with
$\alpha \leq \rho \leq 1$ on $M$.

If $\beta :M\to [0,1)$ is an admissible sub\harm \fn for
$M\setminus \overline U$, then the super\harm \fn $1-\beta $
satisfies, for each $m=1,2,3,\dots $,
$$
\rho _m=\alpha =0\leq 1-\beta \text{ on } (\partial \Omega
_m)\setminus U \qquad\text{and}\qquad \rho _m=\alpha \leq 1=
1-\beta \text{ on } (\partial \Omega _m)\cap U .
$$
It follows that $\rho _m\leq 1-\beta $ on $\Omega _m$ and hence
$0<\rho \leq 1-\beta $ on $M$.

Finally, suppose $\alpha $ has finite energy. Then, since \harm
\fns minimize energy, we have, for each $m$,
$$
\int _{\Omega _m}|d\rho _m|_g^2 \, dV_g\leq \int _{\Omega
_m}|d\alpha |_g^2 \, dV_g \leq \int _M|d\alpha |_g^2 \, dV_g.
$$
Therefore, since $d\rho _m\to d\rho $ uniformly on \cpt sets, we
get
$$
\int _M|d\rho |_g^2 \, dV_g\leq \int _M|d\alpha |_g^2 \, dV_g.
$$
\end{pf}

\section{Special ends and the cup product lemma}\label{Special ends, cup product section}

It will be convenient to have the terminology contained in the
next two definitions.

\begin{defn}\label{Bounded geometry and weakly 1-complete
definition} Let $S$ be a subset of a \cpx manifold $X$ of
dimension $n$.
\begin{enumerate}
\item[(a)] For $g$ a Hermitian metric on $X$ and $k$ a nonnegative
integer, we will say that $(X,g)$ {\it has bounded geometry of
order $k$ along $S$} if, for some constant $C>0$ and for every
point $p\in S$, there is a biholomorphism $\Psi $ of the unit ball
$B=B(0;1)\subset \C ^n$ onto a \nbd of $p$ in $X$ \st $\Psi (0)=p$
and, on $B$, $$ C\inv g_{\C ^n} \leq \Psi ^*g\leq Cg_{\C ^n}
\quad\text{and}\quad | D^m\Psi ^*g| \leq C \text{ for }
m=0,1,2,\dots , k.
$$
\item[(b)] We will say that $X$ is {\it weakly $1$-complete along}
$S$ if there exists a \cont \plsh \fn $\vphi $ on $X$ \st
$$
\setof {x\in S}{\vphi (x)<a}\Subset X \quad \forall \, a\in \R .
$$
\end{enumerate}
\end{defn}
\begin{rmk}
Both (a) and (b) hold if $S\Subset X$.
\end{rmk}

\begin{defn}\label{Special end definition}
We will call an end $E\subset X$ in a \con non\cpt complete
K\"ahler manifold $(X,g)$ {\it special} if $E$ is of at least one
of the following types:
\begin{enumerate}
\item[(BG)] $(X,g)$ has bounded geometry of order $2$ along $E$;
\item[(W)] $X$ is weakly $1$-complete along $E$; \item[(RH)] $E$
is a hyperbolic end and the Green's \fn vanishes at infinity along
$E$; or \item[(SP)] $E$ is a parabolic end, the Ricci curvature of
$g$ is bounded below on $E$, and there exist positive constants
$R$ and $\delta $ \st
$$
\Vol \big( B(p;R)\big) >\delta \quad \forall \, p\in E.
$$
\end{enumerate}
An ends decomposition for $X$ in which each of the ends is special
will be called a {\it special ends decomposition}.
\end{defn}
\begin{rmks}
1. (BG) stands for ``bounded geometry,'' (W) for ``weakly
$1$-complete,'' (RH) for ``regular hyperbolic,'' and (SP) for
``special parabolic.''

\noindent 2. A parabolic end of type (BG) is also of type (SP).

\noindent 3. If $E$ and $E'$ are ends with $E'\subset E$ and $E$
is special, then $E'$ is special

\noindent 4. For our purposes, it is generally enough to replace
the condition (BG) by the condition that $E$ is a hyperbolic end
along which $(X,g)$ has bounded geometry of order $0$. However, we
then lose the condition described in Remark~3 above and so the
existence of a special ends decomposition is no longer determined
by the set of ends $\displaystyle{\lim _{\leftarrow } \pi _0}
(X\setminus K)$.
\end{rmks}

\begin{lem}\label{Special ends to produce a weakly 1-complete end lemma}
Suppose $(X,g)$ is a \con complete K\"ahler manifold with an ends
decomposition
$$
X\setminus K=E_1\cup E_2\cup \cdots \cup E_m;
$$
where $m\geq 2$, $E_1$ is a special end of the type (W), (RH), or
(SP) in Definition~\ref{Special end definition}, and, for
$j=2,3,\dots , m$, $E_j$ is a hyperbolic or special end. Then $X$
is weakly $1$-complete along $E_1$ (i.e.~$E_1$ is of type~(W)).
\end{lem}
\begin{pf}
We may arrange $E_2,\dots ,E_m$ so that, for some $k$ with $1\leq
k\leq m$, $E_2,\dots, E_k$ are ends of type~(W) and $E_{k+1},\dots
, E_m$ are ends which are {\it not} of type (W). We set
$E'=E_2\cup \cdots \cup E_k$ and $E''=E_{k+1}\cup \cdots \cup
E_m$. We have a \cont \plsh \fn $\vphi :X\to \R $ \st $\setof{x\in
E'}{\vphi (x)<a}\Subset X$ for each $a\in \R $ and, by fixing a
domain $\Omega $ with $K\subset \Omega \Subset X$ and replacing
$\vphi $ with the \fn
$$
\left\{
\begin{aligned}
\max (\vphi -\max _{\overline \Omega }\vphi -1,0)
&\quad \text{on } E'\setminus\Omega\\
0&\quad \text {on } {\Omega\cup(X\setminus E')}
\end{aligned}
\right.
$$
we may assume that $\vphi \geq 0$ on $X$ and $\vphi \equiv 0$ on
the \con set
$$
\Omega \cup (X\setminus E')=\Omega \cup E_1\cup E''.
$$
By a theorem of \cite{Nk} and \cite{Dem}, the component $Y$ of
$\setof{x\in X}{\vphi (x)<1}$ containing $\Omega \cup (X\setminus
E')$ admits a complete K\"ahler metric $h$ \st $h=g$ on a \nbd of
$\Omega \cup (X\setminus E')$ (one must modify their proofs
slightly since the associated \plsh \fn $-\log (1-\vphi )$
exhausts $Y$ at the \cpt boundary $\partial Y$ but not entirely
along $Y$). Thus $E'\cap Y$ is a union of ends of $Y$ of type (RH)
and hence we may assume without loss of generality that $E_2,\dots
,E_k$ are of type (RH) (as well as type (W)). Thus, for each
$j=2,\dots ,m$, $E_j$ is hyperbolic or $E_j$ is of type~(SP). If
$E_1$ is of type (SP), then Theorem~2.6 of \cite{NR1} (which is
contained implicitly in the work of Sario, Nakai, and their
collaborators \cite{Nakai1}, \cite{Nakai2}, \cite{SaNa},
\cite{SaNo}, \cite{RS} together with the work of Sullivan
\cite{Sul}) provides a \plh \fn $\rho :X\to \R $ \st $\lim _{x\to
\infty} \rho \restrict{\overline E_1}(x)=\infty $. In particular,
$E_1$ is of type (W) in this case. Thus we may assume that $E_1$
is of type (RH). Moreover, if, for some $j$ with $k+1\leq j\leq
m$, $E_j$ were of type (SP), then $X$ would be weakly $1$-complete
along $E_j$ and hence $E_j$ would have been included in $E'$.
Therefore, $E_{k+1},\dots ,E_m$ must be hyperbolic ends.

Thus we have hyperbolic, and therefore $D$-massive, ends
$E_1,\dots ,E_m$ with $m\geq 2$ and $E_1$ is of type (RH). By
Remark~4 following Definition~\ref{Massive set definition}, there
is a finite energy admissible sub\harm \fn $\alpha :X\to [0,1)$
for $E_1$ \st
$$
\lim _{x\to \infty }\alpha \restrict{\overline E_1}(x)=1
$$
and, applying Proposition~\ref{Massive sets give harmonic
functions prop} ($X\setminus \overline E_1\supset E'\cup E''$ is
massive), we get a finite energy \harm \fn $\rho :X\to (0,1)$ \st
$\alpha \leq \rho <1$ on $X$. The Gaffney theorem \cite{Ga}
implies that $\rho $ is \plh and we have
$$
1>\rho (x)\geq \alpha (x) \to 1  \text{ as }  x\to \infty \text{
in } \overline E_1.
$$
Therefore $X$ is weakly $1$-complete along $E_1$ (with \plsh \fn
$-\log (1-\rho )$ exhausting $\overline E_1$) in this case as
well.

\end{pf}

As described in the introduction, the main goal of this paper is
to obtain a filtered ends version of the following:
\begin{thm}[\cite{Gro1}, \cite{L}, \cite{Gro2}, and Theorem~3.4 of
\cite{NR1}]\label{Old 3 ends theorem} If $(X,g)$ is a \con
complete K\"ahler manifold which admits a special ends
decomposition and $e(X)\geq 3$, then $X$ admits a proper \holo
mapping onto a Riemann surface.
\end{thm}
The following lemma is well known (see, for example, the proof of
Theorem~4.6 of \cite{NR1}):
\begin{lem}\label{Open set to Riemann surface gives global lemma}
Let $(X,g)$ be a \con complete K\"ahler manifold which is \cpt or
which admits a special ends decomposition. If some nonempty open
subset of $X$ admits a surjective proper \holo mapping onto a Riemann
surface, then $X$ admits a surjective proper \holo mapping onto a
Riemann surface.
\end{lem}

\begin{lem}[See Lemma~2.1 of \cite{NR3}]\label{Cup product for compact fiber lemma}
Let $\rho _1$ and $\rho _2$ be two
real-valued \plh \fns on a \con \cpx manifold~$X$. If $\rho _1$
has a nonempty \cpt fiber, then $\partial \rho _1\wedge \partial
\rho _2\equiv 0$ on $X$. Furthermore, if the differentials $d\rho
_1$ and $d\rho _2$ are (globally) linearly independent for some
such pair of \fnsns , then some nonempty open subset of $X$ admits
a proper \holo mapping onto a Riemann surface.
\end{lem}
\begin{rmk}
Two real-valued \plh \fns $\rho _1$ and $\rho _2$ on a \con \cpx
manifold have linearly {\it dependent} differentials $d\rho _1$
and $d\rho _2$ (i.e.~$1$, $\rho _1$, $\rho _2$ are linearly
dependent \fnsns ) if and only if $d\rho _1\wedge d\rho _2\equiv
0$.
\end{rmk}

\begin{lem}\label{Cup product for sublevel of plh function lemma}
Let $(X,g)$ be a \con complete K\"ahler manifold and let $\rho _1$
and $\rho _2$ be two real-valued \plh \fns on a domain $Y\subset
X$. Assume that, for some constant~$a$ with $\inf \rho _1
<a<\sup\rho _1$, some \comp $\Omega $ of $\setof{x\in Y}{a<\rho
_1(x)}$ has the following properties:
\begin{enumerate}
\item[(i)] $\overline \Omega \subset Y$;
\item[(ii)] $|d\rho _1|_g$ is bounded on $\Omega $; and
\item[(iii)] $\int _\Omega |d\rho _j|^2_g\, dV_g <\infty $ for
$j=1,2$.
\end{enumerate}
Then $\partial \rho _1\wedge \partial \rho _2\equiv 0$ on $Y$.

Furthermore, if $d\rho _1$ and $d\rho _2$ are linearly independent
and $(X,g)$ has bounded geometry along $\Omega $, then $\Omega $
admits a proper \holo mapping onto a Riemann surface.
\end{lem}
\begin{rmk}
In this paper, we will only need the fact that some nonempty open
subset of $\Omega $ admits a proper \holo mapping onto a Riemann
surface.
\end{rmk}
\begin{pf*}{Proof of Lemma \ref{Cup product for sublevel of plh function lemma}}
We may assume without loss of generality that $\rho _1$ and $\rho
_2$ are nonconstant. We denote the {\it Levi form} of a $C^2$ \fn
$\vphi $ by
$$
\lev \vphi =\sum _{i,j=1}^n \frac {\partial ^2\vphi }{\partial
z_i\partial \bar z_j} dz_id\bar z_j.
$$
The Hermitian tensor
$$
h=g+\lev {-\log (\rho _1-a)} =g+(\rho _1-a)^{-2} \partial \rho _1
\overline {\partial \rho _1}
$$
is a complete K\"ahler metric on $\Omega $ with $h\geq g$ (see
\cite{Nk} or \cite{Dem}). Moreover, on $\Omega $, we have
$$
|\partial \rho _1|_h^2=\bigl[ 1+(\rho _1-a)^{-2} |\partial \rho
_1|^2_g \bigr] \inv |\partial \rho _1|^2_g \text{, } dV_h=\bigl[
1+(\rho _1-a)^{-2} |\partial \rho _1|^2_g \bigr]dV_g, \text{ and }
|\partial \rho _2|^2_h\leq |\partial \rho _2|^2_g.
$$
For $g=h$ at any point $p\in \Omega $ where $(\partial \rho
_1)_p=0$ while, at any point $p\in \Omega $ where $(\partial \rho
_1)_p\neq 0$, one may get the above by writing $g$ and $h$ in
terms of a $g$-orthonormal basis $e_1, \dots , e_n$ for
$T^{1,0}_p\Omega $ with dual basis
$$
e_1^*=|\partial \rho _1|\inv _g(\partial \rho _1)_p, e_2^*, \dots
, e_n^*.
$$
In particular, we have
$$
|\partial \rho _1|^2_hdV_h=|\partial \rho _1|^2_gdV_g
$$
and hence $\rho _1\restrict\Omega $ has finite energy \wrt $h$ as
well as $g$.

Following \cite{Gro2}, we consider the closed \holo $1$-forms
$$
\alpha _j=\partial \rho _j=\beta _j+i\gamma _j \quad \text{for }
j=1,2,
$$
where
$$
\beta _j=\frac 12(\partial \rho _j+\dbar\rho _j)=\frac 12d\rho _j
\qquad\text{and}\qquad \gamma _j=\frac 1{2i}(\partial \rho
_j-\dbar\rho _j)=\frac 12d^c\rho _j
$$
are closed real $1$-forms (since $\rho _j$ is \plhns ) for
$j=1,2$. Thus we get a closed \holo $2$-form
$$
\alpha _1\wedge \alpha _2=\eta +i\theta
$$
on $Y$ where
$$
\eta =\beta _1\wedge \beta _2-\gamma _1\wedge \gamma _2
\qquad\text{and}\qquad \theta =\beta _1\wedge \gamma _2+\gamma
_1\wedge \beta _2
$$
are closed real $2$-forms. Moreover, $\alpha _1\wedge \alpha _2$
is in $L^2$ \wrt $h$ on $\Omega $ because
$$
|\alpha _1\wedge \alpha _2|^2_h dV_h \leq |\alpha _1|^2_h |\alpha
_2|^2_h dV_h=|\partial\rho _1|^2_h |\partial\rho _2|^2_h dV_h
=|\partial\rho _2|^2_h |\partial\rho _1|^2_gdV_g \leq
|\partial\rho _2|^2_g |\partial\rho _1|^2_gdV_g
$$
while $|\partial\rho _1|_g$ is bounded on $\Omega $ and $\rho
_2\restrict\Omega $ is of finite energy \wrt $g$. Furthermore, the
closed form $\alpha _1\wedge \alpha _2$ is \harm \wrt the complete
K\"ahler metric $h$ (and $g$) because $\dbar ^*_h(\alpha _1\wedge
\alpha _2)=0$ (since $\alpha _1\wedge \alpha _2$ is of
type~$(2,0)$) and $\dbar (\alpha _1\wedge \alpha _2)=0$. Therefore
$$
d(\alpha
_1\wedge \alpha _2)=d\eta =d\theta =0
\text{ and }
d^*_h(\alpha _1\wedge \alpha _2)=d_h^*\eta =d_h^*\theta =0
$$
($L^2$ \harm forms are closed and coclosed by the Gaffney theorem
\cite{Ga}).

For each $R>0$, let $\vphi _R:\Omega \to (0,R]$ and $\psi
_R:\Omega \to [-R,R]$ be the bounded locally Lipschitz \fns given
by, for each $x\in \Omega $,
$$
\vphi _R(x)=\frac 12\min (\rho _1(x)-a,R)
$$
and
$$
\psi _R(x)= \left\{
\begin{aligned}
\frac 12R
& \quad \text{if } \rho _2(x)\leq -R\\
-\frac 12\rho _2(x)
& \quad \text{if } -R<\rho _2(x)<R\\
-\frac 12R & \quad \text{if } \rho _2(x)\geq R
\end{aligned}
\right.
$$
Since $\rho _1$ and $\rho _2$ are nonconstant \plh \fnsns , the
sets $\setof{x\in Y}{\rho _1(x)=a+R}$ and $\setof{x\in Y}{|\rho _2
(x)|=R}$ are sets of measure~$0$ and we have
$$
d\vphi _R= \left\{
\begin{aligned}
\frac 12d\rho _1
& \quad \text{on } \setof{x\in \Omega }{\rho _1(x)<a+R}\\
0 & \quad \text{on } \setof{x\in \Omega }{\rho _1(x)>a+R}
\end{aligned}
\right.
$$
and
$$
d\psi _R= \left\{
\begin{aligned}
-\frac 12d\rho _2
& \quad \text{on } \setof{x\in \Omega }{|\rho _2(x)|<R}\\
0 & \quad \text{on }\setof{x\in \Omega }{|\rho _2(x)|>R}
\end{aligned}
\right.
$$
Thus
$$
\beta _1\wedge \gamma _2=\frac 12d\rho _1\wedge \gamma
_2=d\bigl[\frac 12\rho _1\gamma _2 \bigr] =d[\vphi _R\gamma _2]
\qquad\text{ when } \rho _1<a+R
$$
and
$$
\gamma _1\wedge \beta _2=\gamma _1\wedge \frac 12d\rho
_2=d\bigl[-\frac 12\rho _2\gamma _1 \bigr] =d[\psi _R\gamma _1]
\qquad\text{ when } |\rho _2|<R.
$$

Moreover, $\vphi _R\gamma _2$ is in $L^2$ \wrt $h$ on $\Omega $
because
\begin{align*}
|\vphi _R\gamma _2|^2_h dV_h&\leq \vphi _R^2|\partial\rho _2|^2_h
dV_h \leq\vphi _R^2|\partial\rho _2|^2_g dV_h= \vphi
_R^2|\partial\rho _2|^2_g\bigl( 1+(\rho _1-a)^{-2}|\partial \rho
_1|^2_g  \bigr) dV_g
\\
&=\bigl(\vphi _R^2+\vphi _R^2(\rho _1-a)^{-2}|\partial \rho
_1|^2_g \bigr) |\partial\rho _2|^2_g dV_g \leq \bigl(R^2+|\partial
\rho _1|^2_g \bigr) |\partial\rho _2|^2_g dV_g
\end{align*}
while $|\partial \rho _1|_g$ is bounded on $\Omega $ and $\rho
_2\restrict\Omega $ has finite $g$-energy. The form $\psi _R\gamma
_1$ is also in $L^2$ \wrt $h$ because $|\psi _R|\leq R$ and $\rho
_1\restrict\Omega$ has finite energy \wrt $h$ (as well as $g$).
Thus the form $\lambda _R=\vphi _R\gamma _2+\psi _R\gamma _1$ is
in $L^2$ \wrt $h$ on $\Omega $ and, as $R\to \infty $, $\lambda
_R$ converges pointwise to $\lambda =\frac 12(\rho _1-a)\gamma
_2-\frac 12\rho _2\gamma _1$.

Moreover, $d\lambda _R \to d\lambda =\theta $ in $L^2$ \wrt $h$ as
$R\to \infty $. For $d\lambda _R=\theta $ when we have both $\rho
_1<a+R$ and $|\rho _2|<R$, so $d\lambda _R\to \theta $ pointwise
as $R\to \infty $. We also have (almost everywhere)
$$
|d\lambda _R|_h\leq |d\vphi _R\wedge \gamma _2|_h+|d\psi
_R\wedge\gamma _1|_h \leq |d\vphi _R|_h|\gamma _2|_h+|d\psi
_R|_h|\gamma _1|_h \leq |\beta _1|_h|\gamma _2|_h+|\beta
_2|_h|\gamma _1|_h
$$
and the last expression is in $L^2$ \wrt $h$ (by an argument
similar to that showing that $\alpha _1\wedge \alpha _2$ is in
$L^2$ \wrt $h$). Thus the Lebesgue dominated convergence theorem
gives the claim.

By the Gaffney theorem \cite{Ga}, we get
$$
0=\int _{\Omega } \langle 0,\lambda _R \rangle _h\, dV_h =\int
_{\Omega } \langle d^*_h\theta ,\lambda _R \rangle _h\, dV_h =\int
_{\Omega } \langle \theta ,d\lambda _R \rangle _h\, dV_h \to \int
_{\Omega } |\theta |^2_h\, dV_h
$$
as $R\to \infty $. Thus $\theta \equiv 0$ on $\Omega $ and hence,
since $\alpha _1\wedge \alpha _2=\eta +i\theta $ is a \holo
$2-$form, we get $\alpha _1\wedge \alpha _2\equiv 0$ on $\Omega $
and, therefore, on $Y$.

Assuming now that $d\rho _1$ and $d\rho _2$ are linearly
independent and that $(X,g)$ has bounded geometry along $\Omega $,
the arguments in \cite{Gro2}, \cite{ABR} (see also Chapter~4 of
\cite{ABCKT}), together with some easy observations, give the
required proper \holo mapping of $\Omega $ onto a Riemann surface.
For the convenience of the reader, we include a sketch of the
arguments.

We may assume $n=\dim X>1$.
Because we have $\alpha _1\wedge \alpha _2=\partial \rho _1\wedge
\partial \rho _2\equiv 0$, the \mero \fn $\frac {\alpha _1}{\alpha
_2}:\Omega \to {\mathbb P} ^1$ is actually a \holo map
(i.e.~$\alpha _1/\alpha _2$ has no points of indeterminacy). If
this \fn is equal to a constant $\zeta $, then $\rho _1-\bar\zeta
\rho _2$ is a nonconstant \holo \fn on $\Omega $. In any case, we
get a nonconstant \holo map $f:\Omega \to {\mathbb P}^1$ \st $\rho
_1$ and $\rho _2$ are constant on each level of $f$. In fact, $f$
is locally constant on the (\cpxns ) \anal set
$$
A=\setof{x\in \Omega }{(\alpha _1)_x=(\alpha _2)_x=0}
$$
and the levels of $f\restrict{\Omega \setminus A}$ are precisely
the (smooth) leaves of the \holo foliation determined by $\alpha
_1$ and $\alpha _2$ in $\Omega\setminus A$ (see, for example,
\cite{NR2}, pp.~387--388). If $L$ is a level of $f$, then $\rho
_1$ is equal to a constant $t$ on $L$ with $a<t$ and $\bar
L\subset \overline \Omega \subset Y$. It follows that $\bar
L\subset \Omega $ (since $\Omega $ is a \comp of $\setof{x\in
Y}{\rho _1(x)>a}$) and hence that $\bar L=\bar L\cap \Omega =L$.
Thus $L$ is closed as an \anal subset of $X$.

The coarea formula for the map $\Phi =(\rho _1,\rho _2):\Omega \to
\R ^2$ gives us
$$
\int _{\R ^2}\Vol _{g\restrict{\Phi\inv (t_1,t_2)}}\bigl( \Phi\inv
(t_1,t_2)\bigr) \, dt_1\wedge dt_2= \int _{\Omega }|d\rho _1\wedge
d\rho _2|_g\, dV_g <\infty .
$$
Hence $\Vol \bigl( \Phi\inv (t)\bigr) <\infty $ for almost every
point $t\in \R ^2$. Thus we may fix a regular value $t_0$ in the
interior of $\Phi (\Omega )$, with $t_0$ in the complement of the
countable set $\Phi (A)$, \st $\Vol \bigl( \Phi\inv (t_0)\bigr)
<\infty $. Since $\Phi $ is constant on each leaf of the foliation
in $\Omega \setminus A$ and $\Phi\inv (t_0)\subset \Omega\setminus
A$ is a $\cinf $ submanifold of $\Omega $ (not just $\Omega
\setminus A$) with $\dim _{\R }\Phi\inv (t_0)=2n-2$, we see that a
\comp $L_0$ of $\Phi\inv (t_0)$ is a leaf of the foliation in
$\Omega \setminus A$, $L_0$ is closed in $X$, and $L_0$ is a level
of $f$.

Since $(X,g)$ has bounded geometry along $\Omega $, Lelong's
monotonicity formula (see 15.3 in \cite{Chirka}) shows that there
is a constant $C>0$ \st each point $p\in \Omega $ has a \nbd $U_p$
in $X$ \st $\diam _XU_p<1$ and $\Vol (D\cap U_p)\geq C$ for every
\cpx \anal set $D$ of pure dimension $n-1$ in $X$ with $p\in D$.
Therefore, since $L_0$ has finite volume, $L_0$ must be \cptns .
Thus $f:\Omega \to {\mathbb P}^1$ has a \cpt level $L_0\subset
\Omega $.

It follows that the set $V=\setof{x\in \Omega }{x\text{ lies in a
\cpt level of }f}$ is a nonempty open subset of $\Omega $. To show
that $V$ is also closed relative to $\Omega $, let $V_0$ be a
\comp of $V$, let $\seq xj$ be a sequence in $V_0$ converging to a
point $p\in \overline V_0\cap \Omega $, and, for each $j$, let
$L_j\subset V_0$ be the \cpt level of $f$ through $x_j$. Stein
factoring $f\restrict{V_0}$, we get a proper \holo mapping $\Psi
:V_0 \to W$ with \con fibers of $V_0$ onto a Riemann surface $W$.
We may choose each $x_j$ to lie over a regular value of $f$ and of
$\Psi $. Applying Stokes' theorem as in \cite{Stoll}, we see that
$\Vol (L_j)$ is constant in $j$ and so the above volume estimate
implies that, for some $R\gg 0$, we have $L_j\subset B(p;R)$ for
$j=1,2,3,\dots $. On the other hand, by \cite{Stein} (see
\cite{TW} and Theorem~4.23 in \cite{ABCKT}), a subsequence of
$\seq Lj$ converges to the level $L$ of $f$ through $p$. So we
must have $L\subset \overline {B(p;R)}\cap \Omega $ and hence,
since $L$ is a closed \anal subset of $X$, $L$ must be \cptns .
Thus $p\in \overline V_0\cap V$ and, therefore, $p\in V_0$. It
follows that $V=V_0=\Omega $. Thus every level of $f$ is \cpt and
we get our proper \holo map $\Psi :\Omega \to W$.
\end{pf*}

\section{Filtered ends and mappings to Riemann surfaces}\label{Filtered ends and mappings to Riemann
surfaces section}

Theorem~\ref{Main Theorem from Introduction} is an immediate
consequence of the following theorem which will be proved in this
section:
\begin{thm}\label{General filtered ends theorem}
If $(X,g)$ is a \con complete K\"ahler manifold which admits a
special ends decomposition and $\tilde e(X)\geq 3$, then $X$
admits a proper \holo mapping onto a Riemann surface.
\end{thm}

We first consider two lemmas which will allow us to replace
special ends of type (W) with special ends of type (RH).

\begin{lem}\label{filtered ends of subdomains lemma}
Let $M$ be a \con non\cpt $\cinf$ manifold and let $k\in\N$.
\begin{enumerate}
\item[(a)] Given an end $E$ in $M$ with
$\left[\pi_1(M):\image{\pi_1(E)\to\pi_1(M)}\right]\geq k$, there
exists a \cpt set $D\subset M$ \stns, if $\Omega$ is a domain
containing $D$, then $\Omega\cap E$ is an end of $\Omega$ and, for
any end $F$ of~$\Omega$ contained in $E$, we have
$\left[\pi_1(\Omega):\image{\pi_1(F)\to\pi_1(\Omega)}\right]\geq
k$.

\item[(b)] If $\til e(M)\geq k$, then there exists a \cpt set
$D\subset M$ \stns, for every domain $\Omega$ containing $D$, we
have $\til e(\Omega)\geq k$.

\end{enumerate}
\end{lem}
\begin{pf}
For the proof of (a), we fix a point $x_0\in E$ and a finite set
$A$ of $k$ loops in $M$ based at $x_0$ \stns, if $\alpha,\beta\in
A$ and the product loop $\beta\inv*\alpha$ is homotopic to a loop
in $E$, then $\alpha=\beta$. We may now choose a \cpt set
$D\subset M$ \st $\overset\circ D$ contains $\partial E$, $D$
contains the image $\alpha([0,1])$ of each loop $\alpha\in A$, and
$D\cap E$ is path \con (for example, choosing a $\cinf$ \rel \cpt
domain $U$ in $M$ \st \(\partial E\cup\bigcup_{\alpha\in
A}\alpha([0,1])\subset U\), we may let $D$ be the union of
$\overline U$ with the images of finitely many paths in $E$
joining the boundary components of $U$ contained in $E$).

Suppose $\Omega$ is a domain containing $D$ and $\Omega\sm
K=F_1\cup\cdots\cup F_m$ is an ends decomposition for $\Omega$
with $F=F_1\subset E$. The intersection $\Omega\cap E$ is \conns.
For a path in $\Omega$ with endpoints in $\Omega\cap E$ which
leaves $E$ must meet $D\supset\overset\circ D\supset\partial E$.
Hence the segments between the endpoints and the first and last
points in $D$ together with a path in $D\cap E$ joining these
first and last points yields a path in $\Omega\cap E$ between the
endpoints. Thus $\Omega\cap E$ is an end of $\Omega$. In
particular, we may fix a point $y_0\in F$ and a path $\lambda$ in
$\Omega\cap E$ from $x_0$ to $y_0$, and we may let $B$ be the
finite set of loops in $\Omega$ given by
$B=\setof{\lambda\inv*\alpha*\lambda}{\alpha\in A}$. If
$\alpha,\beta\in A$ and
$[\lambda\inv*\beta*\lambda]\inv*[\lambda\inv*\alpha*\lambda]$ is
homotopic in $\Omega$ to a loop in $F$, then $\beta\inv*\alpha$ is
homotopic in $X$ to a loop in $E$ and hence $\alpha=\beta$. Thus
\[
\left[\pi_1(\Omega):\image{\pi_1(F)\to\pi_1(\Omega)}\right]\geq
\#B=\#A=k.
\]

For the proof of~(b), we fix positive integers
$k_1,\dots,k_m\in\N$ and an ends decomposition $M\sm
K=E_1\cup\cdots\cup E_m$ \st
$\left[\pi_1(M):\image{\pi_1(E_j)\to\pi_1(M)}\right]\geq k_j$ for
each $j=1,\dots,m$ and $\sum k_j\geq k$. By~(a), we may choose a
\cpt set $D\subset M$ \stns, if $\Omega$ is any domain containing
$D$ and $1\leq j\leq m$, then $F_j=E_j\cap\Omega$ is an end of
$\Omega$ and
$\left[\pi_1(\Omega):\image{\pi_1(F_j)\to\pi_1(\Omega)}\right]\geq
k_j$. The claim now follows.
\end{pf}

\begin{lem}\label{Type W to RH lemma}
Let $(X,g)$ be a \con complete K\"ahler manifold, let $E$ be a
special end of type~(W) in $X$, let $k,l\in\N$ with $\til e(X)\geq
k$ and $\left[\pi_1(X):\image{\pi_1(E)\to\pi_1(X)}\right]\geq l$,
and let $D$ be a \cpt subset of $X$. Then there exists a domain
$X'$ in $X$, a complete K\"ahler metric $g'$ on $X'$, and an ends
decomposition $X'\sm K=E_0\cup E_1\cup E_2\cup\cdots\cup E_m$ \st
\begin{enumerate}
\item[(i)] $(X\sm E)\cup D\subset E_0$;

\item[(ii)] On $E_0$, $g'=g$;

\item[(iii)]  For each $j=1,\dots,m$, $E_j$ is a special end of
type~(RH) and (W) satisfying $\left[\pi_1(X):\text{\rm
im}\,[\pi_1(E_j)\to\pi_1(X)]\right]\geq l$;

\item[(iv)] $\til e(X')\geq k$.

\end{enumerate}
\end{lem}
\begin{pf}
By Lemma~\ref{filtered ends of subdomains lemma}, we may assume
without loss of generality that $D$ is \conns; $\partial E\subset
D$; and, if $\Omega$ is any domain in $X$ containing $D$, then
$\til e(\Omega)\geq k$, $\Omega\cap E$ is an end of $\Omega$, and,
for any end $F$ of~$\Omega$ contained in $E$, we have
$\left[\pi_1(\Omega):\image{\pi_1(F)\to\pi_1(\Omega)}\right]\geq
l$.

By hypothesis, there exists a \cont \plsh \fn $\psi $ on $X$ which
exhausts $\overline E$. For constants $r_1$ and $r_2$ with
$\max_D\psi <r_1<r_2$, the \comp $X'$ of $\setof{x\in E}{\psi
(x)<r_2}\cup (X\setminus E)$ containing the \con set $(X\setminus
E)\cup D$ admits a complete K\"ahler metric $g'$ \st $g'=g$ on
$[\setof{x\in E}{\psi (x)<r_1}\cup (X\setminus E)]\cap X'$ (see
\cite{Nk}, \cite{Dem}). Since $\psi\to r_2$ at $\partial X'$, the
\comp $E_0$ of $\setof{x\in E}{\psi (x)<r_1}\cup (X\setminus E)$
containing $(X\setminus E)\cup D$ is contained in $X'$ and the set
$K\equiv X'\sm [E_0\cup\setof{x\in E\cap X'}{\psi (x)>r_1}]$ is
\cptns. Furthermore, by the maximum principle, the \comps
$E_1,\dots,E_m$ of the nonempty set $\setof{x\in E\cap X'}{\psi
(x)>r_1}$ are not \rel \cpt in $X'$. Thus the domain $X'$ and the
ends decomposition $X'\sm K=E_0\cup\cdots\cup E_m$ have the
required properties.
\end{pf}

Several cases of Theorem~\ref{General filtered ends theorem} are
contained in the following:
\begin{lem}\label{Two disjoint massive sets gives map lemma}
Let $(X,g)$ be a \con complete K\"ahler manifold which contains a
special end~$E$ and suppose that $X\setminus \overline E$ contains
two disjoint massive subsets $U_1$ and $U_2$ of $X$ \st $U_1$ is
$D$-massive or $U_1$ has an associated $\cinf $ admissible \plsh
\fn (i.e.~$U_1$ is $\cinf $ pluri-massive). Then some nonempty
open subset of $E$ admits a proper \holo mapping onto a Riemann
surface.
\end{lem}
\begin{rmk}
Lemma~\ref{Two disjoint massive sets gives map lemma} remains true
if we allow the admissible \plsh \fn for $U_1$ to be only \contns.
However, it is then harder to produce a complete K\"ahler metric
on a sublevel and we will only need the $\cinf $ case.
\end{rmk}
\begin{pf*}{Proof of Lemma~\ref{Two disjoint massive sets gives map lemma}}
We first observe that we may assume without loss of generality
that $E$ is a $\cinf $ domain and $E'=X\setminus \overline E$ is
\con (i.e.~$E'$ is an end). For we may choose a $\cinf$ \rel \cpt
domain $\Theta$ \st $\partial E\subset\Theta$ and \st $X\setminus
\Theta$ has no \cpt \compsns , and we may replace $E$ with a \comp
of $E\setminus\overline\Theta $. Observe also that $E'$ is then a
massive, and therefore hyperbolic, end.

Next, we observe that we may assume without loss of generality
that $E$ is a hyperbolic end of type (BG) or type (RH). For, if
$E$ is of type~(W), then we may apply Lemma~\ref{Type W to RH
lemma} and work on a suitable subdomain in place of $X$. If $E$ is
of type (SP) (for example, if $E$ is parabolic of type (BG)),
then, by Lemma~\ref{Special ends to produce a weakly 1-complete
end lemma}, $E$ is also of type (W) and the above applies.

Since $E$ is a hyperbolic end, there exists a finite energy
admissible sub\harm \fn $\alpha _0:X\to (0,1)$ for $E$ \st $\alpha
_0(x_j)\to 1$ as $j\to \infty $ whenever $\seq xj$ is a sequence
in $E$ \st $G(x_j,\cdot )\to 0$ as $j\to \infty $; where $G$ is
the Green's \fn on $X$ (see Remark~4 following Definition
\ref{Massive set definition}). Applying Proposition~\ref{Massive
sets give harmonic functions prop}, we get a finite energy \harmns
, hence \plhns , \fn $\rho _1:X\to (0,1)$ \st $\alpha _0\leq \rho
_1<1$ on $X$ and, for any admissible sub\harm \fn $\beta :X\to
[0,1)$ for $E'$, $0<\rho _1\leq 1-\beta $ on $X$. In particular,
$\rho _1(x_j)\to 1=\sup \rho _1$ whenever $\seq xj$ is a sequence
in $E$ \st $G(x_j,\cdot )\to 0$ as $j\to \infty $. We will produce
a second \plh \fn and apply Lemma~\ref{Cup product for compact
fiber lemma} and Lemma~\ref{Cup product for sublevel of plh
function lemma}. Toward this end, we fix a constant $a$ with $\max
_{\partial E}\rho _1<a<1$ and a \comp $\Omega $ of $\setof{x\in
X}{a<\rho _1(x)}$ contained in $E$.

We have admissible sub\harm \fns $\alpha _1$ and $\alpha _2$ for
$U_1$ and $U_2$, respectively, \st $\sup \alpha _1=\sup \alpha
_2=1$ and \st $\alpha _1$ is of finite energy or $\alpha _1$ is
$\cinf $ \plshns .

Assuming first that $\alpha _1$ has finite energy ($U_1$ is
$D$-massive), we may apply Proposition~\ref{Massive sets give
harmonic functions prop} to get a finite energy \plh \fn $\rho
_2:X\to (0,1)$ \st $\alpha _1\leq \rho _2<1$ on $X$ and, for any
admissible sub\harm \fn $\beta :X\to [0,1)$ for $X\setminus
\overline {U_1}$, $0<\rho _2\leq 1-\beta $ on $X$. The \fns
$1$,~$\rho _1$,~and~$\rho _2$ are then linearly independent on
$X$. To see this, suppose $a_1,a_2,a_3\in \R $ with $a_1\rho
_1+a_2\rho _2+a_3\equiv 0$. Choosing a sequence $\seq xj$ in $U_2$
with $\alpha _2(x_j)\to 1$ as $j\to \infty $, we get
$$
0<\rho _1(x_j),\rho _2(x_j) \leq 1-\alpha _2(x_j)\to 0 \text{ as }
j\to \infty
$$
and it follows that $a_3=0$. Taking a sequence $\seq xj$ in $U_1$
with $\alpha _1(x_j)\to 1$ as $j\to \infty $, we get
$$
1>\rho _2(x_j)\geq \alpha _1(x_j)\to 1 \quad\text{and}\quad 0<\rho
_1(x_j) \leq 1-\alpha _1(x_j)\to 0 \text{ as } j\to \infty
$$
and it follows that $a_2=0$. Thus $a_1\rho _1\equiv 0$ and hence
$a_1=0$. In this case, we also set $Y=X$ and $h=g$.

If $\alpha _1$ is $\cinf $ and \plsh ($U_1$ is $\cinf $
pluri-massive), then we fix a \con \cpt set $H$ \st
$$
\partial E\subset H \qquad\text{and}\qquad \max _H\alpha _i>0\text{
for } i=0,1,2.
$$
Fixing $b$ with $\max _H\alpha _1<b<1$, we see that the \comp $Y$
of $\setof{x\in X}{\alpha _1(x)<b}$ containing $H$ admits the
complete K\"ahler metric
\begin{align*}
h=g+\lev {-\log (b-\alpha _1)} &=g+(b-\alpha _1)\inv \lev {\alpha
_1} +(b-\alpha _1)^{-2}\partial \alpha _1\overline {\partial
\alpha _1} \\
&\geq g+(b-\alpha _1)^{-2}\partial \alpha _1\overline {\partial
\alpha _1}\geq g
\end{align*}
(see \cite{Nk}, \cite{Dem}). We have $\overline \Omega \subset
E\subset \overline E \subset Y$, $\sup _Y\alpha _i\geq \sup
_H\alpha _i
>0$ for $i=0,1,2$, and $h=g$ on $Y\setminus U_1\supset E\cup
(U_2\cap Y) $. Thus, \wrt $(Y,h)$, $\alpha _0\restrict Y$ is a
finite energy admissible sub\harm \fn for $E$, $\alpha _1\restrict
Y$ is a $\cinf $ admissible \plsh \fn for $U_1\cap Y$, and $\alpha
_2\restrict Y$ is an admissible sub\harm \fn for $U_2\cap Y$.
Applying Proposition~\ref{Massive sets give harmonic functions
prop}, we get a finite $h$-energy \harmns , hence \plhns , \fn
$\rho _2:Y\to (0,1)$ \st $\alpha _0\leq \rho _2<1$ on $Y$ and, for
any admissible sub\harm \fn $\beta :Y\to [0,1)$ for $Y\setminus
\overline E$, $0<\rho _2\leq 1-\beta $ on $Y$. To see that the
\fns $1$,~$\rho _1\restrict Y$, and~$\rho _2$ are linearly
independent, suppose $a_1,a_2,a_3\in \R $ with $a_1\rho _1+a_2\rho
_2+a_3\equiv 0$. Choosing a sequence $\seq xj$ in $U_2\cap Y$ with
$\alpha _2(x_j)\to \sup _Y\alpha _2$ as $j\to \infty $, we get
$$
0<\rho _1(x_j),\rho _2(x_j) \leq 1-(\alpha _2(x_j)/\sup _Y\alpha
_2) \to 0 \text{ as } j\to \infty
$$
and it follows that $a_3=0$. Here, we have used the fact that
the \fn given by
$$
\left\{
\begin{aligned}
\alpha _2/\sup _Y\alpha _2
&\quad \text{on } U_2\cap Y \\
0 &\quad \text {on } X\setminus (U_2\cap Y)
\end{aligned}
\right.
$$
is an admissible $g$-sub\harm \fn for $U_2\cap Y\subset
X\setminus\overline E$ in $(X,g)$ (which is the case because $g=h$
on $U_2\cap Y$, $\overline U_2\cap\partial Y=\emptyset$, and
$\alpha _2\equiv 0$ on $X\setminus U_2$). Taking a sequence $\seq
xj$ in $U_1\cap Y$ with $x_j\to x_0\in \partial Y$ as $j\to \infty
$, we get
$$
0<\rho _2(x_j) \leq 1-b\inv\alpha _1(x_j)\to 0
\quad\text{and}\quad \rho _1(x_j)\to \rho _1(x_0)>0 \text{ as }
j\to \infty
$$
and it follows that $a_1=0$. Thus $a_2\rho _2\equiv 0$ and hence
$a_2=0$.

If the end $E$ is of type (RH), then, for any $c$ with $a<c<1=\sup
_\Omega\rho _1$, the fiber $\rho _1\inv (c)\cap \Omega $ of $\rho
_1\restrict\Omega $ is \cptns . Thus, in either of the above
cases, Lemma~\ref{Cup product for compact fiber lemma} implies
that some nonempty open subset of $\Omega \subset E$ admits a
proper \holo mapping onto a Riemann surface.

Suppose the end $E\supset \overline \Omega $ is of type (BG). We
have, in either case, $h=g$ on $E$. Thus
$$
\int _\Omega |d\rho _1|^2_g\, dV_g \leq \int _X |d\rho _1|^2_g\,
dV_g    <\infty \qquad\text{and}\qquad \int _\Omega |d\rho
_2|^2_g\, dV_g =\int _\Omega |d\rho _2|^2_h\, dV_h <\infty .
$$
Moreover, the $L^2/L^\infty $ comparison for \holo $1$-forms on a
bounded geometry K\"ahler manifold shows that $|d\rho _1|_g$ is
bounded on $E\supset \Omega $. Lemma~\ref{Cup product for sublevel
of plh function lemma} now gives the lemma in this case as well.
\end{pf*}

For the proof of Theorem~\ref{General filtered ends theorem}, we
will consider the cases of $e(X)\geq 2$ and $e(X)=1$ separately.
\begin{lem}\label{General filtered ends 2 end case lemma}
Let $(X,g)$ be a \con complete K\"ahler manifold which admits a
special ends decomposition. Assume that $e(X)\geq 2$ and $\tilde
e(X)\geq 3$. Then $X$ admits a proper \holo mapping onto a Riemann
surface.
\end{lem}
\begin{pf}
By Lemma~\ref{Open set to Riemann surface gives global lemma}, it
suffices to find a nonempty open subset of $X$ that admits a
proper \holo mapping onto a Riemann surface.

If $e(X)\geq 3$, then Theorem~\ref{Old 3 ends theorem} provides
the required proper \holo mapping to a Riemann surface.

If $e(X)=2$, then there exists a special ends decomposition
$X\setminus K=E_1\cup E_2$ where $E_1$ and $E_2$ are $\cinf $
domains with $K=\partial E_1=\partial E_2$ and, for some point
$x_0\in E_1$, the image $\Gamma $ of $\pi _1(E_1,x_0)$ in $\pi
_1(X,x_0)$ is a proper subgroup. We may assume without loss of
generality that, for $j=1,2$, $E_j$ is a hyperbolic special end
(of type (BG) or (RH)). For if $E_j$ is of type (W), then we may
apply Lemma~\ref{Type W to RH lemma}, while, if $E_j$ is of type
(SP), then, by Lemma~\ref{Special ends to produce a weakly
1-complete end lemma}, $E_j$ is also of type (W). Finally, a
parabolic end of type (BG) is also of type (SP).

Applying Proposition~\ref{Massive sets give harmonic functions
prop} and Remark 4 following Definition~\ref{Massive set
definition}, we get a finite energy \harmns , hence \plhns , \fn
$\rho :X\to (0,1)$ \st $\limsup _{x\to \infty }\rho
\restrict{\overline E_1}(x)=1$, $\liminf _{x\to \infty }\rho
\restrict{\overline E_2}(x)=0$, and $\lim _{j\to \infty } \rho
(x_j)=1$~$(\lim _{j\to \infty } \rho (x_j)=0)$ for any sequence
$\seq xj$ in $E_1$ (respectively $E_2$) such that $G(x_j,\cdot )
\to 0$ as $j\to \infty $. We may choose a \con covering space
$\Upsilon :\widehat X \to X$ \stns , for some point $y_0\in
\widehat X$, we have $\Upsilon _*\pi _1(\widehat X,y_0)=\Gamma $.
Thus $\Upsilon$ maps the \comp $E$ of $\widehat E_1=\Upsilon\inv
(E_1)$ isomorphically onto $E_1$ and, since $\Gamma $ is a proper
subgroup, $\widehat E_1\setminus E\neq \emptyset $. Since $E_1$ is
a $\cinf $ domain, loops in a small \nbd of $\overline E_1$
homotop into $E_1$. So $\Upsilon $ maps a \nbd of $\overline E$
isomorphically onto a \nbd of $\overline E_1$. Thus $E$ is a
hyperbolic special end in $(\widehat X,\hat g=\Upsilon ^*g)$. Fix
constants $a$ and $b$ with
$$
0<a<\min _{K=\partial E_1=\partial E_2}\rho \leq \max _{K}\rho
<b<1.
$$
In  $\widehat X$,
$$
\alpha _2\equiv \left\{
\begin{aligned}
\max (\rho \circ \Upsilon -b,0)
& \quad \text{on } U_2=\widehat E_1\setminus E \\
0 &\quad \text {on } \widehat X\setminus U_2
\end{aligned}
\right.
$$
is an admissible \plsh \fn for $U_2=\widehat E_1\setminus E$.
Choosing a $\cinf $ \fn $\chi :\R \to \R $ \st $\chi'\geq 0$,
$\chi ''\geq 0$, $\chi (t)=0$ for $t\leq 1-a$, and $\chi '(t)>0$
for $t>1-a$, we get a $\cinf $  admissible \plsh \fn
$$
\alpha _1\equiv \left\{
\begin{aligned}
\chi (1-\rho \circ \Upsilon )
& \quad \text{on } U_1=\Upsilon \inv(E_2) \\
0 &\quad \text {on } \widehat X\setminus U_1
\end{aligned}
\right.
$$
for $U_1\equiv\Upsilon\inv (E_2)$. Applying Lemma~\ref{Two
disjoint massive sets gives map lemma}, we get a proper \holo
mapping of a nonempty open subset of $E\cong E_1$ onto a Riemann
surface; as required.
\end{pf}

\begin{pf*}{Proof of Theorem~\ref{General filtered ends theorem}}
By Lemma~\ref{General filtered ends 2 end case lemma}, it remains
to consider the case $e(X)=1$. The manifold $X$ is itself a
special end in this case, so we have an ends decomposition
$X\setminus K=E_1$ \st $E_1$ is a $\cinf $ domain and, for a point
$x_0\in E_1$, $\Gamma\equiv\image{\pi _1(E_1,x_0)\to\pi_1(X,x_0)}$
is of index $\geq 3$.  We may fix a \con covering space $\Upsilon
:\widehat X\to X$ and a point $y_0\in \widehat X$ \st
$\Upsilon_*\pi _1(\widehat X,y_0)=\Gamma$. Hence $\Upsilon $ maps
the \comp $E$ of $\widehat E_1=\Upsilon\inv (E_1)$ containing
$y_0$ isomorphically onto $E_1$ and, since $\#\Upsilon\inv
(x_0)=[\pi _1(X,x_0):\Gamma ]\geq 3$, we have $\widehat
E_1\setminus E\neq\emptyset$. Again, since $E_1$ is a $\cinf $
domain, $\Upsilon $ maps a \nbd of $\overline E$ isomorphically
onto a \nbd of $\overline E_1$ and hence $E$ is a special end in
$(\widehat X,\hat g=\Upsilon ^*g)$.

We may again assume without loss of generality that $X$ (and hence
any end in $X$) is hyperbolic of type (BG) or (RH). For if $X$ is
of type (W), then we may apply Lemma~\ref{Type W to RH lemma}. If
$X$ is of type (SP), then $E$ is a special end of type (SP) in
$\widehat X$ and any other end in $\widehat X$ is either a
hyperbolic end or a special end of type (SP). Therefore, by
Lemma~\ref{Special ends to produce a weakly 1-complete end lemma},
there exists a \cont \plsh \fn $\hat\psi $ on $\widehat X$ which
exhausts $\overline E$ and hence the \fn
$$
\psi\equiv \left\{
\begin{aligned}
\max (\hat \psi \circ (\Upsilon\restrict E)\inv ,\max _{\partial
E}\hat\psi )
& \quad \text{on } E_1 \\
\max _{\partial E}\hat\psi &\quad \text {on } K
\end{aligned}
\right.
$$
is a \cont \plsh \exh \fn on $X$ and the above applies. Finally,
if $X$ is parabolic of type (BG), then $X$ is of type (SP) and the
above applies. Thus we may assume without loss of generality that
$X$ is hyperbolic of type (BG) or (RH).

In particular, $E_1$ is $D$-massive and there is a finite energy
admissible sub\harm \fn $\beta :X\to [0,1)$ for $E_1$ in $X$ \st
$\beta (x_j)\to 1$ as $j\to \infty $ for any sequence $\seq xj$ in
$E_1$ with $G(x_j,\cdot )\to 0$ as $j\to \infty $. It follows that
$E$ is $D$-massive in $\widehat X$ with finite energy admissible
sub\harm \fn
$$
\left\{
\begin{aligned}
\beta \circ \Upsilon
& \quad \text{on } E \\
0&\quad \text {on } \widehat X\setminus E
\end{aligned}
\right.
$$
and $\widehat E_1\setminus E$ is massive in $\widehat X$ with
admissible sub\harm \fn
$$
\left\{
\begin{aligned}
\beta \circ \Upsilon
&\quad \text{on } \widehat E_1\setminus E \\
0 &\quad \text {on } \widehat X\setminus (\widehat E_1\setminus
E)
\end{aligned}
\right.
$$
Applying Proposition~\ref{Massive sets give harmonic functions
prop}, we get a finite energy \harmns , hence \plhns , \fn $\rho
_1:\widehat X\to (0,1)$ \st $\beta \circ \Upsilon \leq \rho _1 <1$
on $E$ and \stns , for any admissible sub\harm \fn $\gamma
:\widehat X \to [0,1)$ for $\widehat X\setminus \overline E$, we
have $0<\rho _1\leq 1-\gamma $ on $\widehat X$.

Fix constants $t_0$,~$t_1$,~and~$t_2$ with $\max _{\partial E}\rho
_1 <t_2<t_1<t_0<1$ and a $\cinf $ \fn $\chi :\R \to \R $ \st $\chi
'\geq 0$, $\chi ''\geq 0$, $\chi (t)=0$ if $t\leq t_2$, and $\chi
(t)=(t-t_1)/(1-t_1)$ if $t\geq t_0$. Thus
$$
\beta _1=\left\{
\begin{aligned}
\chi \bigl( \rho _1\circ (\Upsilon \restrict E)\inv \bigr)
& \quad \text{on } E_1 \\
0&\quad \text {on } X\setminus E_1
\end{aligned}
\right.
$$
is a finite energy $\cinf $ admissible \plsh \fn for $E_1$ which
is \plh on $\setof{x\in X}{\beta (x)>t_0}$. For we have, on $E_1$,
$$
\big|d\chi \bigl( \rho _1\circ (\Upsilon \restrict E)\inv
\bigr)\big|_g =\chi '\bigl( \rho _1\circ (\Upsilon \restrict
E)\inv \bigr) |d\rho _1|_{\hat g}\circ (\Upsilon \restrict E)\inv
\leq (1-t_1)\inv |d\rho _1|_{\hat g}\circ (\Upsilon \restrict
E)\inv
$$
and $\rho _1$ has finite $\hat g$-energy, so $\chi \bigl( \rho
_1\circ (\Upsilon \restrict E)\inv \bigr) $ has finite $g$-energy.
If $\seq xj$ is a sequence in $E_1$ with $G(x_j,\cdot )\to 0$ as
$j\to \infty $, then
$$
1>\chi \bigl( \rho _1((\Upsilon \restrict E)\inv (x_j))\bigr) \geq
\chi (\beta (x_j)) \to \chi (1)=1,
$$
so $\chi \bigl( \rho _1((\Upsilon \restrict E)\inv (x_j))\bigr)
\to 1$. Finally,
$$
1>\rho _1\geq \chi (\rho _1) \quad (\geq \chi (\beta \circ
\Upsilon )) \quad \text{on }E,
$$
so the relation $\rho _1\geq \beta \circ \Upsilon $ on $E$ is
preserved if we replace $\beta $ with $\beta _1$. Thus we may
assume that $\beta =\chi \bigl( \rho _1\circ (\Upsilon \restrict
E)\inv \bigr)$ on $E_1$.

If $U$ is any \comp of $\widehat E_1=\Upsilon\inv (E_1)$, then the
\fn
$$
\left\{
\begin{aligned}
\beta \circ \Upsilon
& \quad \text{on } U \\
0&\quad \text {on } \widehat X\setminus U
\end{aligned}
\right.
$$
is a $\cinf $ admissible \plsh \fn for $U$ (with finite energy if
$U\to E_1$ is a finite covering). Therefore, if $\widehat E_1$ has
at least 3 \compsns , then Lemma~\ref{Two disjoint massive sets
gives map lemma} provides an open subset of $E\cong E_1$ which
admits a proper \holo mapping onto a Riemann surface. If $\widehat
E_1$ has exactly two \comps $E$ and $E'$ and $E'\to E_1$ is a
finite covering, then we see that $\widehat X\to X$ is a finite
covering and $\widehat X$ has the special ends decomposition
$\widehat X\setminus \widehat K=E\cup E'$, where $\widehat
K=\Upsilon\inv (K)$. Thus $\tilde e(\widehat X)=\tilde e(X)\geq 3$
and $e(\widehat X)\geq 2$, and therefore, by Lemma~\ref{General
filtered ends 2 end case lemma}, $\widehat X$ admits a proper
\holo mapping $\Psi :\widehat X\to S$ with \con fibers onto a
Riemann surface $S$. The \plh \fn $\rho _1$ descends to a \plh \fn
$\tau $ on $S$. Hence any \comp of $\setof{x\in \widehat X}{\rho
_1(x)>\max _{\partial E}\rho _1}$ contained in $E\cong E_1$ admits
a proper \holo mapping onto a \comp of $\setof{\zeta\in S}{\tau
(\zeta )>\max _{\partial E}\rho _1}$.

Thus it remains to consider the case in which $\widehat
E_1=\Upsilon\inv (E_1)$ has exactly two \comps $E$, $E'$ and
$E'\to E_1$ is an infinite covering (i.e. $[\pi _1(X,x_0):\Gamma
]=\infty $). We have the finite energy \plh \fn $\rho _1:\widehat
X\to (0,1)$ with
$$
\beta \circ \Upsilon \leq \rho _1<1 \text{ on }E
\qquad\text{and}\qquad 0<\rho _1\leq 1-\gamma
$$
for every admissible sub\harm \fn $\gamma :X\to [0,1)$ for
$X\setminus \overline E$. By construction, for some constants
$t_0$, $t_1$, and $t_2$ with $\max _{\partial E}\rho _1
<t_2<t_1<t_0<1$, we have on $E$
$$
\beta\circ \Upsilon =\left\{
\begin{aligned}
0
& \quad \text{on } \setof{x\in E}{\rho _1(x)\leq t_2}
\\
\frac{\rho _1-t_1}{1-t_1}&\quad \text {on } \setof{x\in E}{\rho
_1(x)\geq t_0}=\setof{x\in E}{\beta (\Upsilon (x))\geq r_0}
\end{aligned}
\right.
$$
where $r_0=(t_0-t_1)/(1-t_1)$. We will produce a second \plh \fn
and apply Lemma~\ref{Cup product for compact fiber lemma} and
Lemma~\ref{Cup product for sublevel of plh function lemma}.
Lifting, we get $\cinf $ admissible \plsh \fns
$$
\alpha =\left\{
\begin{aligned}
\beta\circ\Upsilon & \quad \text{on } E\\
0&\quad \text {on } \widehat X\setminus E
\end{aligned}
\right.
\qquad\text{and}\qquad\alpha '=\left\{
\begin{aligned}
\beta\circ\Upsilon &\quad \text{on } E'\\
0&\quad \text {on } {\widehat X}\setminus E'
\end{aligned}
\right.
$$
for $E$ and $E'$, respectively. Moreover, $\alpha $ has {\it
finite} energy, $\alpha '$ has {\it infinite} energy, and
$$
\alpha =\frac {\rho _1-t_1}{1-t_1} \text{ on } \setof{x\in
\widehat X}{\alpha (x)\geq r_0}
$$
and
$$\alpha '=\frac {\rho
_1\circ (\Upsilon \restrict E)\inv\circ\Upsilon -t_1}{1-t_1}
\text{ on } \setof{x\in \widehat X}{\alpha '(x)\geq r_0}.
$$

Let $V=\setof{x\in X}{\beta (x)>r_0}\subset E_1$ and let $\widehat
V=\Upsilon\inv(V)\subset\widehat E_1$. If $\widehat V\cap E'$ is
not \conns , then we get two disjoint $\cinf $ plurimassive
subsets of $\widehat X$ contained in $\widehat V\cap E'\subset
\widehat X\setminus\overline E$ and we may apply Lemma~\ref{Two
disjoint massive sets gives map lemma} as before. Thus we may
assume without loss of generality that $\widehat V\cap E'$ is
\conns . In particular, $\widehat V\cap E'\to V$ is a \con
infinite covering. Observe also that $\alpha $ and $\alpha '$ are
\plh on $\widehat V$.

For each $r$ with $r_0<r<1$, the \comp $Y_r$ of $\setof{x\in
\widehat X}{\alpha '(x)<r}$ containing $E$ (and $\overline E$)
admits the complete K\"ahler metric
$$
h_r=\hat g+\lev{-\log (r-\alpha ')}=\hat g+(r-\alpha ')\inv \lev
{\alpha '}+(r-\alpha ')^{-2}\partial \alpha
'\overline{\partial\alpha '}
$$
with $h_r\geq \hat g$ on $\widehat X$ and $h_r=\hat g$ at points
in $\setof{x\in \widehat X}{\alpha '(x)=0}\supset \widehat
X\setminus E'$. Moreover, $\alpha \restrict{Y_r}$ and $\alpha
'\restrict{Y_r}$ are $\cinf $ admissible \plsh \fns for $E=E\cap
Y_r$ and $E'_r=E'\cap Y_r$, respectively, and $\alpha
\restrict{Y_r}$ has finite $h_r$-energy (since $\alpha =0$ on
$\widehat X\setminus E$ and $h_r=\hat g$ on $E$). Applying
Proposition~\ref{Massive sets give harmonic functions prop}, we
get a finite $h_r$-energy \harmns , hence \plhns , \fn $\tau
_r:Y_r\to (0,1)$ \st $\alpha \leq \tau _r<1$ on $Y_r$ and \stns ,
for any admissible $h_r$-sub\harm \fn $\gamma :Y_r\to [0,1)$ for
$Y_r\setminus\overline E$, $0<\tau _r\leq 1-\gamma $ on $Y_r$.

We will show that the \fns $1$, $\rho _1\restrict{Y_r}$, and $\tau
_r$ are linearly independent for some $r$ with $r_0<r<1$. To see
this, suppose that, on the contrary, $1$, $\rho _1\restrict{Y_r}$,
and $\tau _r$ are linearly dependent \fns for every $r$ with
$r_0<r<1$. Then, since $0<\tau _r\leq 1-r\inv\alpha ' \to 0$ at
$\partial Y_r$, we see that $\rho _1$ is constant on $\partial
Y_r\subset \widehat V\cap E'$ for each $r\in (r_0,1)$. Fixing a
regular value $r\in (r_0,1)$ for $\alpha '$ and a point
$p\in\partial Y_r$, we may choose a \holo coordinate \nbd $(W,\Phi
=(\zeta _1,\dots ,\zeta _n))$ mapping $W$ onto a $2n$-dimensional
open rectangle $(r-\epsilon,r+\epsilon )\times
(0,1)\times\cdots\times (0,1)$ in $\R ^{2n}$ with $\Phi
(p)=(r,0,\dots ,0)$ and $\alpha '=\real \zeta _1=u_1$ on $W$. For
each $s$ with $r<s<r+\epsilon $, $\setof{x\in W}{\alpha
'(x)<s}\cong (r-\epsilon,s)\times (0,1)\times\cdots\times (0,1)$
is a \con open subset of $\setof{x\in \widehat X}{\alpha '(x)<s}$
meeting $Y_r\subset Y_s$, so this set is contained in $Y_s$. Thus
$W\cap Y_s=\setof{x\in W}{\alpha '(x)<s}$ and so $\rho _1$ is
constant on $W\cap\partial Y_s=\setof{x\in W}{\alpha '(x)=s}$. It
follows that $d\rho _1\wedge d\alpha '\equiv 0$ on $\setof{x\in
W}{r\leq \alpha '(x)< r+\epsilon }$ and, therefore, since $\rho
_1$ and $\alpha '$ are \plh on the \con open set $\widehat V\cap
E'$, we see that $d\rho _1\restrict{(\widehat V\cap E')}$ and
$d\alpha '\restrict{(\widehat V\cap E')}$ are linearly dependent
forms. On the other hand, $\widehat V\cap E'\to V$ is an infinite
covering. Thus we have, for some constant $C>0$,
$$
\infty >\int _{\widehat V\cap E'}|d\rho _1|^2_{\hat g}\, dV_{\hat
g}=C\int _{\widehat V\cap E'}|d\alpha '|^2_{\hat g}\, dV_{\hat g}
=C\int _{\widehat V\cap E'}|\Upsilon ^*d\beta |^2_{\hat g}\,
dV_{\hat g} =\infty .
$$
Thus we have arrived at a contradiction, so $1$, $\rho
_1\restrict{Y_r}$, and $\tau _r$ must be linearly independent \fns
on $Y_r$ for some $r\in (r_0,1)$.

As in the last part of the proof of Lemma~\ref{Two disjoint
massive sets gives map lemma} (applied to $\rho _1$ and $\tau
_r)$, for $\max _{\partial E}\rho _1<a<1$, one gets a proper \holo
mapping of some open subset of a \comp of $\setof{x\in X}{a<\rho
_1(x)}$ contained in $E\cong E_1$ onto a Riemann surface and the
theorem follows.
\end{pf*}
\begin{rmks}
1. If $M$ is a \con non\cpt manifold with $3\leq\tilde
e(M)<\infty$, then $M$ admits a finite covering space $\widehat M$
with $e(\widehat M)\geq 3$. To see this, we fix an ends
decomposition $M\setminus K=E_1\cup\cdots\cup E_m$ \st the lifting
of $M\setminus K$ to the universal covering $\Upsilon :\widetilde
M\to M$ has $\tilde e(M)$ \compsns . The action of $\pi _1(M)$
permutes these \comps and so we get a homomorphism of $\pi _1(M)$
into the symmetric group on $\tilde e(M)$ objects.  Thus the
kernel $\Gamma$ is a normal subgroup of finite index and hence the
quotient $\widehat M=\Gamma\big\backslash \widetilde M \to M$ is
the desired finite covering. This observation together with
Theorem~\ref{Old 3 ends theorem} gives Theorem~\ref{General
filtered ends theorem} for the case in which $3\leq\tilde
e(X)<\infty$. For the associated finite covering space admits a
proper \holo mapping onto a Riemann surface and any normal \cpx
space which is the image of a \holoconvex \cpx space under a
proper \holo mapping is itself \holoconvexns .  In particular, in
the proof of Theorem~\ref{General filtered ends theorem}, we could
have avoided the argument in the case in which $e(X)=1$, $\widehat
E_1=\Upsilon\inv (E_1)$ has exactly two \comps $E$, $E'$, and
$E'\to E_1$ is a finite covering.

\noindent 2. In general, the number of components of $\Upsilon\inv
(E_1)$ is equal to the number of distinct double cosets of the
image of the fundamental group of $E_1$.
\end{rmks}

Theorem~\ref{General filtered ends theorem} and Lemma~\ref{Two
ends and infinitely generated lemma} give the following:
\begin{thm}\label{General special ends, two ends-infinitely generated maps to Riemann surface theorem}
Let $X$ be a \con complete K\"ahler manifold which admits a
special ends decomposition. Assume that $\tilde e(X)\geq 2$ and
$\pi _1(X)$ is infinitely generated. Then $X$ admits a proper
\holo mapping onto a Riemann surface.
\end{thm}

\section{Mappings of \cpt K\"ahler manifolds to curves}\label{Mappings of Kahler manifolds to
curves section}

In this section, we consider the following consequence of
Theorem~\ref{Main Theorem from Introduction} (and
Theorem~\ref{General filtered ends theorem}).

\begin{thm}\label{Mapping of compact to curve theorem}
If $X$ is a \con \cpt K\"ahler manifold for which there is a \con
infinite covering space $\Upsilon :\widehat X\to X$ with $\tilde
e(\widehat X)\geq 3$, then some finite covering space $X'\to X$
admits a surjective \holo mapping onto a curve of genus $g\geq 2$.
\end{thm}

The main point of the proof is the following
(cf.~Proposition~1.2.11 of \cite{Kollar} and 1.2.3, p.~490, of
\cite{Campana}):
\begin{prop}\label{Kollar normalizer prop}
Let $(X,g)$ be a \con complete K\"ahler manifold with bounded
geometry. Suppose some \con non\cpt covering space $\Upsilon
:\widehat X\to X$ admits a proper surjective \holo mapping
$\hat\Phi :\widehat X\to \widehat S$ onto a Riemann surface
$\widehat S$. Then, for every level $F$ of $\hat\Phi $ over a
regular value $\zeta $ of $\hat\Phi $, the normalizer~$N$ of
$\text{{\rm im}}\, [\pi _1(F)\to \pi _1(X)]$ is of finite index in
$\pi _1(X)$. Furthermore, the associated finite covering space
$X'\to X$ with $\text{{\rm im}}\, [\pi _1(X')\to \pi _1(X)]=N$
admits a surjective proper \holo mapping $\Phi ':X'\to S$ onto a
Riemann surface $S$. In particular, if $X$ is non\cptns, then $X$
admits a proper \holo mapping onto a Riemann surface.
\end{prop}

\begin{rmk}
The last statement is a consequence of the fact that the normal
proper \holo image of a \holoconvex \cpx space is \holoconvexns .
\end{rmk}

The following elegant proof of Proposition~\ref{Kollar normalizer
prop}, which we give in steps, is due to Delzant and Gromov
\cite{DG}. Let $\Phi :X\to S$ be a proper \holo mapping with \con
fibers of a \con \cpx manifold $X$ onto a Riemann surface $S$, let
$X_A=\Phi\inv (A)$ for each set $A\subset S$, and let $X_\zeta
=X_{\{ \zeta \} }$ for each point $\zeta \in S$. If $P$ is the
(discrete) set of critical values of $\Phi $ and, for each point
$p\in P$, $m_p$ is the greatest common divisor of the
multiplicities of the \comps of the divisor $\Phi\inv (p)$ and
$\gamma _p$ is a simple loop tracing the boundary circle of a
coordinate disk $D_p$ in $S$ centered at $p$ with $D_p\cap P=\{
p\} $, then the {\it orbifold fundamental group} is given by
$$
\pi _1^{\text{orb}}(\Phi )\equiv \pi _1 (S\setminus P)/N
$$
where $N$ is the normal subgroup of $\pi _1 (S\setminus P)$
generated by the loops $\{ \gamma _p^{m_p}\} _{p\in P}$. The
following lemma is well known (see \cite{Cat2}, \cite{CatKO},
\cite{Simpson-VHS}):
\begin{lem}\label{Orbifold fundamental group lemma}
For $\Phi :X\to S\supset P$ as above, we have
\begin{enumerate}
\item[(a)] The map $\pi _1(X)\to\pi_1(S)$ is surjective.
\item[(b)] For each point $\zeta \in S\setminus P$, the induced
maps give an exact sequence
$$
\pi _1(X_\zeta )\to \pi _1(X)\to \pi _1^{\text{orb}}(\Phi )\to 1.
$$
\item[(c)] For each pair of points $\zeta _1,\zeta _2\in
S\setminus P$, each choice of a point $x_j\in X_{\zeta _j}$ for
$j=1,2$, and each path $\alpha $ in $X_{S\setminus P}$ from $x_1$
to $x_2$, the isomorphism $\pi _1(X,x_1)\to \pi _1(X,x_2)$ given
by $[\gamma ]\to [\alpha\inv *\gamma *\alpha ]$ restricts to a
surjective isomorphism of $\Gamma _1\equiv \text{{\rm im}}[\pi
_1(X_{\zeta _1},x_1)\to \pi _1(X,x_1)]$ onto $\Gamma _2\equiv
\text{{\rm im}}[\pi _1(X_{\zeta _2},x_2)\to \pi _1(X,x_2)]$.
\item[(d)] If $\Upsilon :X'\to X$ is a \con covering space and,
for some point $\zeta_1\in S\setminus P$, the group $\Gamma
_1\equiv\text{{\rm im}}[\pi _1(X_{\zeta _1 })\to \pi _1(X)]$ is
contained in (is equal to) the group $\Lambda\equiv\text{{\rm
im}}[\pi _1(X')\to \pi _1(X)]$, then this is the case for every
point $\zeta \in S\setminus P$ and we get a commutative diagram
\begin{center}\begin{picture}(230,70)
\put(80,9){$X$} \put(93,13){\vector(1,0){46}} \put(143,9){$S$}
\put(80,49){$X'$} \put(113,54){$\Phi '$}
\put(93,53){\vector(1,0){46}} \put(143,49){$S'$}
\put(152,31){$\theta $} \put(147,47){\vector(0,-1){28}}
\put(72,31){$\Upsilon $} \put(84,47){\vector(0,-1){28}}
\put(113,15){$\Phi $}
\end{picture}
\end{center}
where $S'$ is a Riemann surface, $\Phi '$ is a surjective proper
\holo map with \con fibers, and $\theta $ is a (possibly infinite,
but locally finite) \holo branched covering map with branch locus
$B$ contained in $\setof{p\in P}{m_p>1}$. In particular, if
$\Gamma _1=\Lambda $, then $S'\cong {\mathbb P} ^1$, $\C$, or
$\Delta $ and we get injective homomorphisms of $\pi
_1^{\text{orb}}(\Phi )=\pi _1(X)/\Gamma _1$ into $\text{{\rm
Aut}}\, (X')$ and into $\text{{\rm Aut}}\, (S')$, $\Phi '$ is
equivariant \wrt the action of $\pi _1^{\text{orb}}(\Phi )$, and
$\Upsilon $ and $\theta $ are the quotient maps given by
$$
\Upsilon : X'\to X=\pi _1^{\text{orb}}(\Phi )\lquotient X'
\qquad\text{and}\qquad \theta : S'\to S=\pi _1^{\text{orb}}(\Phi
)\lquotient S'.
$$
\end{enumerate}
\end{lem}

\begin{lem}\label{Finite collection of images of fundamental groups lemma}
Let $(M,g)$ be a \con complete Riemannian manifold.
\begin{enumerate}
\item[(a)] For each point  $p\in M$ and each constant $L>0$, the
set
$$
K(p,L)\equiv \setof{[\alpha ]\in \pi _1(M,p)}{\alpha \text{ is a
piecewise }\cinf \text{ loop in }M \text{ of length }<L}
$$
is finite. \item[(b)] Let $A$ be a path \con \cpt subset of $M$,
let $r>0$ be a lower bound for the injectivity radius at points in
$A$, and let $A_1,\dots ,A_m$ be a covering of $A$ by path \con
subsets which are relatively open in $A$ and which have diameter
$<r$ (with respect to the distance \fn in $M$). Then, for each
point $p\in A$, the group $\Gamma \equiv \text{\rm im}\, [\pi
_1(A,p)\to\pi _1(M,p)] $ is generated by the set
$$
\setof{ [\alpha ]\in\Gamma }{\alpha \text{ is a piecewise }\cinf
\text{ loop in }M \text{ based at }p \text{ of
length}<2(m^2+1)r+1}.
$$
\item[(c)] Let $\Upsilon :\widehat M\to M$ be a \con covering
space, let $\hat g=\Upsilon ^*g$, and let $\seq A\lambda _{\lambda
\in \Lambda}$ be a family of path \con \cpt subsets of $\widehat
M$. Assume that there exist a positive integer $m$ and a positive
constant $r$ \stns , for each $\lambda \in \Lambda $, $r$ is a
lower bound for the injectivity radius in $\widehat M$ at each
point in $A_\lambda $ and there is a covering $A_1^\lambda ,\dots
, A_m^\lambda $ of $A$ by path \con relatively open subsets of $A$
of diameter $<r$ (in $\widehat M$). Then, for each point $p\in M$,
the (possibly empty) collection of subgroups
$$
{\mathcal H}_p\equiv\setof{\text{\rm im}\, [\pi _1(A_\lambda ,\hat
p)\to\pi _1(M,p)]}{\lambda\in \Lambda, \hat p\in \Upsilon\inv
(p)\cap A_\lambda }
$$
is finite.
\end{enumerate}
\end{lem}
\begin{pf*}{Sketch of the proof}
For the proof of (a), we fix a number $r>0$ with $3r$ less than
the injectivity radius at each point in the set
$D\equiv\overline{B(p,L)}$, points $p=p_1,p_2,\dots ,p_k\in D$ \st
the balls $B_1=B(p_1,r),\dots , B_k=B(p_k,r)$ form a covering for
$D$, a Lebesgue number $\delta >0$ for this covering, and a
positive integer $m$ \st $L/m<\delta$. For each pair of indices
$i,j$, we let $\lambda_{ij}=\lambda_{ji}\inv$ be a minimal
geodesic from $p_i$ to $p_j$. Now any piecewise $\cinf$ loop
$\alpha$ of length $<L$ based at $p$ is homotopic to a loop
$\alpha _1 *\alpha _2*\cdots *\alpha _m$ in $D$; where, for each
$\nu=1,\dots ,m$, $\alpha _\nu$ is a piecewise $\cinf$ path of
length $<\delta$ and is, therefore, contained in $B_{i_\nu}$ for
some index $i_\nu$. We may assume that $i_1=i_m=1$. Thus $\alpha$
is homotopic to the loop $\lambda_{i_1i_2}*\lambda_{i_2i_3}*\cdots
*\lambda_{i_{m-1}i_m}$ and the claim follows.

For the proof of (b), we let $I$ be the set of pairs of
indices $(i,j)$ with $1\leq i,j\leq m$ and  $A_i\cap
A_j\neq\emptyset $ and, for each $(i,j)\in I$, we fix a point
$p_{ij}=p_{ji}\in A_i\cap A_j$. If $(i,j), (j,k)\in I$, then
$A_j\subset B_g(p_{ij},r)$ and we get a unique minimal geodesic
$\gamma _{ijk}=\gamma _{kji}\inv $ from $p_{ij}$ to $p_{jk}$.

Fix a point $p\in A$. Given a point $q\in A$, we may form a broken
geodesic $\lambda $ of the form
$$
\lambda =\zeta *\gamma _{i_0i_1i_2}* \gamma _{i_1i_2i_3}*\cdots
*\gamma _{i_{k-2}i_{k-1}i_k}*\eta
$$
from $p$ to $q$ where $p\in A_{i_0}\subset B(p_{i_0i_1},r)$, $q\in
A_{i_k}\subset B(p_{i_{k-1}i_k},r)$, and $\zeta $ and $\eta $ are
the unique minimal geodesics from $p$ to $p_{i_0i_1}$ and from
$p_{i_{k-1}i_k}$ to $q$, respectively. On the other hand, any
broken geodesic $\lambda $ of the above form is homotopic to a
path in $A$. If we choose $\lambda $ so that $k$ is minimal, then
each pair $(i,j)\in I$ can be equal to $(i_{\nu -1},i_\nu )$ for
at most one $\nu \in \{  1,2,3,\dots ,k\} $. Thus we must have
$k\leq m^2$ and hence $l_g(\lambda )<r+(m^2-1)r +r=(m^2+1)r$. Thus
we see that any point $q\in A$ may be joined to $p$ by a piecewise
$\cinf $ path $\lambda $ in $M$ which has length $<(m^2+1)r$ and
which is homotopic (in $M$) to a path in $A$.

Now given a loop $\beta $ in $A$ based at $p$, we may choose a
subdivision $0=t_0<t_1<\dots <t_k=1$ \stns , for $\nu =1,\dots
,k$, $\beta \restrict{[t_{\nu -1,t_\nu ]}}$ is homotopic to a
$\cinf $ path of length $<1$. By the above, for each $\nu =1,\dots
,k-1$, we may also choose a piecewise $\cinf $ path $\lambda _\nu
$ from $p$ to $\beta (t_\nu )$ \st   $l(\lambda _\nu )<(m^2+1)r$
and $\lambda _\nu $ is homotopic to a path in $A$. Thus, in
$\Gamma $, we have
$$
[\beta ]=[\beta _1]\cdots [\beta _k]
$$
where $\beta _1=\beta \restrict{[t_0,t_1]}*\lambda _1\inv $,
$\beta _\nu =\lambda _{\nu -1}*\beta \restrict{[t_{\nu -1,t_\nu
]}}*\lambda _{\nu }\inv $ for $\nu =2,\dots ,k-1$, and $\beta
_k=\lambda _{k-1}*\beta \restrict{[t_{k-1},t_k]}$. Since each of
the above loops is homotopic to a piecewise $\cinf$ loop of length
$<(m^2+1)r+1+(m^2+1)r$, the claim (b) follows.

Finally, for the proof of (c), we fix $p\in M$ and we set
$L=2(m^2+1)r+1$. Applying Part~(b) in the covering space
$(\widehat M,\hat g)$, we see that each element $\Gamma \in
{\mathcal H}_p$ is generated by the set $K(p,L)\cap \Gamma $.
Therefore, since $K(p,L)$ is a finite set by Part~(a), ${\mathcal
H}_p$ must be finite.
\end{pf*}

\begin{pf*}{Proof of Proposition \ref{Kollar normalizer prop} (Delzant-Gromov \cite{DG})}
Stein factoring, we may assume that $\hat\Phi $ has \con fibers.
Furthermore, by passing to the appropriate covering space as in
part (d) of Lemma~\ref{Orbifold fundamental group lemma}, we may
assume that $\pi _1(\widehat X_\zeta )\to \pi _1(\widehat X)$ is
surjective for every regular value $\zeta $ of $\hat\Phi $ and,
therefore, for any $\zeta\in \widehat S$ (since $\pi _1(\widehat
X_\zeta )$ surjects onto $\image{\pi _1(X_U)\to \pi _1(X)}=\pi
_1(X)$ for a sufficiently small \con \nbd $U$ of $\zeta $ in
$\widehat S$). With these additional assumptions, we will show
that the normalizer of $\image{\pi _1(\widehat X)\to \pi _1(X)}$
is of finite index. Equivalently, the normalizer of $\image{\pi
_1(F)\to \pi _1(X)}$ is of finite index for {\it any} (possibly
singular) fiber $F$ of $\hat\Phi $.

Clearly, we may assume without loss of generality that $n=\dim
X>1$ and, since $\widehat X$ is non\cptns , we have $\widehat S=\C
$ or $\Delta $. Fixing a point $\zeta _0\in\widehat S$ and a point
$\hat x_0\in F_0\equiv \widehat X_{\zeta _0}=\hat\Phi\inv (\zeta
_0)$, we get
$$
\widehat \Gamma _0=\pi _1(\widehat X,\hat x_0)=\image{\pi
_1(F_0,\hat x_0)\to \pi _1(\widehat X,\hat x_0)}
\overset{\cong}{\to } \Gamma _0\equiv \Upsilon _*\widehat\Gamma
_0\subset \Lambda \equiv \pi _1(X,x_0),
$$
where $x_0=\Upsilon (\hat x_0)$. We must show that the normalizer
$$
N_0\equiv \setof {\lambda \in \Lambda }{\lambda \Gamma
_0\lambda\inv =\Gamma _0}
$$
is of finite index in $\Lambda $. For this, it suffices to show
that $\Gamma _0$ has only finitely many distinct conjugates in
$\Lambda $.

The collection of conjugates of $\Gamma _0$ in $\Lambda$ is
precisely the collection $\mathcal H$ of subgroups $\Gamma $ of
$\Lambda $ of the form
$$
\Gamma =\image{\pi _1(\widehat X,\hat x)\to \Lambda }=\image{\pi
_1(F,\hat x)\to \Lambda}
$$
where $\hat x\in\Upsilon\inv (x_0)$ and $F$ is the (not
necessarily smooth) fiber of $\hat\Phi $ containing $\hat x$. Let
$\hat g=\Upsilon ^*g$. According to Lemma~\ref{Finite collection
of images of fundamental groups lemma}, to show that $\mathcal H$
is finite, it suffices to find constants $m\in\N $ and $r>0$
\stns, for each fiber $F$ of $\hat\Phi $ meeting $\Upsilon\inv
(x_0)$, $r$ is a lower bound for the injectivity radius in
$\widehat X$ at each point in $F$ and there is a covering
$A_1\dots , A_m$ of $F$ by \con relatively open subsets of $F$ of
diameter $<r$ (in $\widehat X$).

The covering space $\widehat X$ has \bdd geometry because $X$
does. Thus, for some constant $C>0$ and for each point
$p\in\widehat X$, there is a biholomorphism $\Psi _p$ of the unit
ball $B_{\C ^n}(0;1)\subset \C ^n$ onto a \nbd
$B_{E,p}=B_{E,p}(1)$ of $p$ in $\widehat X$ \st $\Psi_p(0)=p$ and
$$
C\inv \Psi _p^*\hat g\leq g_{\C ^n}\leq C\Psi _p^*\hat g
\qquad\text{on }B_{\C ^n}(0;1).
$$
We set $B_{E,p}(r)=\Psi _p(B_{\C ^n}(0;r))$ for each $r\in (0,1)$.
Thus, for constants $R_1$, $R_2$, and $R$ with $1/4\gg
R_1>R_1/2\gg R>R/2 \gg R_2>0$, we get, for each point
$p\in\widehat X$,
$$
B_{E,p}(1/4)\Supset B_{\hat g}(p,R_1)\Supset B_{\hat
g}(p,R_1/2)\Supset B_{E,p}(R)\Supset B_{E,p}(R/2)\Supset B_{\hat
g}(p,R_2).
$$

If $F$ is a fiber of $\hat\Phi $ meeting $\Upsilon\inv (x_0)$,
then we may choose points $p_1,\dots ,p_s\in F$ \st
$$
F\subset B_{\hat g}(p_1,R_1) \cup \cdots\cup B_{\hat g}(p_s,R_1)
$$
but
$$
p_j\in F\setminus [B_{\hat g}(p_1,R_1) \cup \cdots\cup B_{\hat
g}(p_{j-1},R_1)] \qquad\text{for } j=2,\dots ,s.
$$
In particular, the balls $B_{\hat g}(p_1,R_1/2),\dots ,B_{\hat
g}(p_s,R_1/2)$ are disjoint and, therefore, the balls
$B_{E,p_1}(R),\dots ,B_{E,p_s}(R)$ are disjoint. According to
\cite{Stoll}, since $\widehat X$ is K\"ahler, the fibers of
$\hat\Phi $ have equal volume $v_0$ (counting multiplicities).
But, by Lelong's monotonicity formula (see estimate~15.3 in
\cite{Chirka}), there is a constant $\delta >0$ \st any \anal set
$A$ in $B(0,1/2)$ of pure dimension $n-1$ passing through $0$
satisfies $\text{vol}_{\C ^n}(A\cap B_{\C ^n}(0,R))\geq\delta $
(recall that $R<1/4$). It follows that
\begin{align*} v_0=\text{vol}_{\hat g}(F)&\geq \sum
_{j=1}^s\text{vol}_{\hat g}(F\cap B_{E,p_j}(R)) \geq \sum _{j=1}^s
C^{-n}\text{vol}_{\C ^n}(\Psi _{p_j}\inv (F\cap
B_{E,p_j}(1/2))\cap B_{\C ^n}(0,R))
\\
&\geq C^{-n}\cdot s\cdot \delta .
\end{align*}
Thus we have the uniform bound $s\leq C^nv_0/\delta $ for $s$.
Consequently, we also get the uniform bound $\diam F\leq
2R_1C^nv_0/\delta $ for the diameter (with respect to the distance
in $\widehat X$) of $F$.

It follows that the union $Z$ of the images in $X$ of the fibers
of $\hat\Phi $ meeting $\Upsilon\inv (x_0)$
is a bounded set. Hence we may choose a constant
$r>0$ so that $r$ is a lower bound for the injectivity radius in
$X$ at each point in $Z$ and, therefore, for the injectivity
radius in $\widehat X$ at each point in any fiber of $\hat\Phi $
which meets $\Upsilon\inv (x_0)$. We may also choose each of the
\nbds ${\seq B{E,p} }_{p\in \widehat X}$ so that $\diam B_{E,p}<r$
for each $p\in\widehat X$ (although the above proof of the
boundedness of $Z$ involved the choice of ${\seq B{E,p} }_{p\in
\widehat X}$, this boundedness property is clearly independent of
the choice of ${\seq B{E,p} }_{p\in \widehat X}$).

Now suppose that $F$ is again a fiber of $\hat\Phi $
meeting $\Upsilon\inv (x_0)$ and, in the
above notation,  $1\leq j\leq s$ and $A$ is a \con \comp of $F\cap
B_{E,p_j}$ which meets $B_{\hat g}(p_j;R_1)$. Choosing a point
$a\in A\cap B_{\hat g}(p_j,R_1)$ and applying the above volume
estimate then gives
$$
\text{vol}_{\hat g}(A)\geq  C^{-n}\text{vol}_{\C ^n}((\Psi
_{p_j}\inv (A)\cap B_{\C ^n}(\Psi _{p_j}\inv (a);1/2))\cap B_{\C
^n}(\Psi _{p_j}\inv (a);R))
\\
\geq C^{-n}\cdot \delta .
$$
It follows that $F\cap B_{E,p_j}$ can have at most $C^nv_0/\delta
$ such components. Combining this uniform bound with the uniform
bound for $s$, we see that, if we fix a positive integer $m\geq
(C^nv_0/\delta )^2$, then, for any such $F$, we may choose sets
$A_1,\dots ,A_m$ so that $F=A_1\cup \cdots \cup A_m$ (possibly
with some repetition) and, for each $\nu =1,\dots ,m$, $A_\nu $ a
\comp of $F\cap B_{E,p_j}$ (which meets $B_{\hat g}(p_j,R_1)$) for
some $j$. The finiteness of $\mathcal H$, and hence of $[\Lambda
:N_0]$, now follows.

Thus we get a commutative diagram
\begin{center}\begin{picture}(230,70)
 \put(139,9){$X$}
\put(80,49){$\widehat X$} \put(110,56){$\Upsilon '$}
\put(90,53){\vector(1,0){46}} \put(139,49){$X'=\Lambda
'\big\backslash \widehat X$} \put(147,31){$\Upsilon '' $}
\put(144,47){\vector(0,-1){27}} \put(110,20){$\Upsilon $}
\put(85,47){\vector(2,-1){55}}
\end{picture}
\end{center}
in which $\Upsilon '$ is a Galois covering map, $\Upsilon ''$ is a
finite covering map, and $\Lambda '=N_0'/\Gamma _0'$; where
$N_0'\equiv \pi _1(X',x_0')\overset {\cong }\to N_0$ for
$x_0'=\Upsilon '(\hat x_0)$ and $\Gamma _0'=\Upsilon '_*\pi
_1(\widehat X,\hat x_0)$.

It remains to show that $X'$ admits a proper \holo mapping onto a
Riemann surface. But each automorphism $\sigma\in\Lambda '$ maps
fibers (of $\hat\Phi $) to fibers, because $\hat\Phi $ is constant
on every \con \cpt \anal set. Thus $\sigma $ descends to an
automorphism of $\widehat S$ and $\Lambda '$ acts properly
discontinuously on $\widehat S$ (in other words, the image of the
homomorphism $\Lambda '\to \text{Aut}\, (\widehat S)$ acts
properly discontinuously and the kernel is finite). We therefore
get a commutative diagram
\begin{center}\begin{picture}(230,70)
\put(80,9){$X'$} \put(93,13){\vector(1,0){45}}
\put(140,9){$S'\equiv \Lambda '\big\backslash \widehat S$}
\put(80,49){$\widehat X$} \put(110,55){$\widehat\Phi $}
\put(93,53){\vector(1,0){45}} \put(140,49){$\widehat S$}
\put(149,31){$\theta $} \put(144,47){\vector(0,-1){28}}
\put(70,31){$\Upsilon '$} \put(84,47){\vector(0,-1){28}}
\put(110,15){$\Phi $}
\end{picture}
\end{center}
where $\theta $ is a branched infinite covering map, $S'$ is a
Riemann surface, and $\Phi $ is a surjective proper \holo map with
\con fibers.
\end{pf*}

\begin{lem}\label{Covering to hyperbolic Riemann surface gives map to hyperbolic compact Riemann surface}
Let $(X,g)$ be a \con \cpt
K\"ahler manifold . Suppose some \con covering space $\Upsilon
:\widehat X\to X$ admits a surjective proper \holo mapping
$\hat\Phi :\widehat X\to \widehat S$ onto a Riemann surface
$\widehat S$ whose universal covering is the unit disk~$\Delta $.
Then some finite covering space of $X$ admits a surjective \holo
mapping onto a curve of genus $g\geq 2$.
\end{lem}
\begin{pf}
Clearly, we may assume that $\Upsilon :\widehat X\to X$ is an
infinite covering and, by Stein factorization, we may assume that
$\hat\Phi $ has \con fibers.  After passing to the appropriate
covering space, we may also assume that $\pi _1(F)$ surjects onto
$\pi _1(\widehat X)$ for every fiber $F$ of $\hat\Phi $ and hence
that $\widehat S=\Delta $ (in the above, we have used part~(d) of
Lemma~\ref{Orbifold fundamental group lemma} and the fact that,
for a surjective \holo map of Riemann surfaces $S^*\to \widehat
S$, we get a lifting $\widetilde S^*\to\Delta $, where $\widetilde
S^*$ is the universal covering of $S^*$, and hence $\widetilde
S^*=\Delta $ by Liouville's theorem). Furthermore, applying (the
proof of) Proposition~\ref{Kollar normalizer prop}, we get a
commutative diagram of \holo mappings
\begin{center}\begin{picture}(230,70)
\put(57,43){$\Upsilon $}
\put(65,23){$X$}\put(49,12){\vector(1,1){15}}\put(56,11){$\beta $}
\put(0,0){$\check\Lambda\big\backslash\widehat X=\check X$}
\put(55,3){\vector(1,0){80}} \put(140,0){$\check S=\check\Theta
\big\backslash \Delta $} \put(43,60){$\widehat X$}
\put(50,57){\vector(1,-2){13}} \put(90,66){$\hat\Phi $}
\put(55,63){\vector(1,0){80}} \put(140,60){$\widehat S=\Delta$}
\put(147,35){$\theta $} \put(144,58){\vector(0,-1){45}}
\put(35,35){$\alpha $} \put(45,58){\vector(0,-1){45}}
\put(90,5){$\check\Phi $}
\end{picture}
\end{center}
where $\beta :\check X\to X$ is the finite covering with $\beta
_*\pi _1(\check X)$ equal to the normalizer of $\Upsilon _*\pi
_1(\widehat X)$; $\alpha $ is the corresponding intermediate
covering map; $\check\Lambda =\pi _1(\check X)/\check\Gamma $ is
the quotient by the normal subgroup $\check\Gamma =\alpha _*\pi
_1(\widehat X)$; $\check\Theta $ is the image of $\check\Lambda $
under the homomorphism $\check\Lambda\to\text{Aut}\, (\Delta )$;
and $\theta $ is the corresponding branched covering map of the
quotient Riemann surface $\check S$. Since $\check X$ is \cpt
(and, therefore, $\check\Lambda $ is finitely generated),
Selberg's lemma provides a finite index torsion free normal
subgroup $\Theta '$ of $\check\Theta $ and hence a commutative
diagram of \holo mappings
\begin{center}\begin{picture}(230,100)
\put(60,59){$\Upsilon $} \put(125,52){\vector(-3,-2){68}}
\put(65,23){$X$}\put(49,12){\vector(1,1){15}}\put(50,22){$\beta $}
\put(0,0){$\check\Lambda\big\backslash\widehat X=\check X$}
\put(55,3){\vector(1,0){148}} \put(208,0){$\check S=\check\Theta
\big\backslash \Delta $} \put(43,85){$\widehat X$}
\put(50,82){\vector(1,-3){16}} \put(130,91){$\hat\Phi $}
\put(55,88){\vector(1,0){148}} \put(208,85){$\widehat S=\Delta$}
\put(240,50){$\theta $} \put(237,83){\vector(0,-1){70}}
\put(35,50){$\alpha $} \put(45,83){\vector(0,-1){70}}
\put(130,5){$\check\Phi $} \put(85,55){$\Lambda
'\big\backslash\widehat X=X'$} \put(90,52){\vector(-1,-1){18}}
\put(144,58){\vector(1,0){26}}\put(172,55){$S'=\Theta
'\big\backslash\Delta $}\put(152,61){$\Phi '$}
\put(180,53){\vector(1,-2){22}}\put(206,84){\vector(-1,-1){21}}
\put(53,83){\vector(2,-1){33}}\put(187,73){$\theta '$}
\end{picture}
\end{center}
where $\Lambda '\subset \check\Lambda $ is the inverse image of
$\Theta '$. Clearly, the map $\Phi '$ is surjective and the map
$\theta '$ is an unramified covering map, so the proof is
complete.
\end{pf}

\begin{pf*}{Proof of Theorem~\ref{Mapping of compact to curve theorem}}
By Theorem~\ref{Main Theorem from Introduction}, $\widehat X$
admits a proper \holo mapping $\hat\Phi$ with \con fibers onto a
(non\cptns ) Riemann surface $\widehat S$. For an ends
decomposition $\widehat X\setminus K=E_1\cup\cdots\cup E_m$ and a
generic fiber $F$ of $\hat\Phi$ contained in $E_1$, we get a
commutative diagram of \holo maps
\begin{center}\begin{picture}(335,100)
\put(159,52){\vector(-3,-2){68}}
\put(100,24){$X$}\put(84,13){\vector(1,1){15}}\put(84,23){$\Upsilon$}
\put(6,0){$\pi_1^{\text{orb}}(\hat\Phi )\big\backslash\widetilde
X=\widehat X$} \put(90,3){\vector(1,0){148}} \put(243,0){$\widehat
S=\pi _1^{\text{orb}}(\hat\Phi )\big\backslash \widetilde S $}
\put(78,85){$\widetilde X$} \put(85,82){\vector(1,-3){16}}
\put(165,91){$\tilde\Phi $} \put(90,88){\vector(1,0){148}}
\put(243,85){$\widetilde S=\Delta$ or $\C $} \put(275,50){$\theta
$} \put(272,83){\vector(0,-1){70}} \put(70,45){$\widehat\Upsilon$}
\put(80,83){\vector(0,-1){70}} \put(165,5){$\hat\Phi $}
\put(120,55){$\check\Lambda \big\backslash\widetilde X=\check X$}
\put(125,52){\vector(-1,-1){18}}
\put(179,58){\vector(1,0){26}}\put(207,55){$\check
S=\check\Lambda\big\backslash\widetilde S
$}\put(187,61){$\check\Phi$}
\put(215,53){\vector(1,-2){22}}\put(241,84){\vector(-1,-1){21}}
\put(88,83){\vector(2,-1){33}}
\put(105,77){$\alpha$}\put(125,20){$\beta$} \put(226,75){$\rho$}
\put(226,33){$\tau$}
\end{picture}
\end{center}
where $\widehat\Upsilon :\widetilde X\to\widehat X$ is the \con
Galois covering space with
$$
\widehat\Upsilon _*\pi _1(\widetilde X)=\Gamma \equiv \image{\pi
_1(F)\to\pi _1(\widehat X)}=\kernel{\pi _1(\widehat X)\to\pi
_1^{\text{orb}}(\hat\Phi )};
$$
$\beta :\check X\to \widehat X$ is the \con covering space with
$$
\beta _*\pi _1(\check X)= \check\Gamma\equiv\image{\pi
_1(E_1)\to\pi _1(\widehat X)}\supset\Gamma ;
$$
$\check\Lambda$ is the group given by
$$\check\Lambda=\check \Gamma
/\Gamma \subset \pi _1(\widehat X)/\Gamma \cong\pi
_1^{\text{orb}}(\hat\Phi );
$$
$\alpha :\widetilde X\to \check X$ is the \con Galois covering
space with
$$
\alpha _*\pi _1(\widetilde X)= \image{\pi _1(\check F)\to\pi
_1(\check X)}=(\beta _*)\inv (\Gamma )\cong \Gamma
$$
for a fiber $\check F$ of $\check\Phi$ which $\beta $ maps
isomorphically onto $F$; the maps $\theta $, $\rho $, and $\tau $
are branched covering maps (with branch locus mapping into the set
of critical values of $\hat\Phi $); and the maps $\tilde\Phi$ and
$\check\Phi$ (and $\hat\Phi$) are surjective proper \holo mappings
with \con fibers. Here, we identify $\pi _1^{\text{orb}}(\hat\Phi
)$ with the corresponding discrete subgroup of $\text{Aut}\,
(\widetilde X)$ or $\text{Aut}\, (\widetilde S)$, depending on the
context. According to Lemma~\ref{Covering to hyperbolic Riemann
surface gives map to hyperbolic compact Riemann surface}, it now
suffices to show that $\widetilde S=\Delta $ (in fact, this will
imply that one may get the desired finite covering $X'$  of $X$ by
applying Lemma~\ref{Covering to hyperbolic Riemann surface gives
map to hyperbolic compact Riemann surface} to $X$ and a suitable
covering $\widetilde X$ of the given covering $\widehat X$).
Observe that this is the case if some intermediate infinite
covering space $X^*$ between $\widetilde X$ and $X$ admits a
proper \holo mapping $\Phi ^*$ with \con fibers onto a Riemann
surface $S^*$ whose universal covering is $\Delta $. For the
fibers of $\Phi ^*$ lift to a union of fibers of $\tilde\Phi$
(since $\Phi ^*$ is constant on the image of each fiber of
$\tilde\Phi$ in $X^*$), so the covering map $\widetilde X\to X^*$
descends to a surjective \holo mapping $\widetilde S\to S^*$.
Lifting to a map $\widetilde S \to\Delta $, we see that
$\widetilde S\neq\C$ and hence $\widetilde S=\Delta $.

In particular, it suffices to consider the cases $\widehat
S=\C$~or~$\C ^*$. Since $e(\widehat X)=e(\widehat S)$ and $\tilde
e(\widehat X)\geq 3$, we may choose the ends decomposition
$\widehat X\setminus K=E_1\cup\cdots\cup E_m$ so that
$m=e(\widehat X)=1$~or~$2$ and $\check\Gamma\neq\pi _1(\widehat X)$. It
follows that $e(\check X)\geq 2$ and therefore, since
$e(\widetilde X)=e(\widetilde S)=1$, $\check\Lambda =\check\Gamma
/\Gamma$ is infinite (i.e. $\alpha :\widetilde X\to \check X$ is
an infinite covering). On the other hand, we have $[\pi
^{\text{orb}}_1(\hat\Phi):\check\Lambda
 ]=[\pi _1(\widehat X):\check\Gamma ]$. If this index is finite (i.e. $\beta: \check X\to \widehat X$ is a
finite covering), then $\tilde e(\check X)=\tilde e(\widehat
X)\geq 3$ and hence, in this case, we may replace $\widehat X$ by
$\check X$. Thus we may assume without loss of generality that
$\check\Lambda$ is of infinite index in $\pi
^{\text{orb}}_1(\hat\Phi)$ or $e(\widehat X)=2$ (i.e. $\widehat
X=\C ^*$).

We now consider the possible properties of $\pi
_1^{\text{orb}}(\hat\Phi )$. We first observe that, if $\pi
_1^{\text{orb}}(\hat\Phi )$ contains an infinite cyclic subgroup
of finite index, then every infinite subgroup of $\pi
_1^{\text{orb}}(\hat\Phi )$ is of finite index. In particular, we
get $e(\widehat X)=2$ and $\beta :\check X\to\widehat X$ is a
finite covering and, therefore, $e(\check X)\geq 3$. Therefore,
$\check S\neq \C$~or~$\C ^*$ and hence $\widetilde S=\Delta$ in
this case. In general, if $\widehat S=\C$, then $\pi
^{\text{orb}}_1(\hat\Phi)$ is either trivial or a free product of
a countable collection of nontrivial finite cyclic groups while,
if $\widehat S=\C^*$, then $\pi ^{\text{orb}}_1(\hat\Phi)$ is
either $\Z$ or the free product of $\Z$ and a countable collection
of nontrivial finite cyclic groups. The group $\Z _2*\Z _2$
contains the finite index infinite cyclic subgroup $< a_1a_2>$,
where $a_1$ and $a_2$ are the generators for the first and second
copy of $\Z _2$, respectively. Thus, by the above remarks, we need
only consider cases in which the free product representation for
$\pi ^{\text{orb}}_1(\hat\Phi)$ has at least three nontrivial
factors or at least two nontrivial factors with one of order $>2$.
But in these cases, $\pi ^{\text{orb}}_1(\hat\Phi)$ contains a
non-Abelian free group and is therefore Fuchsian.
\end{pf*}

\begin{rmk}
Lemma~\ref{Two ends and infinitely generated lemma} and
Theorem~\ref{Mapping of compact to curve theorem} together give
part~(b) of Theorem~\ref{Gromov-Schoen Theorem}. A proof giving
both parts (a) (Gromov and Schoen) and (b) simultaneously appears
in the next section.
\end{rmk}

\section{Amalgamations and mappings to Riemann surfaces}\label{Amalgamations and mappings to Riemann
surfaces section}

Theorem~\ref{Mapping of compact to curve theorem} together with
standard facts from geometric group theory (see
Proposition~\ref{Action on tree gives 3 filtered ends prop} below)
give Theorem~\ref{Gromov-Schoen Theorem}, which is equivalent to
the following:
\begin{thm}\label{Action on tree version of Gromov-Schoen theorem}
Suppose $X$ is a \con \cpt K\"ahler manifold whose fundamental
group $\Lambda =\pi _1(X)$ induces a minimal action without
inversion on a (simplicial) tree~$T$ which is not a line or a
point. Then some finite covering space $X'\to X$ admits a
surjective \holo mapping onto a curve of genus $g\geq 2$.
\end{thm}
\begin{rmks}
1. Equivalently, if $X$ is a \con \cpt K\"ahler manifold whose
fundamental group $\Lambda =\pi _1(X)$ is the fundamental group of
a minimal reduced graph of groups for which the universal covering
tree is not a line or a point, then some finite covering space
admits a surjective \holo mapping onto a curve of genus $\geq 2$.

\noindent 2. We will only consider group actions on {\it
simplicial} trees.
\end{rmks}

A brief discussion of the required facts from Bass-Serre Theory
\cite{Se} will be provided for the convenience of the reader. For
a graph $Y$, we will denote the set of {\it vertices} by $\Vert
(Y)$ and the set of {\it edges} by $\Edge (Y)$, the {\it origin}
and {\it terminus} maps by
$$
\alpha :\Edge (Y)\to\Vert (Y)\qquad\text{and}\qquad \omega :\Edge
(Y)\to\Vert (Y),
$$
respectively, and the {\it edge inversion} by $e\mapsto\bar e$. We
often identify each edge $e\in\Edge (Y)$ with the corresponding
map $e:[0,1]\to Y$. In a metric space $Z$, for $r>0$, we denote
the ball of radius $r$ centered at a point $z\in Z$ by $B_Z(z;r)$
and the $r$-\nbd of a subset $A\subset Z$ by $N_Z(A;r)$. In
particular, for a \con graph~$Y$, a vertex $v\in\Vert (Y)$, and a
number $r\in (0,1]$, we have
$$
B_Y(v;r)=\bigcup _{e\in\Edge (Y),\, \alpha (e)=v} e([0,r)).
$$

Suppose $\Lambda $ is a group which acts without inversion on a
tree~$T$ (i.e. $\lambda e\neq \bar e$ for every edge $e\in\Edge
(T)$). The action is called {\it minimal} if $T$ contains no
proper subtree that is invariant under the action. In general, if
$\Lambda$ is finitely generated, then one can form a minimal
$\Lambda$-invariant subtree $T_{\min}$. For we may take $T_{\min}$
to be a vertex~$v$ of $T$ fixed by~$\Lambda$ if such a vertex
exists (i.e. if the action is elliptic).  If no vertex is fixed,
then we may take $T_{\min}$ to be the intersection of all
$\Lambda$-invariant subtrees.

Observe also that, if the action of $\Lambda$ on $T$ is minimal
and $T$ is not a single point (i.e. no vertex is fixed), then the
valence of every vertex is at least~$2$ (i.e. there are no ``dead
ends''). For the collection of edges having an endpoint of
valence~$1$ in $T$ is invariant and hence, after removing such
edges, one gets an invariant subtree which, by minimality, must be
equal to $T$. Thus no such edges can exist. Equivalently, each
edge of $T$ is contained in a line.

Theorem~\ref{Action on tree version of Gromov-Schoen theorem} is
an immediate consequence of Theorem~\ref{Mapping of compact to
curve theorem} and the following fact:
\begin{prop}\label{Action on tree gives 3 filtered ends prop}
Suppose $M$ is a \con \cpt \cinf manifold whose fundamental group
$\Lambda =\pi _1(M)$ induces a minimal action without inversion on
a tree~$T$ which has a vertex $v_0$ of valence at least~$3$. Then
the \con covering $\widehat M\to M$ for which $\text{{\rm im}}\,
[\pi _1(\widehat M)\to\pi _1(M)]$ is equal to the isotropy
subgroup $\Lambda _{v_0}$ at $v_0$ is an infinite covering (i.e.
$[\Lambda :\Lambda _{v_0}]=\infty $) and $\tilde e(\widehat M)\geq
3$.
\end{prop}
\begin{rmk}
For $\Lambda =\pi _1(X)=\Gamma _1 * _\Gamma \Gamma _2$ as in the
theorem of Gromov and Schoen (Theorem~\ref{Gromov-Schoen
Theorem}), $\Lambda$ induces a minimal action on a tree $T$ with
fundamental domain $e\in\Edge(T)$ \st $\Gamma$ is the stabilizer
of $e$, $\Gamma_1=\Lambda _{\alpha (e)}$, and $\Gamma _2=\Lambda
_{\omega (e)}$, and the index of $\Gamma$ in $\Gamma _1$ and
$\Gamma _2$ is equal to the valence of $\alpha (e)$ and $\omega
(e)$, respectively.
\end{rmk}

For the proof of the proposition, we first consider the following
standard fact:

\begin{lem}\label{Action on tree isotropy subgroup subgraph lemma}
Given a minimal co\cpt action without inversion of a finitely
generated group~$\Lambda$ on a tree~$T$ which is not a point, let
$\tilde\theta\colon T\to Y$ be the quotient map to the finite
quotient graph $Y=\Lambda\big\backslash T$, let $v_0$ be a vertex
in $T$, let $\theta\colon T\to\widehat
T\equiv\Lambda_{v_0}\big\backslash T$ be the quotient by the
isotropy subgroup $\Lambda_{v_0}$, let $\hat\theta\colon\widehat
T\to Y$ be the induced map, let $y_0=\tilde\theta (v_0)$, let
$\hat v_0=\theta (v_0)$, and let $T_0$ be the subtree of $T$ with
$$
\Edge (T_0)=\setof{e\in\Edge (T)}{\alpha (e)=v_0\text{ or }\omega
(e)=v_0} $$ and
$$
\Vert (T_0)=\{ v_0\}\cup\setof{\omega (e)}{e\in\Edge (T)\text{ and
}\alpha (e)=v_0}.
$$
Then we have the following:
\begin{enumerate}
\item[(a)] The image $Y_0\equiv\tilde\theta (T_0)$ (i.e. $\Vert
(Y_0)=\tilde\theta (\Vert (T_0))$, $\Edge (Y_0)=\tilde\theta
(\Edge (T_0))$) is the (finite) subgraph of $Y$ with
$$
\Edge (Y_0)=\setof{e\in\Edge (Y)}{\alpha (e)=y_0\text{ or }\omega
(e)=y_0} $$ and
$$
\Vert (Y_0)=\{ y_0\}\cup\setof{\omega (e)}{e\in\Edge (Y)\text{ and
}\alpha (e)=y_0}.
$$

\item[(b)] The graph $\widehat T$ is a tree.

\item[(c)] $\Lambda _{v_0}$ acts on $T_0$ and on $T\setminus T_0$,
$\widehat T_0\equiv\Lambda _{v_0}\big\backslash T_0$ is the finite
subtree of $\widehat T$ with
$$
\Edge (\widehat T_0)=\setof{e\in\Edge (\widehat T)}{\alpha
(e)=\hat v_0\text{ or }\omega (e)=\hat v_0} $$ and
\[
\Vert (\widehat T_0)=\{ \hat v_0\}\cup\setof{\omega (e)}{e\in\Edge
(\widehat T)\text{ and }\alpha (e)=\hat v_0}, \] \(T_0=\theta\inv
(\widehat T_0)\), and there exists a finite subtree $T_1$ of $T_0$
which $\theta $ maps isomorphically onto $\widehat T_0$.

\item[(d)] If $R:[0,\infty )\to T$ is a ray in $T$ with vertex
$R(0)=v_0$, then $\theta$ maps $R$ isomorphically onto a ray
$\theta (R)$ in $\widehat T$.

\end{enumerate}

\end{lem}
\begin{pf}
For the proof of (a), we observe that if $e\in\Edge (Y_0)$ with
$\alpha (e)=y_0=\tilde\theta (v_0)$, then there is an edge
$f\in\Edge (T)$ with $\tilde\theta (f)=e$. Since $\tilde\theta
(\alpha (f))=y_0$, there is an element $\lambda\in\Lambda$ with
$\alpha (\lambda\cdot f)=\lambda\cdot\alpha (f)=v_0$. Thus
$\lambda\cdot f\in\Edge (T_0)$ and $\tilde\theta (\lambda\cdot
f)=e$.

For the proof of (b), suppose $\widehat T$ contains a circuit.
Then (equivalently) $\widehat T$ contains a loop $\hat\gamma=\hat
e_1*\hat e_2*\cdots *\hat e_k$; where $\hat e_1,\hat e_2,\dots
,\hat e_k$ are edges \st $\alpha (\hat e_1)=\omega (\hat
e_k)=\theta (v_0)=\hat v_0$ and, for some $i\in\{\, 1,\dots
,k\,\}$, we have $\hat e_i\neq\hat e_j$~and~$\hat
e_i\neq\overline{\hat e_j}$ for $j=1,\dots ,\hat i,\dots k$. We
may then lift $\hat\gamma$ to a path $\gamma =e_1*e_2*\cdots *e_k$
with $\alpha (e_1)=v_0$.  Since $\omega (e_k)\in\Lambda
_{v_0}\cdot v_0=\{ v_0\}$, $\hat\gamma$ must be a loop. But we
also have $e_i\neq e_j$~and~$e_i\neq\bar e_j$ for $j=1,\dots ,\hat
i,\dots k$, which is impossible since $T$ is a tree. Thus
$\widehat T$ is a tree.

For (c), it is clear that $\Lambda _{v_0}$ acts on $T_0$ and on
$T\setminus T_0$ and the proof of (a) shows that $\widehat T_0$ is
the subtree of $\widehat T$ as described. We may form the finite
subtree $T_1$ of $T_0$ as follows. For each edge $e\in \Edge
(Y_0)$ with $\alpha (e)=y_0$, we may choose an edge $f_e\in\Edge
(T_0)$ with $\alpha (f_e)=v_0$ and $\tilde\theta (f_e)=e$. If $e$
is a loop edge, then we may choose an element $\lambda
_e\in\Lambda$ \st $\omega (\lambda _ef_e)=\lambda _e\cdot\omega
(f_e)=v_0$. Hence $\alpha (\overline{\lambda _ef_e})=v_0$,
$\tilde\theta (\overline{\lambda _ef_e})=\tilde\theta (\lambda _e
\overline{f_e})=\tilde\theta (\overline{f_e})=\bar e$, and
$f_e=\overline{\lambda\inv _e(\overline{\lambda _ef_e})}$. Thus we
may choose $f_e$, $f_{\bar e}$, $\lambda _e$, and $\lambda _{\bar
e}$ so that $f_{\bar e}=\overline{\lambda _ef_e}$ and $\lambda
_{\bar e}=\lambda _e\inv $ for each loop edge $e\in\Edge (Y_0)$.
Note also that $f_{\bar e}\neq f_e$ since $\Lambda $ acts without
inversion. We now define the subtree $T_1$ by
$$
\Edge (T_1)=\setof{f_e}{e\in\Edge (Y_0)\text{ and }\alpha (e)=y_0}
\cup \setof{\overline {f_e}}{e\in\Edge (Y_0)\text{ and }\alpha
(e)=y_0}
$$
and
$$
\Vert (T_1)=\{ v_0\}\cup\setof{\omega (f_e)}{e\in\Edge (Y_0)\text{
and }\alpha (e)=y_0}.
$$
We have $\Lambda _{v_0}\cdot T_1=T_0$. For, if $f\in\Edge (T_0)$
with $\alpha (f)=v_0$, then the edge $e=\tilde\theta (f)$ is an
edge in $Y_0$ with initial point~$y_0$ and hence $\tilde\theta
(f)=e=\tilde\theta (f_e)$. Thus $f=\lambda f_e$ for some
$\lambda\in\Lambda$. In particular, $\lambda v_0=\lambda\alpha
(f_e)=\alpha (f)=v_0$, so $\lambda\in\Lambda _{v_0}$. Thus
$T_0=\Lambda _{v_0}\cdot T_1$. On the other hand, if $f,f'\in\Edge
(T_1)$ with $\lambda f=f'$ for some $\lambda\in\Lambda _{v_0}$,
then $f=f'$. For we may assume without loss of generality that
$\alpha (f)=\alpha (f')=v_0$ (otherwise, we replace the pair by
$\bar f, \bar f'$). We then have $\tilde\theta (f)=\tilde\theta
(f')=e\in\Edge (Y_0)$ and hence $f=f'=f_e$ by the construction of
$T_1$. Thus $\theta$ maps $T_1$ isomorphically onto $\widehat
T_0$. In particular, $\widehat T_0=\Lambda _{v_0}\big\backslash
T_0$ is a finite subtree of the tree $\widehat T=\Lambda
_{v_0}\big\backslash T$.

Finally, for (d), observe that if $v=R(j)$ and $w=R(k)$ are
vertices in $R$ ($j,k\in\N $) and $\lambda\in\Lambda _{v_0}$ with
$\lambda v=w$, then $\lambda \cdot [v_0,v]=[v_0,w]\subset R$.
Since we then have $j=l([v_0,v])=l([v_0,w])=k$, we get $v=w$.

\end{pf}

\begin{pf*}{Proof of Proposition~\ref{Action on tree gives 3 filtered ends prop}}
Let $\widetilde\Upsilon :\widetilde M\to M$ be the universal
covering of $M$ and let $B=\widetilde M\times _\Lambda T\to M$ be
the associated bundle. Here, we will consider the left action of
$\Lambda$ on $\widetilde M$ so that
$$
B=\widetilde M\times T/[(x,t)\sim (\lambda\cdot x,\lambda\cdot
t)\quad\forall\, \lambda\in\Lambda ]
$$
(in terms of the right action, $(\lambda\cdot x,\lambda\cdot
t)=(x\cdot \lambda\inv,\lambda\cdot t)\sim (x,t)$). Since the
fiber $T$ is contractible, $B\to M$ admits a continuous
section~$\sigma$. Thus we get a \cont equivariant map
$\widetilde\Psi :\widetilde M\to T$ given by
$$
\widetilde\Psi (x)=t \iff [(x,t)]=\sigma (\widetilde\Upsilon (x)).
$$
Furthermore, the minimality of the action implies that
$\widetilde\Psi$ is surjective. To see this, we let $T'$ be the
subgraph of $T$ for which $\Edge (T')$ is the set of edges
$e\in\Edge (T)$ which lie entirely (including the endpoints) in
$\widetilde\Psi (\widetilde M)$ and $\Vert (T')=\Vert (T)\cap
\widetilde\Psi (\widetilde M)$. Then $T'$ is $\Lambda$-invariant
because $\widetilde\Psi (\widetilde M)$ is $\Lambda$-invariant.
Furthermore, $T'$ is \con and, therefore, $T'$ is a subtree. For
if $p$ and $q$ are distinct points in $T'$, then we may form a
geodesic $\eta $ in $T$ from $p$ to $q$ which is contained in the
path \con set $\widetilde\Psi (\widetilde M)$. Clearly, each
vertex in $T$ which $\eta$ meets will lie in $T'$. If $e$ is an
edge whose interior meets $\eta$, then either $p$ or $q$ is in the
interior or $e$ is a segment of $\eta$. In either case, we get
$e\in\Edge (T')$. Thus $\eta$ is contained in $T'$ and hence $T'$
is a subtree. Therefore, by minimality, we have $T'=T$ and hence
$\widetilde\Psi (\widetilde M)=T$. Thus we get a commutative
diagram of surjective \cont mappings
\begin{center}\begin{picture}(290,100)
\put(119,44){\vector(-3,-2){56}}
\put(6,0){$\Lambda\big\backslash\widetilde M=M$}
\put(65,3){\vector(1,0){168}}
\put(236,0){$Y\equiv\Lambda\big\backslash T$}
\put(53,85){$\widetilde M$} \put(140,91){$\tilde\Psi $}
\put(65,88){\vector(1,0){168}} \put(236,85){$T$}
\put(243,45){$\tilde\theta $} \put(239,83){\vector(0,-1){70}}
\put(45,45){$\widetilde\Upsilon$} \put(55,83){\vector(0,-1){70}}
\put(140,5){$\Psi $} \put(69,45){$\Lambda
_{v_0}\big\backslash\widetilde M\equiv\widehat M$}
\put(137,48){\vector(1,0){26}}\put(166,45){$\widehat
T\equiv\Lambda _{v_0}\big\backslash
T$}\put(145,50){$\widehat\Psi$} \put(174,42){\vector(2,-1){58}}
\put(233,83){\vector(-2,-1){56}} \put(63,83){\vector(2,-1){56}}
\put(90,71){$\Upsilon$}\put(90,11){$\widehat\Upsilon$}
\put(201,71){$\theta$} \put(201,11){$\hat\theta$}
\end{picture}
\end{center}
In particular, the quotient graph $Y$ is \cptns ; that is, finite.
Thus we may form $y_0=\tilde\theta (v_0)$, $\hat v_0=\theta
(v_0)$, $T\supset T_0\supset T_1$, $\widehat
T_0=\Lambda_{v_0}\big\backslash T_0\subset\widehat T$, and
$Y_0=\tilde\theta (T_0)=\hat\theta (\widehat T_0)$ as in
Lemma~\ref{Action on tree isotropy subgroup subgraph lemma}.

Let
$$
C_0\equiv\overline{B_T(v_0;1/4)}=\bigcup _{e\in\Edge (T),\, \alpha
(e)=v_0} e([0,1/4])
$$
and
$$
C_1\equiv\overline{B_{T_1}(v_0;1/4)}=\bigcup _{e\in\Edge (T_1),\,
\alpha (e)=v_0} e([0,1/4]).
$$
We then have
$$
C_0=\bigcup _{e\in\Edge (T_0),\, \alpha (e)=v_0}
e([0,1/4])=\overline{B_{T_0}(v_0;1/4)}=\Lambda _{v_0}\cdot C_1.
$$
Moreover, $\theta$ maps $C_1$ isomorphically onto the \con \cpt
set
$$
\widehat C_0\equiv\theta (C_0)=\overline{B_{\widehat T}(\hat
v_0;1/4)}=\overline{B_{\widehat T_0}(\hat v_0;1/4)}.
$$

We will pull back components of $\widehat T\setminus\widehat C_0$
to get ends in $\widehat M$. For this, we first observe that the
set $K_0\equiv\widehat\Psi\inv (\widehat C_0)$ is \cptns . For,
given a point
$$
x_0\in\widehat\Upsilon (K_0)=\Psi\inv (\hat\theta (\widehat
C_0))=\Psi\inv (\tilde\theta (C_0))=\Psi\inv (\tilde\theta (C_1)),
$$
we may choose \con \nbds $B$ and $B'$ \st $x_0\in B\Subset B'$,
$B'$ is contractible, and $\diam \Psi (B')<1/4$. Suppose $\widehat
B_1$ and $\widehat B_2$ are two \comps of $\widehat\Upsilon\inv
(B)$ which meet $K_0$. Then, for $i=1,2$, we may choose a \comp
$\widetilde B_i$ of $\Upsilon\inv (\widehat B_i)$ so that
$\widetilde \Psi (\widetilde B_i)$ meets $C_1$. We then have
$$
\tilde\theta (\widetilde \Psi (\widetilde B_i))=\Psi (B)\subset
N_Y(\tilde\theta (C_1);1/4)\subset B_Y(y_0;1/2).
$$
Thus $\widetilde\Psi (\widetilde B_i)$ is a \con subset of
$N_T(\Lambda\cdot v_0;1/2)$ which meets $C_1\subset B_T(v_0;1/2)$
and hence $\widetilde\Psi (\widetilde B_i)\subset B_T(v_0;1/2)$.
On the other hand, for some $\lambda\in\Lambda$, $\widetilde
B_2=\lambda\cdot \widetilde B_1$ and so $B_T(v_0;1/2)\cap
B_T(\lambda v_0;1/2)\neq\emptyset$. Thus $\lambda\in \Lambda
_{v_0}$ and hence $\widehat B_1=\Upsilon (\widetilde B_1)=\Upsilon
(\widetilde B_2)=\widehat B_2$. Therefore, a unique \comp
$\widehat B_1\Subset \widehat M$ of $\widehat\Upsilon\inv (B)$
meets $K_0$. Covering the \cpt set $\widehat\Upsilon
(K_0)=\Psi\inv (\tilde\theta (C_1))$ by finitely many such
sets~$B$, we see that $K_0$ is \cpt (in fact, we have shown that
$\widehat\Upsilon$ maps $K_0$ homeomorphically onto
$\widehat\Upsilon (K_0)$).

We may now choose an ends decomposition $\widehat M\setminus
K=E_1\cup\cdots\cup E_m$ in which $K$ is a \con \cpt set with
$K_0\subset K$ and we may set $C=\widehat\Psi (K)\supset\widehat
C_0$. If $R:[0,\infty )\to\widehat T$ is a ray with vertex
$R(0)=\hat v_0$, then, since $C$ is a \cpt \con subset of the tree
$\widehat T$ and $\widehat C_0\subset C$, we have $C\cap
R=R([0,b])$ for some $b\geq 1/4$. Similarly, for each $j=1,\dots
,m$, $\widehat\Psi (E_j)\cap R$ is a singleton or a (possibly
unbounded) interval whose closure meets $C\cap R$. But then
$$
R((b,\infty ))=R\setminus C=R\setminus\widehat\Psi
(K)\subset\widehat\Psi (\widehat M\setminus K)=\widehat\Psi
(E_1)\cup\cdots\cup\widehat\Psi (E_m)
$$
(since $\widehat\Psi$ is surjective) and hence $R\setminus
C\subset\widehat\Psi (E_j)$ for some~$j$. In other words, for each
ray $R$ in $\widehat T$ with vertex $\hat v_0$, there is an index
$j\in\{\, 1,\dots ,m\,\}$ \st the unbounded interval $R\setminus
C$ is contained in $\widehat\Psi (E_j)$. Note also that, if $S$ is
the \comp of $\widehat T\setminus\{ v_0\}$ containing $R((0,\infty
))$, then $\widehat\Psi (E_j)\subset S$ because
$$
\widehat\Psi (E_j)\subset\widehat\Psi (\widehat M\setminus
K)\subset\widehat\Psi (\widehat M\setminus K_0)=\widehat
T\setminus \widehat C_0\subset\widehat T\setminus \{ v_0\} $$ and
$\widehat\Psi (E_j)$ meets~$R((0,\infty ))$.

Now, by hypothesis, there exist distinct edges
$u_1,u_2,u_3\in\Edge (T)$ and rays $R_1,R_2,R_3$ with vertex $v_0$
\stns , for $\nu =1,2,3$, $\alpha (u_\nu )=v_0$ and $u_\nu =R_\nu
\restrict{[0,1]}$ (see the remarks preceding the statement of the
proposition). We may choose (not necessarily distinct) indices
$i_1,i_2,i_3\in\{\, 1,\dots ,m\,\}$ \stns , for $\nu =1,2,3$,
$$
\widehat R_\nu\setminus C\subset\widehat\Psi (E_{i_\nu})\subset
\widehat S_\nu ,
$$
where $\widehat R_\nu$ is the ray given by $\widehat R_\nu=\theta
(R_\nu )$ and $\widehat S_\nu$ is the \comp of $\widehat
T\setminus \{ \hat v_0\}$ containing $\widehat R_\nu \setminus \{
\hat v_0\}$. We then have, for some $b_\nu\geq 1/4$, $\widehat
R_\nu\cap C=\widehat R_\nu\restrict{[0,b_\nu ]}$, $\widehat
R_\nu\setminus C=\widehat R_\nu\restrict{(b_\nu ,\infty )}$,
$R_\nu\cap \theta\inv (C)=R_\nu\restrict{[0,b_\nu ]}$, and
$R_\nu\setminus\theta\inv (C)=R_\nu\restrict{(b_\nu ,\infty )}$.
We may also choose a \comp $F_\nu$ of $\Upsilon\inv (E_{i_\nu})$
\st $\widetilde\Psi (F_\nu )$ meets $R_\nu\setminus\theta\inv
(C)$. We have
$$
\theta (\widetilde\Psi (F_\nu ))=\widehat\Psi (\Upsilon (F_\nu
))=\widehat\Psi (E_{i_\nu})\subset\widehat S_\nu ,
$$
so $\widetilde\Psi (F_\nu )$ is a \con subset of $\theta\inv
(\widehat S_\nu )$ meeting $R_\nu\setminus\{ v_0\}$. Thus
$\widetilde\Psi (F_\nu )$ is contained in the \comp $S_\nu$ of
$T\setminus\{ v_0\}$ containing $R_\nu\setminus\{ v_0\}$ (note
that $\theta\inv (\hat v_0)=\{ v_0\}$). But
$S_1$,~$S_2$,~and~$S_3$ are disjoint because the sets $e_\nu
((0,1])\subset S_\nu$ for $\nu =1,2,3$ are in different \comps of
$T\setminus\{ v_0\}$. Therefore, $F_1$,~$F_2$,~and~$F_3$ are
disjoint and it follows that $\tilde e(\widehat M)\geq 3$.
\end{pf*}

We now complete the proof of Theorem~\ref{HNN and Thompson not
Kahler theorem from intro}.  Suppose $\Lambda$ is a properly
ascending HNN extension with base group $\Gamma$ and stable
letter~$\tau$. In other words, for some isomorphism $\vphi$ of
$\Gamma$ onto a proper subgroup of $\Gamma$, we have
\[
\Lambda=\langle\Gamma,\tau;\, \vphi(\gamma)=\tau\inv\gamma\tau
\text{ for }\gamma\in\Gamma\rangle.
\]
\begin{rmks} 1. We have $\Gamma\hookrightarrow\Lambda$ and
\[
\cdots\varsubsetneqq\tau^{-2}\Gamma\tau^2\varsubsetneqq\tau\inv\Gamma\tau
\varsubsetneqq\Gamma\varsubsetneqq\tau\Gamma\tau\inv\varsubsetneqq\tau^2\Gamma\tau^{-2}
\varsubsetneqq\cdots.
\]

\noindent 2. The group
$\hat\Lambda=\bigcup_{m\in\Z}\tau^m\Gamma\tau^{-m}$ is an
infinitely generated normal subgroup of $\Lambda$. The quotient
group $\Lambda\big\slash\hat\Lambda$ is infinite cyclic with
generator $\tau\hat\Lambda$.

\noindent 3. Clearly, if $\Gamma$ is finitely generated, then
$\Lambda$ is finitely generated. According to a theorem of Bieri
and Strebel \cite{Bieri-Strebel-HNN extension almost finitely
presented}, if $\Lambda$ is finitely presented, then one express
$\Lambda$ as an HNN extension with finitely generated base group.
However, it may be impossible to express $\Lambda$ as a {\it
properly ascending} HNN extension with {\it finitely generated}
base group. For example, if $\Lambda=\pi_1(S)=\langle
\alpha_1,\beta_1,\alpha_2,\beta_2;
[\alpha_1,\beta_1]\cdot[\alpha_2,\beta_2]=1\rangle$ for a \cpt
Riemann surface $S$ of genus~$2$; $\Gamma_0$ is the subgroup
generated by $\alpha_1,\beta_1,\alpha_2$, $\tau=\beta_2$;
$\Theta_1=\langle\alpha_2\rangle$;
$\Theta_2=\langle\alpha_1\beta_1\alpha_1\inv\beta_1\inv\alpha_2\rangle$;
and $\vphi\colon\Theta_1\overset{\cong}\to\Theta_2$ is the
isomorphism
$\alpha_2^m\mapsto(\alpha_1\beta_1\alpha_1\inv\beta_1\inv\alpha_2)^m$,
then
$\Lambda=\langle\tau,\Gamma_0;\tau\lambda\tau\inv=\vphi(\lambda)\,\forall\lambda\in\Theta_1\rangle$
is an HNN extension with (finitely generated) base group
$\Gamma_0$ and stable letter $\tau$. Moreover, for $\Gamma$ the
infinitely generated subgroup of $\Lambda$ generated by elements
$\tau^{-m}\gamma\tau^m$ for $\gamma\in\Gamma_0$ and $m\in\Z_{\geq
0}$, $\Lambda$ is a properly ascending HNN extension with base
group $\Gamma$ and stable letter~$\tau$. However, $\Lambda$ cannot
be expressed as a {\it properly ascending} HNN extension with {\it
finitely generated base group} (for example, by part~(ii) of
Theorem~\ref{HNN and Thompson not Kahler theorem from intro}).

\noindent 4. $\Lambda$ induces a minimal left action without
inversion on a tree $T$ \st the quotient graph
$\Lambda\big\backslash T$ is a single loop edge. Under the
embedding $\Gamma\hookrightarrow\Lambda$, we may identify the base
group $\Gamma$ with the isotropy subgroup $\Lambda_{v_0}$ for some
vertex $v_0$. For some unique edge~$e_0$, we have
$v_0=\alpha(e_0)$ and $v_1\equiv\omega(e_0)=\tau\cdot v_0$.  The
edges $\tau^me_0$, $m\in\Z$, form a line $l\colon\R\to T$, the
axis for~$\tau$, with $l\restrict{[m,m+1]}=\tau^m\cdot e_0$ for
each $m\in\Z$. For each $m\in\Z$, we set $v_m=\tau^m\cdot
v_0=l(m)$. We then have
$\Lambda_{v_m}=\tau^m\Lambda_{v_0}\tau^{-m}$
$\Lambda_{v_m}\supset\Lambda_{v_0}$ if $m\geq 0$ (while
$\Lambda_{v_m}\subset\Lambda_{v_0}$ if $m\leq 0$), so
$\Lambda_{v_0}$ must fix each point in the ray
$l\restrict{[0,\infty)}$. In particular, $\Lambda_{v_0}$ fixes the
edge $e_0$ and $\Lambda_{v_0}$ acts transitively on the remaining
edges with initial vertex $v_0$. For if $f$ is an edge not equal
to $e_0$ or $\tau\inv\cdot\bar e_0$ with $\alpha(f)=v_0$, then
there exists a $\lambda\in\Lambda$ with $\lambda\tau\inv\cdot\bar
e_0=f$ or $\bar f$. If the former, then $\lambda\in\Lambda_{v_0}$.
If the latter, then $f=\lambda\tau\inv\cdot e_0$ with
$\lambda\tau\inv\in\Lambda_{v_0}$, contradicting the above.

\noindent 5. Let $m\in\Z_{\geq 0}$ and let $D$ be the end in $T$
which is the component of $T\sm\set{v_{m+1}}$ containing~$v_0$.
Then $\Lambda_{v_m}=\tau^m\Gamma\tau^{-m}$ is precisely the set of
elements $\lambda\in\hat\Lambda$ with $\lambda\cdot v_0\in D$. For
each element of $\Lambda_{v_{m+1}}$ maps $D$ onto a component of
$T\sm\set{v_{m+1}}$. If
$\lambda\in\Lambda_{v_m}\subset\Lambda_{v_{m+1}}$, then
$\lambda\cdot\tau^m\cdot e_0=\tau^m\cdot e_0$, and hence
$\lambda\cdot D=D$. Conversely, if
$\lambda\in\hat\Lambda\sm\Lambda_{v_m}$, then we have
$\lambda\in\Lambda_{v_{k+1}}$ for some minimal $k\geq m$. In
particular, since $\tau^k\cdot e_0$ has endpoints $v_k$ and
$v_{k+1}$, $\lambda$ must map the component $F$ of
$T\sm\set{v_{k+1}}$ containing $\tau^k\cdot e_0$ onto a different
component. Hence, since $e_0\subset D\subset F$, $D$ and
$\lambda\cdot D$ must be disjoint.
\end{rmks}

\begin{lem}\label{HNN covering finitely gen lem}
Let $M$ be a \con \cpt \cinf manifold whose fundamental group
$\Lambda=\pi_1(M)$ is a properly ascending HNN extension with
finitely generated base group $\Gamma$ and stable letter $\tau$,
let $\hat\Lambda=\bigcup_{m\in\Z}\tau^m\Gamma\tau^{-m}$, and let
$\what\Upsilon\colon\what M\to M$ be a \con (Galois) covering
space with $\what\Upsilon_*\pi_1(\what M)=\hat\Lambda$. Then there
is an end $E$ for $\what M$ \st $\image{\pi_1(E)\to\pi_1(\what
M)}$ is finitely generated. In fact, for any \cinf \rel \cpt
domain $\Omega'$ in $\what M$ containing $\partial E$, the image
of the fundamental group of the end $E'=E\cup\Omega'$ in
$\pi_1(\what M)$ will be finitely generated.
\end{lem}
\begin{pf}
As in the proof of Proposition~\ref{Action on tree gives 3
filtered ends prop}, for $\wtil\Upsilon\colon\wtil M\to M$ the
universal covering and $T$ the universal covering tree as above,
we get a surjective \cont equivariant map $\wtil\Psi\colon\wtil
M\to T$ and a commutative diagram
\begin{center}\begin{picture}(290,100)
\put(119,44){\vector(-3,-2){56}}
\put(6,0){$\Lambda\big\backslash\widetilde M=M$}
\put(65,3){\vector(1,0){168}}
\put(236,0){$Y\equiv\Lambda\big\backslash T$}
\put(53,85){$\widetilde M$}
\put(140,91){$\tilde\Psi $}
\put(65,88){\vector(1,0){168}} \put(236,85){$T$}
\put(243,45){$\tilde\theta $} \put(239,83){\vector(0,-1){70}}
\put(45,45){$\widetilde\Upsilon$} \put(55,83){\vector(0,-1){70}}
\put(140,5){$\Psi $}
\put(69,45){$\hat\Lambda\big\backslash\widetilde M\equiv\widehat
M$} \put(133,48){\vector(1,0){26}}\put(166,45){$\widehat
T\equiv\hat\Lambda\big\backslash T$}\put(140,50){$\widehat\Psi$}
\put(174,42){\vector(2,-1){58}} \put(233,83){\vector(-2,-1){56}}
\put(63,83){\vector(2,-1){51}}
\put(90,71){$\Upsilon$}\put(90,11){$\widehat\Upsilon$}
\put(201,71){$\theta$} \put(201,11){$\hat\theta$}
\end{picture}
\end{center}
The quotient $\what T$ is a line onto which $\theta$ maps the line
$l$ isomorphically. We set $\hat l=\theta\circ l$. Furthermore,
the map $\what\Psi$ is proper (since
$\Lambda\big\slash\hat\Lambda\cong\Z$ acts properly
discontinuously on $\what M$ and $\what T$) and surjective and
$e(\what M)=2$. For some point $\til x_0\in\til M$ with
$\wtil\Psi(\til x_0)=v_0$ and for $\hat x_0=\Upsilon(\til x_0)$,
$\hat v_0=\theta(v_0)=\what\Psi(\hat x_0)$,
$x_0=\wtil\Upsilon(\til x_0)=\what\Upsilon(\hat x_0)$, and
$y_0=\til\theta (v_0)=\hat\theta(\hat v_0)=\Psi(x_0)$, we may
identify $\Lambda=\pi_1(M)$ with $\pi_1(M,x_0)$, $\pi_1(\what M)$
with $\pi_1(M,\hat x_0)$, and $\hat\Lambda$ with
$\what\Upsilon_*\pi_1(\what M,\hat x_0)$.

By the arguments in the proof of Proposition~\ref{Action on tree
gives 3 filtered ends prop}, some end of $\what M$ is a \comp of
$\what\Psi\inv(\hat l((-\infty,1)))$. Forming the union of this
end with a large relatively \cpt domain in $\what M$, we get an
end $E_0$ of $\what M$ \st $\hat x_0\in E_0$ and
$\what\Psi(E_0)\subset\hat l((-\infty,m))$ for some positive
integer~$m$. We may now choose a \cinf \rel \cpt domain $\Omega$
in $\what M$ \st $\hat x_0\in\Omega$, $\partial E_0\subset\Omega$,
$\Omega\cap E_0$ is \conns, and $\image{\pi_1(\Omega,\hat
x_0)\to\pi_1(M,x_0)}\supset\Lambda_{v_m}=\tau^m\Gamma\tau^{-m}$
(here we have used the fact that $\Gamma$ is finitely generated).

We now show that, for $E=E_0\cup\Omega$, $\image{\pi_1(E,\hat
x_0)\to\pi_1(M,x_0)}$ is finitely generated. Since $\Omega$ is a
\cinf domain, we may choose a domain $\Theta$ \st
$\Omega\Subset\Theta\Subset\what M$ and $\overline\Omega$ is a
strong deformation retract of $\Theta$. Given a loop $\beta$ in
$E$ based at $\hat x_0$, we may choose a partition
$0=s_0<t_0<s_1<t_1<\cdots<s_k<t_k=1$ \st
$\beta(s_j),\beta(t_j)\in\overline\Omega$ and
$\beta([s_j,t_j])\subset\Theta$ for $j=0,\dots,k$, and
$\beta((t_{j-1},s_j))\subset E\sm\Omega=E_0\sm\Omega$ for
$j=1,\dots,k$. For each $j=1,\dots,k$, we may choose a path
$\delta_j$ in $E_0\cap\overline\Omega$ from $\beta(t_{j-1})$ to
$\beta(s_j)$ (since these points lie in
$E\cap\partial\Omega=E_0\cap\partial\Omega$ and
$E_0\cap\overline\Omega$ is \conns), and, for each $j=0,\dots,k$,
we may choose a path $\epsilon_j$ in $E_0\cap\overline\Omega$ from
$x_0$ to $\beta(t_j)$. The loop $\beta$ is homotopic to the loop
$\eta_0*\kappa_1*\eta_1*\kappa_2*\cdots*\eta_{k-1}*\kappa_k*\eta_k$,
where, $\eta_0=\beta\restrict{[s_0,t_0]}*\epsilon_0\inv$,
$\kappa_j=\epsilon_{j-1}*\beta\restrict{[t_{j-1},s_j]}*
\delta_j\inv*\epsilon_{j-1}\inv$ for $j=1,\dots,k$,
$\eta_j=\epsilon_{j-1}*\delta_j*\beta\restrict{[s_j,t_j]}*
\epsilon_j\inv$ for $j=1,\dots,k-1$, and
$\eta_k=\epsilon_{k-1}*\delta_k*\beta\restrict{[s_k,t_k]}$. For
each $j=0,\dots,k$, $\eta_j$ is contained in $\Theta$ and,
therefore, homotopic to a loop in $\Omega$. For each
$j=1,\dots,k$, $\kappa_j$ is a loop in $E_0$ and hence the lifting
$\til\kappa_j$ to $\wtil M$ with $\til\kappa_j(1)=\til x_0$ lies
in $\wtil\Psi\inv(D)$; where $D\subset\theta\inv(\hat
l((-\infty,m+1)))$ is the \comp of $T\sm\set{v_{m+1}}$ containing
$v_0$. Therefore, the element
$\lambda_j=[\what\Upsilon(\kappa_j)]\in\pi_1(M,x_0)$ satisfies
\(\lambda_j\cdot v_0=\wtil\Psi(\lambda_j\cdot\til
x_0)=\wtil\Psi(\til\kappa_j(0))\in D\). Thus, by the remarks
preceding this proof, we have $\lambda_j\in\Lambda_{v_m}$ and
hence, by construction, $\kappa_j$ is homotopic to a loop in
$\Omega$. Since $\Omega\subset E$, we get
\[
\image{\pi_1(\Omega,\hat
x_0)\to\pi_1(M,x_0)}\subset\image{\pi_1(E,\hat
x_0)\to\pi_1(M,x_0)}\subset\image{\pi_1(\Omega,\hat
x_0)\to\pi_1(M,x_0)}.
\]
Thus we get equality of the above finitely generated groups.

Finally, given a \rel \cpt \cinf domain $\Omega'$ in $\what M$
containing $\partial E$, we may choose a \cinf \rel \cpt domain
$\Omega''$ in $\what M$ \st $\Omega\Subset\Omega''$,
$E_0\cap\Omega''$ is \conns, and $\Omega''\sm E=\Omega'\sm E$. The
above argument applied to the ends $E_0\subset
E'=E\cup\Omega'=E_0\cup\Omega''$ gives finite generation of
$\image{\pi_1(E',\hat x_0)\to\pi_1(M,x_0)}$.
\end{pf}

\begin{rmk}
If $M=S$ is a \cpt Riemann surface of genus $2$ and
$\Lambda=\pi_1(M)$ is expressed as a properly ascending HNN
extension with (infinitely generated) base group $\Gamma$ as in
the remarks preceding the above lemma, then the corresponding
covering $\what M$ is, topologically, an infinite tube with an
infinite sequence of handles attached; and $\Z$ acts transitively
and freely on this collection of handles. In particular, for any
end $E$, the image of $\pi_1(E)$ is {\it not} finitely generated.
This illustrates how the lemma and the theorem fail if the base
group is not assumed to be finitely generated.
\end{rmk}

\begin{pf*}{Completion of the proof of Theorem~\ref{HNN and Thompson not
Kahler theorem from intro}} It remains to show that any properly
ascending HNN extension $\Lambda$ with finitely generated base
group $\Gamma$ and stable letter $\tau$ is {\it not} K\"ahler. For
this, we assume that $\Lambda=\pi_1(X)$ for some \con \cpt
K\"ahler manifold $X$ and reason to a contradiction.

Let $\hat\Lambda=\bigcup_{m\in\Z}\tau^m\Gamma\tau^{-m}$ and let
$\Upsilon\colon\what X\to X$ be a \con (Galois) covering space
with $\Upsilon_*\pi_1(\what X)=\hat\Lambda$. According to
Lemma~\ref{HNN covering finitely gen lem},
$\image{\pi_1(E)\to\pi_1(\what X)}$ is finitely generated for some
end $E$ of $\what X$. In particular, this image group is of
infinite index and hence, by Theorem~\ref{Main Theorem from
Introduction}, there exists a proper \holo mapping $\hat\Phi$ with
\con fibers of $\what X$ onto a Riemann surface $\what S$, with
$e(\what S)=2$. The action of
$\Z\cong\Lambda\big\slash\hat\Lambda$ on $\what X$ descends to a
properly discontinuous action of $\Z$ on $\what S$ and, since $\Z$
is torsion-free, the action is free. Thus we get a commutative
diagram
\begin{center}\begin{picture}(230,70)
\put(40,9){$\Z\big\slash\what X=X$} \put(93,13){\vector(1,0){45}}
\put(140,9){$S\equiv \Z\big\backslash\widehat S$}
\put(80,49){$\widehat X$} \put(110,55){$\widehat\Phi $}
\put(93,53){\vector(1,0){45}} \put(140,49){$\widehat S$}
\put(149,31){$\theta $} \put(144,47){\vector(0,-1){28}}
\put(70,31){$\Upsilon$} \put(84,47){\vector(0,-1){28}}
\put(110,15){$\Phi$}
\end{picture}
\end{center}
where $\theta$ is an (unbranched) infinite covering map, $S$ is a
\cpt Riemann surface, and $\Phi $ is a surjective proper \holo map
with \con fibers. In particular, the singular fibers of $\hat\Phi$
are precisely the liftings of the singular fibers of $\Phi$.

By the second conclusion of Lemma~\ref{HNN covering finitely gen
lem}, we may assume that $E=\what\Phi\inv(F)$ for some end $F$ of
$\what S$ with smooth boundary. Hence
$\image{\pi_1(F)\to\pi_1(\what S)}$ is finitely generated and,
therefore, the Riemann surface $F$ is of finite type. Thus $F$
contains an end $D$ which is isomorphic to either a punctured disk
or an annulus.

Moreover, $\image{\pi_1^{\text{\rm orb}}(\what\Phi\restrict
E)\to\pi_1^{\text{\rm orb}}(\what\Phi)}$ is finitely generated. If
$F$ contains a point $p$ \st the greatest common divisor of the
multiplicities of the components of the divisor $\hat\Phi\inv(p)$
is greater than $1$, then $F$ contains infinitely many such points
(since the singular fibers are liftings of the singular fibers of
$\Phi$). Hence $\pi_1^{\text{\rm orb}}(\hat\Phi\restrict E)$
contains an infinite free product of finite cyclic groups which
injects into $\pi_1^{\text{\rm orb}}(\hat\Phi)$. But this
contradicts the finite generation property. Thus $F$ cannot
contain such points and therefore, applying the action of $\Z$ to
move the fibers over points in $\what S\sm F$ into $F$, we see
that, for each fiber of $\hat\Phi$ or $\Phi$, the greatest common
divisor of the multiplicities of the components is $1$. Thus
$\pi_1(\what S)=\pi_1^{\text{\rm orb}}(\hat\Phi)$ is infinitely
generated and hence, in particular, $\what S\neq\C^*$. Hence $S$
is a curve of genus $g\geq 2$ and therefore, since a punctured
disk end for $\what S$ would imply the existence of a parabolic
element in $\pi_1(S)$, the end $D\subset F$ must be isomorphic to
an annulus.

On the other hand, $\what X$ admits a proper \plh \fn
$\alpha\colon\what X\to\R$ (for example, one may take $\alpha$ to
be the integral of the lifting of a closed real \harm $1$-form on
$X$ which integrates to $1$ on $\tau$, and $0$ on $\hat\Lambda$).
The \fn $\alpha$ descends to a proper (pluri)\harm \fn
$\beta\colon\what S\to\R$. This implies that $\what S$ is
parabolic. For assuming, as we may, that $\pm 1$ is a regular
value for $\beta$, the \harm \fn $|\beta|/R$ on the set
$\setof{x\in\what S}{1<|\beta(x)|<R}$ vanishes on $\beta\inv(\pm
1)$ and is equal to $1$ on $\beta\inv (\pm R)$. Letting
$R\to\infty$, we see that the \harm measure of the ideal boundary
of $\what S$ with respect to $\setof{x\in\what S}{|\beta(x)|<1}$
vanishes. Thus we have again arrived at a contradiction and,
therefore, $\Lambda$ cannot be K\"ahler.
\end{pf*}

\section{Principal functions and the Evans-Selberg potential}\label{Principal fns Evans-Selberg section}

For the convenience of the reader, the proof of Sario's existence
theorem of principal functions \cite{RS} and Nakai's construction
of the Evans-Selberg potential \cite{Nakai1}, \cite{Nakai2},
\cite{SaNo} are provided in this section. These facts were applied
in \cite{NR1}. However, it is difficult to find proofs for a
general oriented Riemannian manifold in a convenient form in the
literature. This section will not appear in the version submitted
for publication. Throughout this section, $(M,g)$ will denote a
\con non\cpt oriented Riemannian manifold of dimension~$n>2$.

\subsection*{\bf A. Principal functions}

Throughout this subsection, we will assume that $(M,g)$ is
parabolic and $M_0$ will denote a $\cinf $ \rel \cpt domain
in~$M$. In this subsection, we recall the following theorem of
Sario~\cite {RS}:

\begin{thm}[Sario]\label{Appendix-Sario Principal functions theorem}
If $u$ is a \cont function on~$M\setminus M_0$ which is \harm on
$M\setminus \overline M_0$ and satisfies the flux condition:
$$
\int _{\partial M_0} \ndof u \, =0,
$$
then there exists a \harm function~$v$ on~$M$ \st $u-v$ is bounded
on~$M\setminus M_0$.
\end{thm}
\begin{rmks} 1. This theorem is also true in the hyperbolic case, but we
will only need it for the parabolic case.

\noindent 2. The flux condition gives
$$
\int _{\partial \Omega } \ndof u \, =0
$$
for some (hence for every) $\cinf$ \rel \cpt domain~$\Omega $
containing $\overline M_0$.
\end{rmks}
For the rest of this subsection, we will assume, as we may, that
$M\sm\overline{M_0}$ has no \rel \cpt \concompsns.

\begin{lem}\label{Appendix-Linear map Dirichlet lemma}
There exists a linear map
$$
L : C^0 (\partial M_0)\to \text {{\rm Harm}}\, (M\setminus
\overline M_0) \cap C^0 (M\setminus M_0)
$$
\stns, for every \cont \fn $\alpha $ on $\partial M_0$, we have:
\begin{enumerate}
\item[(i)] $(L\alpha )\restriction _{\partial M_0}=\alpha $,

\item[(ii)] $\min _{\partial M_0} \alpha \leq L\alpha \leq \max
_{\partial M_0} \alpha $, and

\item[(iii)] $\int _{\partial M_0} \ndof {(L\alpha )} =0$.
\end{enumerate}
\end{lem}

\begin{rmk}
Such an operator~$L$ is called a {\it normal operator}
for~$M\setminus M_0$.
\end{rmk}

\begin{pf*}{Proof of Lemma~\ref{Appendix-Linear map Dirichlet lemma}} Let $\seq Mk _{k=1} ^{\infty }$ be
a fixed \exh of~$M$ by $\cinf $ \rel \cpt domains
containing~$\overline M_0$, and, for each positive integer~$k$,
let~$v_k$ be the \hm of~$\partial M_k$ \wrt $M_k\setminus
\overline M_0$. Given a \cont function~$\alpha $ on~$\partial
M_0$, for each positive integer~$k$, let $w_k \in \Harm
(M_k\setminus \overline M_0) \cap C^0(\overline M_k\setminus M_0)$
be the function which vanishes on~$\partial M_k$ and is equal
to~$\alpha $ on~$\partial M_0$.  Since~$M$ is parabolic, $v_k
\searrow 0$ uniformly on \cpt subsets of~$M\setminus M_0$. Since
the sequence of nonnegative functions $\{ w_k -(\min _{\partial
M_0} \alpha )(1-v_k) \} $ is bounded and nondecreasing, Harnack's
principle (Lemma~1.3) implies that the sequence~$\seq wk $
converges uniformly on \cpt subsets of~$M\setminus M_0$ to a
function $w\in \Harm (M\setminus \overline M_0) \cap
C^0(M\setminus M_0)$.

We set $L\alpha \equiv w$. It remains to verify the properties
(i),~(ii),~and~(iii). The property~(i) is clear. The property~(ii)
follows from the minimum principle for super\harm functions on
parabolic manifolds (see Sect.~1). Finally, to verify the
property~(iii), we observe that
\begin{align*} \int _{\partial M_1}
\ndof {w_k} &=\int _{\partial M_k} \ndof {w_k} =\int _{\partial
M_k} v_k\ndof {w_k}\\
&=\int _{\partial (M_k\setminus M_1)} v_k\ndof {w_k}
+ \int _{\partial M_1} v_k\ndof {w_k} \\
&=\int _{\partial (M_k\setminus M_1)} \ndof {v_k}w_k
+ \int _{\partial M_1} v_k\ndof {w_k} \\
&=-\int _{\partial M_1} \ndof {v_k}w_k + \int _{\partial M_1}
v_k\ndof {w_k}\to 0 \quad \text {as} \quad k \to \infty.
\end{align*}
Therefore
$$
\int _{\partial M_0} \ndof w =\int _{\partial M_1} \ndof w =\lim
_{k\to \infty } \int _{\partial M_1} \ndof {w_k} =0.
$$
Thus $L$ is a normal operator.
\end{pf*}

Next we recall the following important consequence of Harnack's
principle:

\begin{lem}[$q$-lemma]\label{Appendix q-lemma} Given a \cpt subset~$K$ of an oriented
Riemannian manifold~$N$, there exists a constant~$q\in (0,1)$ \st
$$
q\inf _N u \leq u(x) \leq q\sup _Nu \quad \forall \, x\in K
$$
for every \harm function~$u$ on~$N$ which changes sign on~$K$.
\end{lem}

\begin{pf*}{Proof of Theorem~\ref{Appendix-Sario Principal functions theorem}}
Since $L(u\restriction _{\partial
M_0})$ is bounded, $u-v$ will be bounded for a given function~$v$
if and only if $u-L(u\restriction _{\partial M_0})-v$ is bounded.
Hence we may assume without loss of generality that~$u$ vanishes
on~$\partial M_0$.

Fix a $\cinf $ domain $\Omega $ with $M_0\Subset \Omega \Subset M$
and let $K : C^0 (\partial \Omega ) \to \Harm (\Omega ) \cap C^0
(\overline \Omega )$ be the linear map which associates to each
(finitely) \cont \fnns~$\alpha $ on~$\partial \Omega $ the \cont
function on~$\overline \Omega $ which is \harm on~$\Omega $ and
equal to~$\alpha $ on~$\partial \Omega $. Clearly, it suffices to
find a \harm function~$v$ on~$M$ such that
$$
u-v=-L((K(v\restriction _{\partial \Omega })) \restriction
_{\partial M_0})
$$
on~$M\setminus M_0$; since the right-hand side is bounded. For
this, we need only find a solution $\alpha \in C^0(\partial \Omega
)$ to
\[
(1)\qquad\qquad\qquad\qquad\qquad\qquad\qquad\qquad  (I-J)\alpha =
u\restriction _{\partial
\Omega}\qquad\qquad\qquad\qquad\qquad\qquad\qquad\qquad
\]
where
$J : C^0(\partial \Omega ) \to C^0(\partial \Omega )$ is the \cont
linear operator defined by
$$
J\alpha =(L((K\alpha )\restriction {\partial M_0})) \restriction
_{\partial \Omega } \quad \forall \, \alpha \in C^0(\partial
\Omega ) .
$$
For the function~$v$ defined by
$$
v\restriction _{\overline \Omega } \equiv K\alpha \quad \text
{and} \quad v\restriction _{(M\setminus M_0)} \equiv u+L((K\alpha
)\restriction {\partial M_0})
$$
will then have the required properties. For this, we will prove
uniform convergence of the series
$$
\sum _{m=0}^{\infty } J^m (u\restriction _{\partial \Omega }).
$$
The sum $\alpha \in C^0(\partial \Omega )$ will then be a solution
to equation~(1).

We first prove two identities. Let $w$ be the \hm of $\partial
\Omega $ \wrt $\Omega \setminus \overline M_0$. Then
\begin{enumerate}
\item[(a)] $\int _{\partial M_0} \beta \ndof {w}
=\int _{\partial \Omega } (L\beta )\ndof {w} \quad \forall \,
\beta \in C^0(\partial M_0) $; and

\item[(b)] $\int _{\partial M_0} (K\beta ) \ndof {w} =\int
_{\partial \Omega } \beta \ndof {w} \quad \forall \, \beta \in
C^0(\partial \Omega)$.
\end{enumerate}
For the proofs, it suffices to consider~$\cinf $ functions since
$\cinf (\partial M_0)$ and $\cinf (\partial \Omega )$ are dense
in~$C^0(\partial M_0)$ and $C^0 (\partial \Omega )$, respectively.
If $\beta \in \cinf (\partial M_0) $, then $L\beta $ is \harm
on~$M\setminus \overline M_0$ with~$\cinf $ boundary data.
Therefore $L\beta \in \cinf (M\setminus M_0)$ and
\begin{align*}
\int _{\partial M_0} \beta \ndof {w} =\int _{\partial M_0} (L\beta
) \ndof {w} &=\int _{\partial \Omega } (L\beta )\ndof {w} -\int
_{\partial \Omega } \ndof {(L\beta )}w
+\int _{\partial M_0} \ndof {(L\beta )} w \\
&=\int _{\partial \Omega } (L\beta )\ndof {w}
\end{align*}
since $w\equiv 0$ on~$\partial M_0$, $w\equiv 1$ on~$\partial
\Omega $, and $L\beta $ satisfies the flux condition~(iii) of
Lemma~\ref{Appendix-Linear map Dirichlet lemma}. Thus the
identity~(a) is proved. The proof of~(b) is similar.

Next we show that the function $KJ^m(u\restriction _{\partial
\Omega })$ changes sign on~$\partial M_0$ for every nonnegative
integer~$m$ (which will allow us to apply the $q$-lemma). In fact,
we show that
$$
\int _{\partial M_0} KJ^m(u\restriction _{\partial \Omega })\ndof
{w} =0.
$$
Since $\ndof {w}>0$ on~$\partial M_0$, it will follow that
$KJ^m(u\restriction _{\partial \Omega })$ changes sign
on~$\partial M_0$. For $m=0$, the identity~(b) implies that
$$
\int _{\partial M_0} (K(u\restriction _{\partial \Omega })) \ndof
{w} =\int _{\partial \Omega } u \ndof {w} = \int _{\partial M_0} u
\ndof {w} +\int _{\partial \Omega } \ndof {u} w -\int _{\partial
M_0} \ndof {u} w=0
$$
since $u\equiv w\equiv 0$ on~$\partial M_0$, $w\equiv 1$
on~$\partial \Omega $, and $u$ satisfies the flux condition. If
$m>0$, then, by the identities,
\begin{align*}
\int _{\partial
M_0} (KJ^m(u\restriction _{\partial \Omega })) \ndof {w} &=\int
_{\partial \Omega }
J^m(u\restriction _{\partial \Omega }) \ndof {w}\\
&=\int _{\partial \Omega } L((KJ^{m-1}(u\restriction _{\partial
\Omega }) \restriction _{\partial M_0})
\ndof {w} \\
&=\int _{\partial M_0} KJ^{m-1}(u\restriction _{\partial \Omega })
\ndof {w}.
\end{align*}
The claim now follows by induction on~$m$. Thus
$KJ^m(u\restriction _{\partial \Omega })$ is a \harm \fn
on~$\Omega $ which changes sign on the \cpt set~$\partial M_0$.
Therefore, by the $q$-lemma, there exists a number~$q\in (0,1)$
(independent of~$m$) \stns , for each nonnegative integer~$m$,
\begin{align*}
q\min _{\partial \Omega }(KJ^m(u\restriction _{\partial
\Omega }))
&=q\inf _{\Omega } (KJ^m(u\restriction _{\partial \Omega })) \\
&\leq KJ^m(u\restriction _{\partial \Omega })\restriction _{\partial M_0} \\
&\leq q \sup _{\Omega } (KJ^m(u\restriction _{\partial \Omega })) \\
&=q\max _{\partial \Omega }(KJ^m(u\restriction _{\partial \Omega
})).
\end{align*}
Hence
$$
q\min _{\partial \Omega }(J^m(u\restriction _{\partial \Omega }))
\leq (KJ^m(u\restriction _{\partial \Omega }))\restriction
_{\partial M_0} \leq q\max _{\partial \Omega }(J^m(u\restriction
_{\partial \Omega }))
$$
since $KJ^m(u\restriction _{\partial \Omega }) =J^m(u\restriction
_{\partial \Omega })$ on~$\partial \Omega $. By the property~(ii)
of the normal operator~$L$ stated in Lemma~\ref{Appendix-Linear
map Dirichlet lemma}, we get
\begin{align*} \max _{\partial \Omega }(J^{m+1}(u\restriction _{\partial
\Omega })) &=\max _{\partial \Omega }
(L((KJ^m(u\restriction _{\partial \Omega }))\restriction _{\partial M_0}))\\
&\leq
\max _{\partial M_0}(KJ^m(u\restriction _{\partial \Omega })) \\
&\leq q\max _{\partial \Omega }(J^m(u\restriction _{\partial
\Omega })).
\end{align*}
Similar inequalities hold for the corresponding minima. Therefore,
by induction on~$m$, we get
$$
q^m\min _{\partial \Omega } u \leq \min _{\partial \Omega
}(J^m(u\restriction _{\partial \Omega })) \leq \max _{\partial
\Omega }(J^m(u\restriction _{\partial \Omega })) \leq q^m\max
_{\partial \Omega } u.
$$
Since $0<q<1$, the series
$$
\sum _{m=0}^{\infty } J^m (u\restriction _{\partial \Omega })
$$
converges uniformly on~$\partial \Omega $ and the proof is
complete.
\end{pf*}

\subsection*{B. Green's potentials and the energy principle}

In this subsection, we recall the facts concerning Green's
potentials which are used in Nakai's construction of the
Evans-Selberg potential. We include sketches of some of the
proofs. For more details the reader may refer to \cite{Ancona} and
\cite{Maeda}. Throughout this subsection, we will assume that
$(M,g)$ is hyperbolic with Green's function~$G$.

A nonconstant nonnegative super\harm \fn~$\vphi $ on~$M$ is called
a {\it potential} if, for every harmonic \fnns~$u$ with $0\leq
u\leq\vphi$ on~$M$, we have $u\equiv 0$. For example, $G_x$ is a
potential for each point~$x\in M$.

Given a positive regular Borel measure~$\mu $ supported in a \cpt
subset~$K$ of~$M$, the function $G_\mu : M \to [0,+\infty )$ given
by
$$
G_\mu (x)= \int _M G(x,y) \, d\mu (y) \quad \forall \, x\in M
$$
is a potential called the {\it Green's potential of}~$\mu $.
Moreover, $G_\mu $ is \harm on an open subset of~$U$ if and only
if $\mu (U)=0$.

We will need two well-known facts concerning Green's potentials.
The first is the following:

\begin{lem}\label{Appendix potential gives borel measure lemma}
Suppose $\vphi $ is a potential on~$M$ which is harmonic on the
complement of some \cpt subset~$K$ of Lebesgue measure zero
in~$M$. Then there exists a positive regular Borel measure~$\mu$
on~$M$ (supported in~$K$) such that $\vphi =G_{\mu}$.
\end{lem}

\begin{rmk} In fact, every potential is a Green's potential, but we will
only need to consider this special case. For the general case, see
\cite{Ancona}. The general fact is the key element in the proof
that the distributional Laplacian of a sub\harm \fn is a positive
measure. Moreover, one can prove that the measure~$\mu$ is unique.
\end{rmk}

\begin{pf*}{Proof of Lemma \ref{Appendix potential gives borel measure lemma}}
We first show that we may assume without loss of generality that
$M$ is a $\cinf$ \rel \cpt domain in an oriented Riemannian
manifold~$(M',g')$, $g=g'\restrict M$, and the restriction of
$\vphi $ to $M\setminus K$ vanishes smoothly at~$\partial M$. Let
$\seq Mk$ be an \exh of $M$ by $\cinf $ \rel \cpt domains
containing~$K$, and, for each~$k$, let $G_k$ be the Green's
function on~$M_k$. Given an upper semi\cont function~$\alpha $ and
a $\cinf $ domain~$\Omega $ in~$M$ with Green's
function~$\widetilde G$, we denote by $\rho _\Omega (\alpha )$ the
\harm function given by
$$
\rho _{\Omega }(\alpha )(x) =\int _{\partial \Omega } \ndof
{\widetilde G_x} (y) \alpha (y) d\sigma (y).
$$
Since $\vphi $ is harmonic and therefore of class~$\cinf $ on
$M\setminus K$, for each~$k$, the function
$$
\psi _k \equiv \vphi -\rho _{M_k}(\vphi )
$$
is a potential which is $\cinf $ on~$\overline M_k\setminus K$,
harmonic on~$M_k\setminus K$, and identically equal to~$0$
on~$\partial M_k$. Moreover, since
$$
\rho _{M_{k+1}}(\vphi ) \leq \vphi =\rho _{M_k}(\vphi )
$$
on $\partial M_k$ and the functions $\rho _{M_{k+1}}(\vphi )$ and
$\rho _{M_k}(\vphi )$ are harmonic, we have
$$
0\leq \rho _{M_{k+1}}(\vphi ) \leq \rho _{M_k}(\vphi )\leq \vphi
\quad \text {and}\quad \psi _{k+1} \geq {\psi _k}\geq 0
$$
on~$\overline M_k$. By Harnack's principle, the sequence of
functions~$\{ \rho _{M_k}(\vphi )\} $ converges uniformly on \cpt
subsets of~$M$ to some \harm function~$h$ satisfying $0\leq h\leq
\vphi $. Since $\vphi $ is a potential, $h$ must vanish
identically and therefore $\vphi -\psi _k \searrow 0$ uniformly on
\cpt sets as $k \to \infty $.

Now suppose that, for each~$k$, $\psi _k=(G_k)_{\mu _k}$ for some
positive regular Borel measure~$\mu _k$. Then $\mu _k$ is
supported in~$K$ since $\psi _k$ is \harm on~$M\setminus K$.
Fixing any point~$x_0\in M\setminus K$, we get
\begin{align*} \mu _k(M)&=\int _K \, d\mu _k (y) \leq \biggl( \min
_K(G_k)_{x_0})\biggr) \inv
\int _K G_k(x_0,y)\, d\mu _k (y) \\
&=\biggl( \min _K(G_k)_{x_0})\biggr) \inv \psi _k(x_0) \leq
\biggl( \min _K(G_k)_{x_0})\biggr) \inv \vphi (x_0) < +\infty .
\end{align*}
Hence the sequence $\{ \mu _k(K) \} $ is bounded. Thus, by passing
to a subsequence, we may assume that the sequence of measures
$\seq {\mu } {k}$ converges {\it weakly} to some positive regular
Borel measure~$\mu $ supported in~$K$; that is,
$$
\int _M \alpha \, d\mu _k \to \int _M \alpha \, d\mu
$$
for every \cont function $\alpha $ on $M$. On the other hand, for
every point~$x\in M\setminus K$, $G_x$ is \cont on~$K$ and $G-G_k
\to 0$ uniformly on compact subsets of~$M\times M$. Therefore
\begin{align*} G_\mu (x) &=\int _K G(x,y) \, d\mu (y)
=\lim _{k \to \infty }\int _M G(x,y) \, d\mu _k(y) \\
&=\lim _{k \to \infty } \biggl[ \int _K (G(x,y) -G_k(x,y))\, d\mu
_k(y) +\int _K G_k(x,y)\, d\mu _k(y)
\biggr] \\
&=\lim _{k \to \infty } (G_k)_{\mu _k} (x) =\lim _{k \to \infty }
\psi _k(x) =\vphi (x).
\end{align*}
Therefore $G_\mu =\vphi $ on $M\setminus K$. But $G_\mu $ and
$\vphi $ are superharmonic and the set~$K$ has Lebesgue measure
zero. Therefore $G_\mu \equiv \vphi $ on~$M$. Thus we may assume
that $M$ is a $\cinf $ \rel \cpt domain in an oriented Riemannian
manifold~$(M',g')$, $g=g'\restriction _M$, $\vphi $ is of
class~$\cinf $ on $\overline M\setminus K$, and $\vphi \equiv 0$
on~$\partial M$.

Next, we approximate~$\vphi $ by Lipschitz \cont potentials. Let
$\seq \Omega m$ be a sequence of~$\cinf $ \rel \cpt open sets
in~$M$ such that
$$
\Omega _m \Supset \Omega _{m+1} \quad \forall \, m =1,2,3, \dots
\quad \text {and} \quad \bigcap _{m=1}^{\infty } \Omega _m =K.
$$
For each $m$,let $\vphi _m$ be the Lipschitz \cont potential
defined by
$$
\vphi _m(x) =(\cal P_{\Omega _m} (\vphi ))(x) =\left\{
\begin{aligned} \vphi (x)
&\qquad \text {if } x\in M\setminus \Omega _m\\
(\rho _{\Omega _m}(\vphi ))(x)&\qquad \text {if } x\in \Omega _m.
\end{aligned}
\right.
$$
Then $\vphi _m$ is \harm on~$M\setminus \partial \Omega _m$,
$\vphi _m$ is smooth up to the boundary on $\overline M\setminus
\Omega _m$ and on $\overline \Omega _m$, and $\vphi _m\leq \vphi
_{m+1}\leq \vphi $ on~$M$. In particular, the sequence $\seq \vphi
m$ converges to a super\harm function on~$M$. Since this function
is equal to~$\vphi $ on~$M\setminus K$ and~$K$ has Lebesgue
measure zero, we get $\vphi _m \nearrow \vphi $ on~$M$. We will
prove that each function $\vphi _m$ is a Green's potential and
then pass to the limit to obtain a Green's potential equal
to~$\vphi $.

For each positive integer~$m$, let $u_m=\rho _{\Omega _m}(\vphi
)=\vphi _m \restriction _{\overline \Omega _m}$. Then
$$
\nd (u_m-\vphi ) \geq 0 \quad \text {on } \partial \Omega _m
$$
because $u_m-\vphi \leq 0$ on $\overline \Omega _m$ and $u_m-\vphi
\equiv  0$ on~$\partial \Omega _m$. Let $\tau_{n-1}$ denote the
volume of the unit sphere in $\R^n$. It follows that the
distributional Laplacian of the function~$-\vphi _m /((n-2)\tau
_{n-1})$ determines a positive regular Borel measure
$$
d\mu _m \equiv \frac {1}{(n-2)\tau _{n-1}} \biggl[ \ndof {u_m}
-\ndof {\vphi } \biggr] \, d\sigma _m ,
$$
where $d\sigma _m$ is the volume element on~$\partial \Omega _m$.
For if~$\alpha $ is a~$\cinf $ function on $\overline M$ which
vanishes on~$\partial M$, then, since $u_m\equiv \vphi $
on~$\partial \Omega _m$,
\begin{align*} \int _M(\lap \alpha )
\vphi _m \, dV &=\int _{M\setminus \overline \Omega _m}(\lap
\alpha ) \vphi  \, dV
+\int _{\Omega _m}(\lap \alpha ) u_m \, dV \\
&=-\int _{\partial \Omega _m} \ndof {\alpha } \vphi \, d\sigma _m
+ \int _{\partial \Omega _m} \alpha  \ndof {\vphi } \, d\sigma _m \\
&\qquad\qquad +\int _{\partial \Omega _m} \ndof {\alpha } u_m \,
d\sigma _m
- \int _{\partial \Omega _m} \alpha  \ndof {u_m} \, d\sigma _m \\
&=-\int _{\partial \Omega _m} \alpha \biggl[ \ndof {u_m} -\ndof
{\vphi } \biggr] \, d\sigma _m
\end{align*}
Moreover, $\vphi _m=G_{\mu _m}$. For if $x\in M\setminus \overline
\Omega _m$ and $G^*$ is the Green's function on~$M\setminus
\overline \Omega _m$, then, since $G_x $ is harmonic on~$\Omega
_m$, $G_x -G^*_x $ and $\vphi $ are harmonic on~$M\setminus
\overline \Omega _m$, $G_x$, $G^*_x $, and $\vphi $ are equal
to~zero on~$\partial M$, and $G^*_x =0$ and $\vphi =u_m$
on~$\partial \Omega _m$, we have
\begin{align*}
(n-2)\tau _{n-1}
G_{\mu _m }(x)
&= (n-2)\tau _{n-1} \int _M G(x,y) \, d\mu _m (y) \\
&= \int _{\partial \Omega _m} G(x,y)
\biggl[ \ndof {u_m} (y)-\ndof {\vphi } (y)\biggr] \, d\sigma _m(y)\\
&= \int _{\partial \Omega _m} \ndof {G_x} u_m  \, d\sigma _m
-\int _{\partial \Omega _m} G_x  \ndof {\vphi }   \, d\sigma _m \\
&=\int _{\partial \Omega _m} \ndof {G_x} \vphi   \, d\sigma _m
-\int _{\partial \Omega _m} G_x  \ndof {\vphi }
\, d\sigma _m  \\
&=\int _{\partial \Omega _m} \ndof {(G_x-G^*_x)} \vphi   \,
d\sigma _m -\int _{\partial \Omega _m} (G_x -G^*_x )\ndof {\vphi }
\, d\sigma _m  \\
&\qquad\qquad  +\int _{\partial \Omega _m}
\ndof {G^*_x} \vphi  \, d\sigma _m  \\
&=\int _{\partial M} \ndof {(G_x-G^*_x)} \vphi   \, d\sigma _m
-\int _{\partial M} (G_x -G^*_x )\ndof {\vphi }
\, d\sigma _m  \\
&\qquad\qquad  -\int _{\partial (M\setminus \Omega _m)}
\ndof {G^*_x} \vphi  \, d\sigma _m  \\
&=-\int _{\partial (M\setminus \Omega _m)}
\ndof {G^*_x} \vphi  \, d\sigma _m  \\
&=(n-2)\tau _{n-1}\vphi (x) =(n-2)\tau _{n-1}\vphi _m(x).
\end{align*}
Therefore $G_{\mu _m}\equiv \vphi \equiv \vphi _m$ on $\overline
M\setminus \overline \Omega _m$. Similarly (by working with the
Green's function of~$\Omega _m$) one can show that $G_{\mu
_m}\equiv u_m\equiv \vphi _m$ on $\Omega _m$. Since $G_{\mu _m}$
and $\vphi _m$ are superharmonic and $\partial \Omega _m$ has
Lebesgue measure zero, $G_{\mu _m}$ and $\vphi _m$ are equal
on~$\overline M$. In particular, $G_{\mu _m}=\vphi _m \nearrow
\vphi $ as $m \to \infty $.

Since $u_m$ is \harm on $\Omega _m$ and $\vphi $ is \harm
on~$M\setminus K\Supset \Omega _1\setminus \Omega _m$, we have
\begin{align*} \mu _m(M)&=\int _{\partial \Omega _m} \, d\mu _m = \frac
{1}{(n-2)\tau _{n-1}}\int _{\partial \Omega _m} \biggl[ \ndof
{u_m} -\ndof {\vphi }
\biggr] \, d\sigma _m  \\
&=-\frac {1}{(n-2)\tau _{n-1}}\int _{\partial \Omega _m}
\ndof {\vphi } \, d\sigma _m  \\
&=-\frac {1}{(n-2)\tau _{n-1}}\int _{\partial \Omega _1} \ndof
{\vphi } \, d\sigma _m .
\end{align*}
Thus $\{ \mu _m (M)\} $ is a bounded (in fact, constant) sequence
and therefore, by replacing $\seq \mu m $ by a subsequence (if
necessary), we may assume that $\seq \mu m $ converges weakly to a
positive regular Borel measure~$\mu $ on~$M$ supported
in~$K=\bigcap _m  \Omega _m$. If $x\in M\setminus K$, then, for
$m_0$ a sufficiently large positive integer, the function~$G_x$ is
finite and \cont on $\overline \Omega _{m_0}$ and, for $m>m_0$,
$$
\vphi _m(x)=G_{\mu _m}(x) =\int _M G(x,y)\, d\mu _m (y) =\int
_{\overline \Omega _{m_0}} G(x,y)\, d\mu _m (y).
$$
Passing to the limit as $m \to \infty $, we get
$$
\vphi (x) =\int _{\overline \Omega _{m_0}} G(x,y)\, d\mu (y).
$$
Thus $\vphi \equiv G_{\mu }$ on $M\setminus K$. Since $K$ has
Lebesgue measure zero and these functions are superharmonic,
$\vphi =G_{\mu }$ on~$M$. Thus the lemma is proved.
\end{pf*}

The second fact which is needed in Nakai's construction of the
Evans-Selberg potential (\cite{Nakai1} and \cite{Nakai2}) is the
energy principle, which is an analogue of the Schwarz inequality.
Given two compactly supported positive regular Borel measures
$\lambda$ and $\mu$ on~$M$, the nonnegative numbers
$$
\langle \mu ,\lambda  \rangle \equiv \int _M \int _M G(x,y) \,
d\mu (x)d\lambda (y) =\int _M G_{\mu }\, d\lambda \quad \text
{and} \quad \| \mu \|^2\equiv \langle \mu ,\mu \rangle
$$
are called the {\it mutual energy} of~$\mu $ and $\lambda $ and
the {\it energy} of~$\mu $, respectively.
\begin{lem}[Energy principle]\label{Appendix Energy Principle Lemma} For all compactly
supported positive regular Borel measures $\mu $ and $\lambda $,
we have
$$
\langle \mu ,\lambda  \rangle \leq \| \mu \| \cdot \| \lambda \| .
$$
\end{lem}
\begin{pf} Recall that $M$ admits a {\it minimal heat kernel};
that is, a positive $\cinf $ function
$$
P: M\times M \times (0,+\infty ) \to (0,+\infty )
$$
\st for all $x,y\in M$,
$$
P(x,y,\cdot )=P(y,x,\cdot ) \quad \text {and} \quad G(x,y)=\int
_0^{\infty } P(x,y,t)\, dt .
$$
Moreover, for all $t,s>0$,
$$
P(x,y,t+s)=\int _M P(x,z,t)P(z,y,s)\, dV(z)
$$
(see \cite {Chavel-Book}). Therefore
\begin{align*} \langle \mu ,\lambda
\rangle
&=\int _M \int _M G(x,y)\, d\mu (x)d\lambda (y) \\
&=\int _0^{\infty } \int _M \int _M P(x,y,t)\, d\mu (x)d\lambda (y) dt \\
&=\int _0^{\infty } \int _M \biggl[ \biggl( \int _M P(x,z,t/2)\,
d\mu (x) \biggr) \biggl( \int _M P(y,z,t/2)\, d\lambda (y) \biggr)
\biggr] \, dV(z) dt.
\end{align*}
In particular,
$$
\| \mu \| ^2 =\int _0^{\infty } \int _M \biggl( \int _M
P(x,z,t/2)\, d\mu (x) \biggr)^2 \, dV(z) dt
$$
and similarly for $\| \lambda \| ^2$. Therefore, by the Schwarz
inequality,
\begin{align*} \langle \mu ,\lambda  \rangle &\leq
\biggl[ \int _0^{\infty } \int _M \biggl( \int _M P(x,z,t/2)\,
d\mu (x) \biggr) ^2 \, dV(z)dt
\biggr] ^{1/2} \\
&\qquad\qquad \cdot \biggl[ \int _0^{\infty } \int _M \biggl( \int
_M P(y,z,t/2)\, d\lambda (y) \biggr) ^2 \, dV(z) dt
\biggr] ^{1/2}\\
&=\| \mu \| \cdot \| \lambda \| .
\end{align*}
\end{pf}
\begin{rmk} For a different proof see \cite{Maeda}.
\end{rmk}

\subsection*{C. Evans-Selberg potential}

In this subsection, we recall Nakai's construction of a harmonic
\exh function on the closure of an end of a parabolic Riemannian
manifold \cite{Nakai1} and \cite{Nakai2} (see also \cite{SaNo}).
We begin with two useful facts:

\begin{lem}\label{Appendix H(x,y) harmonic in one var is cont lem} Suppose $M$ is an oriented Riemannian
manifold, $N$ is a topological space, and $H(x,y)$ is a positive
\fn on $M\times N$ which is harmonic in~$x$ and \cont in~$y$.
Then~$H$ is \cont on $M\times N$.
\end{lem}
\begin{pf} Let $(x_0,y_0)$ be a point in $M\times N$. By the
Harnack inequality, there exists a \cont \fn $\delta : M\times M
\to [0,+\infty )$ such that $\delta $ vanishes on the diagonal and
$$
(1+\delta (x,x_0))\inv H(x_0,y)\leq H(x,y) \leq (1+\delta
(x,x_0))H(x_0,y)
$$
for all points $x\in M$ sufficiently close to~$x_0$ and all points
$y\in N$. Hence
$$
|H(x,y)-H(x_0,y)|\leq \delta (x,x_0)(1+\delta (x,x_0))H(x_0,y).
$$
Since the function $H(x_0,\cdot )$ is \cont on~$N$, $H(x_0,y) \to
H(x_0,y_0)$ as $y \to y_0$. It now follows easily that
$$
|H(x,y)-H(x_0,y_0)| \leq |H(x,y)- H(x_0,y)|+|H(x_0,y)-H(x_0,y_0)|
\to 0
$$
as $(x,y) \to (x_0,y_0)$. Thus $H$ is \cont at~$(x_0,y_0)$.
\end{pf}

\begin{lem}\label{Appendix H(x,y) harmonic wrt product lem} Suppose $M$~and~$N$ are two oriented non\cpt
\con Riemannian manifolds and $H(x,y)$ is a positive \fn on
$M\times N$ which is harmonic in each variable. Then~$H$ is
Harmonic with respect to the product metric on $M\times N$. In
particular, $H$ is of class~$\cinf $.
\end{lem}
\begin{pf} By the previous lemma, $H$ is continuous. Moreover,
given a smooth compactly supported function $\alpha $ on~$M\times
N$, we have
\begin{align*} \int _{M\times N}H(x,y)\lap \alpha
(x,y) \, dV(x,y) &=\int _{M} \biggl[ \int _{N}H(x,y)\lap \alpha
_x(y) \, dV(y) \biggr]
dV(x) \\
&\qquad + \int _{N} \biggl[ \int _{M}H(x,y)\lap \alpha _y(x) \,
dV(x) \biggr] dV(y) =0;
\end{align*}
where, in each integral, $\lap $ and $dV$ denote the appropriate
Laplacian and volume element, and, for each $x\in M$ and $y\in N$,
$\alpha _x\equiv \alpha (x,\cdot )$ and $\alpha _y\equiv \alpha
(\cdot ,y)$. Thus $H$ is harmonic on $M\times N$.
\end{pf}

{\it Stone-\v Cech compactification}. Let $X$ be a topological
space. We will call a \cont map $f :X \to [-\infty , +\infty ]$ a
{\it continuous function} on~$X$ and we will denote the space of
\cont functions by~$\cal C(X)$. If $f(X)\subset (-\infty , +\infty
)$, then we will call~$f$ a {\it finitely \cont function}. We will
denote the space of finitely \cont functions by $C^0(X)$.

If~$X$ is locally \cpt Hausdorff, then the {\it Stone-\v Cech
compactification} $\check X$ is the unique Hausdorff
compactification of~$X$ to which every \cont \fn on~$X$ extends
continuously. The set $\Gamma =\check X\setminus X$ is the {\it
Stone-\v Cech boundary} of~$X$.

{\it Green's function on the Stone-\v Cech compactification}. For
the rest of this subsection, we will assume that $(M,g)$ is {\it
parabolic}. $\check M$ will denote the Stone-\v Cech
compactification of~$M$, $\Gamma $ will denote the Stone-\v Cech
boundary of~$M$, and $M_0$ will denote a fixed (nonempty) $\cinf $
\rel \cpt domain in~$M$ {\it with connected boundary}.

In particular, the manifold $M\setminus \overline M_0$ is
connected and hyperbolic. For if $M_1$ and $\Omega $ are~$\cinf $
domains in~$M$ with
$$
M_0 \Subset M_1 \Subset M \quad \text {and} \quad \Omega \Subset
M\setminus \overline M_1,
$$
\(u\) is the \hmibns ~$M\setminus \overline M_0$ \wrt the
complement of $\overline \Omega $, and $v$ is the \cont function
on~$\overline M_1\setminus M_0$ which is \harm on~$M_1\setminus
\overline M_0$, equal to~$0$ on~$\partial M_1$, and equal to~$1$
on~$\partial M_0$, then $v\leq u\leq 1$ on~$M_1\setminus \overline
M_0\subset (M\setminus \overline M_0)\setminus \overline \Omega $.
Hence $u\not \equiv 0$.

We will denote the Green's function on~$M\sm\overline{M_0} $
by~$G(x,y)$. In the above, we may extend~$u$ continuously
to~$\partial M_0$ by the constant~$1$. Thus every sequence in
$M\setminus \overline M_0$ approaching~$\partial M_0$ is a regular
sequence. Hence if we extend~$G$ to a function on $(M\sm
M_0)\times (M\sm M_0)$ by setting
$$
G(x,y) =0 \quad \text {if } x\in \partial M_0 \text { or } y\in
\partial M_0,
$$
then, for each point $x_0\in  M\sm M_0$, the function $G_{x_0}$ is
\cont on $M\setminus M_0$. In fact, by Lemma~\ref{Appendix
potential gives borel measure lemma}, $G$ is \cont on the set
$$
[(M\setminus M_0)\times (M\setminus M_0)]\setminus (\partial M_0
\times \partial M_0 ).
$$

Nakai's main observation is that the Green's function~$G(x,y)$ may
be continuously extended, in each variable, to the Stone-\v Cech
boundary. More precisely, we have the following:

\begin{prop}[Nakai]\label{Appendix Green's fn extension lemma} The Green's function
$G$ on $M\sm\overline{M_0} $ extends to a function
$$
\check G :(\check M\sm M_0)\times (\check M\sm M_0) \to [0,+\infty
]
$$
given by the double limit
$$
\check G (x_0,y_0)\equiv \lim _{x \to x_0} \biggl(\lim _{y \to
y_0} G(x,y)\biggr) \quad \forall \, x_0,y_0 \in \check M\sm M_0,
$$
where, in the above limits, $x,y\in M\sm\overline{M_0}$. Moreover,
this function has the following properties:
\begin{enumerate} \item[(i)] $\check G(x,y)=\check G(y,x) \quad
\forall \, (x,y) \in (M\sm M_0)\times (\check M\sm M_0) $;

\item[(ii)] $\check G$ is \cont on $[(M\sm M_0) \times (\check
M\sm M_0)]\setminus (\partial M_0 \times
\partial M_0)$ and finitely \cont on $(M\sm M_0) \times \Gamma $;

\item[(iii)] For each point $y\in \check M\sm M_0$, the function
$\check G _y \equiv \check G (\cdot ,y)$ is \cont on $\check M\sm
M_0$, \harm on $(M\sm\overline{M_0} )\setminus \{ y \} $, and
equal to~$0$ on~$\partial M_0$;

\item[(iv)]  For each point $y\in \check M\sm\overline{M_0} $
($y\in
\partial M_0$), $\check G_y >0$ ($\, \equiv 0$) on $\check
M\sm\overline{M_0} $; and

\item[(v)] For each point $y\in \check M\sm\overline{M_0} $,
$$
\int _{\partial M_0} \ndof {\check G_y} (x) \, d\sigma (x)
=(n-2)\tau _{n-1}.
$$
\end{enumerate}
\end{prop}
\begin{rmk} It is not clear that $\check G$ is symmetric
on $(\check M\sm M_0) \times (\check M\sm M_0) $ since, for
$x\in\Gamma$, it is not clear that the function $y\mapsto\check
G(x,y)$ is \cont on $\check M\sm M_0$.
\end{rmk}

\begin{pf*}{Proof of Proposition~\ref{Appendix Green's fn extension lemma}} For each point $x\in\set$, the
function
$$
G_x=G(x,\cdot )=G(\cdot ,x): M\sm M_0\to [0,+\infty ]
$$
is \cont  and hence extends to a \cont \fn on~$\check M\sm M_0$.
Thus the function $G^* (x,y)$ defined by
$$
G^* (x,y)=G^*_{y}(x)\equiv \lim _{z \to y} G(x,z) \quad \forall
(x,y)\in (M\sm M_0)\times (\check M\sm M_0) ,
$$
with $z\in M\sm\overline{M_0} $ is an extension of~$G$ which is
\cont in~$y$ for each fixed\linebreak $x\in M\sm M_0$. We show
that, for each fixed $y_0\in \Gamma $, the function~$G^* _{y_0}$
is \cont on~$M\sm M_0$, positive and \harm on $M\sm\overline{M_0}
$, and identically equal to~$0$ on~$\partial M_0$.

Suppose $\Omega _1$ and $\Omega _2$ are domains with
$$
\Omega _1 \Subset \Omega _2 \Subset M\sm\overline{M_0} .
$$
Then, by the Harnack inequality (Lemma~1.2), there exists a \cont
\fn $\delta : \Omega _2 \times \Omega _2  \to [0,+\infty ) $ \st
$\delta $ vanishes on the diagonal and, for every positive \harm
\fnns~$u$ on $\Omega _2$,
$$
|u(x) -u(x_0)| \leq\delta (x,x_0) \max (u(x),u(x_0)) \quad \forall
\, x,x_0\in \Omega _1.
$$
Next, observe that
$$
0\leq G(x,y) \leq a \equiv \max _{\overline \Omega _1 \times
\partial \Omega _2}G < +\infty \quad \forall \, (x,y) \in
\overline \Omega _1 \times ((M\sm M_0) \setminus \Omega _2).
$$
For if $\Omega $ is a~$\cinf $ domain with
$$
\Omega _2 \Subset \Omega \Subset M\sm\overline{M_0} ,
$$
and $G'$ is the Green's function of~$\Omega $, then, for each
point $x\in \overline \Omega _1$, $G'_x \restriction _{(\Omega
\setminus \overline \Omega _1)} $ is a positive \harm \fn which
vanishes continuously at~$\partial \Omega $. Hence
$$
G'_x \restriction _{(\Omega \setminus \overline \Omega _2)} \leq
\max _{\partial \Omega _2} G_x' \leq \max _{\overline \Omega
_1\times \partial \Omega _2} G' .
$$
Since $G$ is the pointwise limit of an increasing limit of such
Green's functions, the inequality follows. It also follows that
the collection of positive {\it finitely} \cont functions $\cal G
\equiv \{ \, G_y \restriction _{\overline \Omega _1} \mid y\in
(M\sm M_0) \setminus \Omega _2 \, \} $ is precompact
in~$C^0(\overline \Omega _1)$. For if $y\in (M\sm M_0) \setminus
\Omega _2$, then $G_y$ is bounded by~$a$ and \harm on~$\Omega _2$.
Hence
$$
|G_y(x)-G_y(x_0)|\leq a\delta (x,x_0) \quad \forall \, x,x_0 \in
\Omega _1.
$$
Since $\delta (x,x_0) \to 0$ as $x \to x_0$, $\cal G$ is bounded
and equi\contns , hence precompact, by Ascoli's Theorem.

Now let $\seq y \alpha $ be a net in~$(M\sm\overline{M_0} )
\setminus \overline\Omega _2$ which converges to the given
point~$y_0\in \Gamma $. Then, by passing to a subnet, we may
assume that $\{ \bigl( G_{y_{\alpha }} \bigr) \restriction
_{\overline \Omega _1} \} $ converges uniformly on~$\overline
\Omega _1$ to the function $G^*_{y_0}\restriction _{\overline
\Omega _1}$. In particular, $G^* _{y_0}$ is nonnegative and \harm
on~$\Omega _1$.

Finally, for a fixed point~$x_0\in \Omega _1$, the \fn $G_{x_0}$
is \contns , positive, and  superharmonic on~$M\sm\overline{M_0}
$. Hence, if $\Omega _0$ is a~$\cinf $ domain with
$$
M_0 \Subset \Omega _0 \Subset M\setminus \overline \Omega _1,
$$
then, by the minimum principle for parabolic manifolds,
$$
\inf _{M\setminus \Omega _0} G_{x_0} =\min _{\partial \Omega
_0}G_{x_0} \geq \min _{\overline \Omega _1 \times \partial \Omega
_0} G(x,y) \equiv b >0.
$$
Therefore
$$
G^*_{y_0}(x_0)=\lim _{\alpha } G_{y_\alpha }(x_0) =\lim _{\alpha }
G _{x_0} (y_\alpha ) \geq b >0.
$$
Therefore $G^*_{y_0}>0$ on~$\Omega _1$. Since the domain $\Omega
_1 \Subset M\sm\overline{M_0} $ is arbitrary, the function
$G^*_{y_0}$ is positive and \harm on~$M\sm\overline{M_0} $.

Next we verify that $G^*_{y_0}$ is continuous at~$\partial M_0$.
It is clear from the definition that~$G^*_{y_0}$ vanishes
on~$\partial M_0$. Let $\Omega $ be a~$\cinf $ domain with
$$
M_0 \Subset \Omega \Subset M.
$$
Then $\partial \Omega $ is a \cpt subset of~$M\sm\overline{M_0} $.
By the above discussion, every net $\{ y_\alpha \} _{\alpha \in
A}$ in $M\sm\overline{M_0} $ converging to the given point~$y_0\in
\Gamma $ admits a subnet~$\{ y_\beta \} _{\beta  \in B}$ in
$M\setminus \overline \Omega $ \st $\seq {G}{y_\beta }$ converges
uniformly on~$\partial \Omega $ to $G^*_{y_0}$. For all $\beta
,\gamma \in B$ and $x\in \overline \Omega \setminus M_0$, we have
$$
\bigl| G_{y_\beta }(x)-G_{y_\gamma }(x) \bigr| \leq \max
_{\overline \Omega \setminus M_0} \bigl| G_{y_\beta }-G_{y_\gamma
} \bigr| \leq \max _{\partial \Omega } \bigl| G_{y_\beta
}-G_{y_\gamma } \bigr| ,
$$
because $G_{y_\beta }$ and $G_{y_\gamma }$ are \harm on~$\Omega
\setminus \overline M_0$ and equal to~$0$ on~$\partial M_0$. Since
$\seq {G}{y_\gamma }$ converges to~$G^*_{y_0}$ pointwise
on~$\overline \Omega \setminus M_0$ and uniformly on~$\partial
\Omega $, we get
$$
\bigl| G_{y_\beta }(x)-G^*_{y_0}(x) \bigr| \leq \max _{\partial
\Omega } \bigl| G_{y_\beta }-G^*_{y_0} \bigr| .
$$
Hence
$$
\sup _{\overline \Omega \setminus M_0} \bigl| G_{y_\beta
}-G^*_{y_0}\bigr| \leq \max _{\partial \Omega } \bigl| G_{y_\beta
}-G^*_{y_0} \bigr| \to 0.
$$
Therefore the \cont functions~$\seq {G}{y_\beta }$ converge
uniformly on~$\overline \Omega \setminus M_0$ to~$G^*_{y_0}$ and
hence $G^*_{y_0}$ vanishes continuously at~$\partial M_0$.

Thus, for each point $y_0\in \check M\sm M_0$, the function
$G^*_{y_0}$ is \cont on~$M\setminus M_0$. Hence $G^*_{y_0}$
extends to a \cont \fn $\check G_{y_0}$ on~$\check M\sm M_0$. We
may therefore define $\check G (x_0,y_0)$ for each pair of points
$x_0,y_0 \in \check M\sm M_0$ by
$$
\check G (x_0,y_0)\equiv \check G_{y_0}(x_0)\equiv \lim _{x \to
x_0} \biggl(\lim _{y \to y_0} G(x,y)\biggr) =\lim _{x \to x_0}
\biggl( G^*_{y_0}(x)\biggr) ;
$$
where, in the above limits, $x,y\in M\sm\overline{M_0} $.

We now verify that $\check G$ has the properties (i)--(v). Given
points $x\in M\sm M_0$ and $y\in \check M\sm M_0$, we have
\begin{align*} \check G(y,x) &= \lim _{z \to y} \biggl(\lim _{w
\to x} G(z,w)\biggr) =\lim _{z \to y} G(z,x)
=\lim _{z \to y} G(x,z)\\
&=G^*_y(x) =\check G_y\restriction _{M\sm M_0} (x) =\check G(x,y).
\end{align*}
with $w,z\in M\sm\overline{M_0} $. Thus (i) is proved.

The properties (iii) and (iv) follow by construction.

By the above discussion, given a cpt set~$K$ in~$M\sm M_0$ and a
point~$y_0\in \Gamma $, every net $\seq y\alpha $
in~$M\sm\overline{M_0} $ converging to~$y_0$ admits a subnet~$\seq
y\beta $ \st $\seq {G}{y_\beta }$ converges to~$G^*_{y_0}$
uniformly on~$K$. It follows that $G_y \to G^*_{y_0}$ uniformly
on~$K$ as $y\in M\sm\overline{M_0} $ approaches~$y_0$. Moreover,
we have shown that the functions $G_y$ are uniformly bounded
on~$K$ for $y\in M\sm\overline{M_0} $ near~$y_0$. Hence $\check
G_{y_0} =G^*_{y_0}$ is bounded on~$K$.

Since $\check G=G$ on $((M\sm M_0)\times (M\sm M_0) )\setminus
(\partial M_0 \times \partial M_0)$, in order to prove~(ii), it
remains to show that $\check G$ is finitely \cont at each point
$(x_0,y_0)\in (M\sm M_0)\times \Gamma $. By the above remarks,
$\check G(x_0,y_0)<+\infty $ and for points $x\in M\sm M_0$
near~$x_0$ and $y\in M\sm\overline{M_0} $ near~$y_0$, we have
\begin{align*}
\bigl| \check G(x,y) -\check G(x_0,y_0) \bigr| &\leq \bigl|
G(x,y)-\check G(x,y_0) \bigr|
+\bigl| \check G(x,y_0) -\check G(x_0,y_0) \bigr| \\
&=\bigl| G_y(x) -G^*_{y_0}(x) \bigr| +\bigl| G^*_{y_0}(x)
-G^*_{y_0}(x_0) \bigr| \to 0
\end{align*}
as $(x,y) \to (x_0,y_0)$. Thus we need only consider points in
$(M\sm M_0) \times \Gamma $ approaching~$(x_0,y_0)$. Moreover, the
term $\check G(x,y_0) -\check G(x_0,y_0) $ is independent of~$y$
and approaches~$0$ as~$x$ approaches~$x_0$. Thus it remains to
show that
$$
\bigl| \check G(x,y) -\check G(x,y_0)\bigr| \to 0 \quad \text {as}
\quad (x,y) \to (x_0,y_0) \quad \text {with} \quad (x,y) \in (M\sm
M_0) \times \Gamma .
$$
If $z$ is a point in~$M\sm\overline{M_0} $, then
$$
\bigl| \check G(x,y) -\check G(x,y_0)\bigr| \leq \bigl| \check
G(x,y) -\check G(x,z)\bigr| +\bigl| \check G(x,z) -\check
G(x,y_0)\bigr| .
$$
Since $\check G_z =G_z$ converges uniformly to $G^*_{y_0}$ on a
\rel \cpt \nbd of~$x_0$ as $z\in M\sm\overline{M_0} $
approaches~$y_0$, given a positive real number~$\epsilon $ there
exists a \nbdns~$U$ of~$y_0$ in~$\check M$ such that
$$
\bigl| \check G(x,z) -\check G(x,y_0)\bigr| <\frac {\epsilon }{2}
\quad \forall \, z\in (M\sm\overline{M_0} )\cap U.
$$
Hence, given a point $y\in U\cap \Gamma $ and a point $z\in
(M\sm\overline{M_0} )\cap U$ so close to~$y$ that $\bigl| \check
G(x,y) -\check G(x,z)\bigr| < \epsilon  /2$ for all points $x\in
M\sm M_0$ near~$x_0$, we get
$$
\bigl| \check G(x,y) -\check G(x,y_0)\bigr| <\epsilon .
$$
Thus the claim is proved.

For the proof of~(v), let $y_0\in M\sm\overline{M_0}$, let $\seq
Mk$ be an \exh of~$M$ by $\cinf $ \rel \cpt domains containing
$\overline M_0\cup \{ y_0 \} $, let $\Omega $ be a $\cinf $ \rel
\cpt domain in $M_1\setminus \overline M_0$ containing~$y_0$, and,
for each positive integer~$k$, let $G_k$ be the Green's function
on~$M_k\setminus \overline M_0$ and let~$v_k$ be the \cont \fn on
$\overline M_k\setminus (M_0 \cup \Omega )$ which is \harm
on~$M_k\setminus (\overline M_0 \cup \overline \Omega )$, equal
to~$1$ on~$\partial M_0\cup \partial \Omega $, and equal to~$0$
on~$\partial M_k$. Then the function~$1-v_k$ is nonnegative and
not greater than the restriction of the \hm of~$\partial M_k$ \wrt
$M_k\setminus \overline M_0$. Therefore, since~$M$ is parabolic,
$v_k \to 1$ uniformly on \cpt sets in~$M\setminus (M_0 \cup \Omega
)$. Consequently, $\grad v_k \to 0$ uniformly on \cpt subsets
of~$M\setminus (\overline M_0 \cup \overline \Omega )$.

The function $G_{y_0}-(G_k)_{y_0}$ is positive and \cont
on~$\overline M_k\setminus M_0$, harmonic on~$M_k\setminus
\overline M_0$, and equal to~$0$ on~$\partial M_0$. Moreover,
$G_{y_0}-(G_k)_{y_0} \to 0$ uniformly on \cpt subsets of~$M\sm
M_0$. For $\{ (G_k)_{y_0} \} $ dominates a sequence of Greens
functions on \rel \cpt domains exhausting~$M\sm\overline{M_0} $
and hence $G_{y_0}-(G_k)_{y_0} \to 0$ uniformly on \cpt subsets
of~$M\sm\overline{M_0} $.  By applying the maximum principle and
the fact that $G_{y_0}$ and the functions $\{ (G_k)_{y_0}\} $
vanish on~$\partial M_0$, we get uniform convergence
near~$\partial M_0$ as well. In particular, $\grad G_{y_0}-\grad
(G_k)_{y_0} \to 0$ uniformly on \cpt subsets
of~$M\sm\overline{M_0} $. Therefore
\begin{align*} \dir
_{M_k\setminus (\overline M_0\cup \overline \Omega )}
((G_k)_{y_0}, v_k) &= \int _{M_k\setminus (\overline M_0\cup
\overline \Omega )}
\langle \grad (G_k)_{y_0}, \grad v_k \rangle \\
&=\int _{\partial (M_k\setminus (\overline M_0\cup \overline
\Omega ))}
\ndof {(G_k)_{y_0}}v_k \\
&=\int _{\partial M_k} \ndof {(G_k)_{y_0}}v_k -\int _{\partial
M_0} \ndof {(G_k)_{y_0}}v_k -\int _{\partial \Omega }
\ndof {(G_k)_{y_0}}v_k \\
&= -\int _{\partial M_0} \ndof {(G_k)_{y_0}} -\int _{\partial
\Omega }
\ndof {(G_k)_{y_0}} \\
&= -\int _{\partial M_0'} \ndof {(G_k)_{y_0}} -\int _{\partial
\Omega } \ndof {(G_k)_{y_0}}
\end{align*}
where $M_0'$ is any $\cinf $ \rel \cpt domain in~$M_1 \setminus
\overline \Omega $ containing~$\overline M_0$. Hence
\begin{align*} \dir _{M_k\setminus (\overline M_0\cup \overline
\Omega )} ((G_k)_{y_0}, v_k) & \to
 -\int _{\partial M_0'}
\ndof {G_{y_0}} -\int _{\partial \Omega } \ndof {G_{y_0}} = -\int
_{\partial M_0} \ndof {G_{y_0}} -\int _{\partial \Omega } \ndof
{G_{y_0}}
\end{align*}
as $k \to \infty $. On the other hand, we have
\begin{align*} \dir
_{M_k\setminus (\overline M_0\cup \overline \Omega )}((G_k)_{y_0},
v_k) &=\int _{\partial (M_k\setminus (\overline M_0\cup \overline
\Omega ))}
(G_k)_{y_0}\ndof {v_k} \\
&=\int _{\partial M_k} (G_k)_{y_0}\ndof {v_k} -\int _{\partial
M_0} (G_k)_{y_0}\ndof {v_k} -\int _{\partial \Omega }
(G_k)_{y_0}\ndof {v_k} \\
&= -\int _{\partial \Omega }
(G_k)_{y_0}\ndof {v_k} \\
&= -\int _{\partial \Omega '} (G_k)_{y_0}\ndof {v_k} +\int
_{\partial \Omega '} \ndof {(G_k)_{y_0}}v_k - \int _{\partial
\Omega }
\ndof {(G_k)_{y_0}}v_k \\
&\to \int _{\partial \Omega '} \ndof {G_{y_0}} - \int _{\partial
\Omega } \ndof {G_{y_0}}=0 \quad \text {as} \quad k \to \infty ;
\end{align*}
where $\Omega '$ is any $\cinf $ \rel \cpt domain in~$M_1
\setminus \overline M_0$ containing~$\overline \Omega $. Hence if
$\widehat G $ is the Green's function on~$\Omega $, then
\begin{align*} \int _{\partial M_0} \ndof {G_{y_0}} &=-\int
_{\partial \Omega } \ndof {G_{y_0}} =-\int _{\partial \Omega }
\ndof {\widehat G _{y_0}} + \int _{\partial \Omega }
\ndof {(G _{y_0}-\widehat G _{y_0})} \\
&=-\int _{\partial \Omega } \ndof {\widehat G _{y_0}} \cdot 1
=(n-2)\tau _{n-1},
\end{align*}
since $G _{y_0}-\widehat G _{y_0}$ is a $\cinf $ function
on~$\overline \Omega $ which is \harm on the interior.

Finally, suppose $y_0\in \Gamma $ and let $M_1$ be a $\cinf $ \rel
\cpt domain in~$M$ containing~$\overline M_0$. Then, as $y \to
y_0$ with $y\in M\setminus \overline M_1$, $G_y \to \check
G_{y_0}$ uniformly on~$\overline M_1\setminus M_0$; and therefore
$$
(n-2)\tau _{n-1} =\int _{\partial M_0} \ndof {\check G_{y}} =\int
_{\partial M_1} \ndof {\check G_{y}} \to \int _{\partial M_1}
\ndof {\check G_{y_0}} =\int _{\partial M_0} \ndof {\check
G_{y_0}} .
$$
The claim (v) follows.
\end{pf*}

\noindent {\it Transfinite Diameter and Tchebycheff constant}. Let
$K$ be a \cpt subset of $\check M\sm\overline{M_0} $. Given an
integer~$m>1$ and points $x_1,\dots ,x_m\in K$, let $D_m(K,x_1,
\dots ,x_m)=D_m(x_1, \dots ,x_m)$ and $E_m(K,x_1, \dots ,x_m)$ be
the following numbers:
$$
D_m(x_1, \dots ,x_m) \equiv {\binom m2}\inv  \sum _{1\leq i<j\leq
m} \check G(x_i,x_j)
$$
and
$$
E_m(K,x_1, \dots ,x_m) \equiv \frac 1m \biggl( \inf _{x\in K} \sum
_{i=1}^m \check G (x,x_i) \biggr) .
$$
We then define $D_m(K)$ and $E_m(K)$ by
$$
D_m(K)\equiv \inf _{x_1,\dots ,x_m \in K} D_m(x_1, \dots ,x_m)
$$
and
$$
E_m(K)\equiv \sup _{x_1,\dots ,x_m \in K} E_m(K,x_1, \dots ,x_m).
$$
The numbers
$$
D(K) \equiv \sup _{m>0} D_m(K) \quad \text {and} \quad E(K) \equiv
\sup _{m>0} E_m(K)
$$
are called the {\it transfinite diameter} of~$K$ and the {\it
Tchebycheff constant} of~$K$, respectively.

\begin{lem}\label{Appendix D(K) and E(K) lemma} For every \cpt subset~$K$ of $\check
M\sm\overline{M_0}$,
\begin{enumerate}
\item[(i)] $D_m(K)
\nearrow D(K)$ and $E_m(K) \to E(K)$ as $m \to \infty $; and

\item[(ii)] $0\leq D(K) \leq E(K) \leq +\infty $.
\end{enumerate}
\end{lem}
\begin{pf} Given a positive integer~$m$ and points~$x_1,\dots,x_{m+1}\in K$,
we have, for each~$k=1, \dots ,m+1$,
\begin{align*} \binom {m+1}{2} D_{m+1}(x_1,\dots ,x_{m+1})
&=\sum _{1\leq i<k}\check G(x_i,x_j)+\sum _{k<j\leq m+1}\check G(x_i,x_j) \\
&\qquad \qquad + \sum _{1\leq i<j\leq m+1;\, i,j\neq k}\check G(x_i,x_j) \\
&=\sum _{1\leq i<k}\check G(x_i,x_j)+\sum _{k<j\leq m+1}\check G(x_i,x_j) \\
&\qquad \qquad +\binom {m}{2} D_m(x_1,\dots ,\hat x_k,\dots ,x_{m+1}) \\
&\geq \sum _{1\leq i<k}\check G(x_i,x_j)+\sum _{k<j\leq m+1}\check G(x_i,x_j)\\
&\qquad \qquad +\binom {m}{2} D_m(K),
\end{align*}
where $(x_1,\dots ,\hat x_k,\dots ,x_{m+1})$ denotes the $m$-tuple
obtained by removing the $k^{\text {th}}$ term from the
$(m+1)$-tuple~$(x_1,\dots ,x_{m+1})$. Summing over~$k$, we get
\begin{align*} (m+1)\binom {m+1}{2} D_{m+1}(x_1,\dots ,x_{m+1})
\geq
2\binom {m+1}{2} &D_{m+1}(x_1,\dots ,x_{m+1})\\
&+(m+1)\binom {m}{2} D_m(K).
\end{align*}
It follows that $D_{m+1}(x_1,\dots ,x_{m+1})\geq D_m(K)$.
Therefore
$$
D_m(K)\leq D_{m+1}(K) \to D(K) \quad \text {as} \quad m \to
\infty.
$$
Given positive integers~$m$~and~$l$ and
points~$x_1,\dots,x_{m+l}\in K$, we have
\begin{align*}
(m+l)E_{m+l}(K,x_1,\dots ,x_{m+l}) &=\inf _{x\in K} \biggl( \sum
_{i=1}^m \check G (x,x_i) +\sum _{i=m+1}^{m+l}\check G (x,x_i)
\biggr) \\
&\geq \inf _{x\in K} \biggl( \sum _{i=1}^m \check G (x,x_i)
\biggr) + \inf _{x\in K} \biggl( \sum _{i=m+1}^{m+l}\check G
(x,x_i)
\biggr) \\
&=mE_{m}(K,x_1,\dots ,x_{m})+lE_{l}(K,x_{m+1},\dots ,x_{m+l}).
\end{align*}
Therefore $(m+l)E_{m+l}(K)\geq mE_{m}(K)+lE_{l}(K)$. Let~$q$
and~$r$ be integers satisfying
$$
m=ql+r \quad \text {and} \quad 0\leq r<l.
$$
By the above, we have $qlE_{ql}(K)\geq
(q-1)lE_{(q-1)l}(K)+lE_l(K)$. Proceeding inductively, we get
$qlE_{ql}(K)\geq qlE_l(K)$. Therefore
$$
mE_m(K) =(ql+r)E_{ql+r}(K)\geq qlE_{ql}(K)+rE_{r}(K) \geq
qlE_{l}(K).
$$
Hence
$$
E(K)\geq E_m(K) \geq \biggl( \frac {ql}{ql+r} \biggr) E_l(K).
$$
For fixed~$l$, we have $q \to +\infty $ as $m \to +\infty $
and~$0\leq r<l$. Thus
$$
E(K)=\sup_{m>0}E_m(K)\geq \limsup_{m \to\infty}E_m(K) \geq
\liminf_{m \to \infty }E_m(K) \geq E_l(K)
$$
for every positive integer~$l$ and hence
$$
E(K)= \lim _{m \to \infty }E_m(K).
$$
Thus~(i) is proved.

Given a positive integer~$m$, there exist points $x_1,\dots ,x_m
\in K$ \st
$$
(1)\qquad\qquad\qquad\qquad iE_i(K,x_{m-i+1}, \dots ,x_m)= \sum
_{j=m-i+1}^m \check
G(x_{m-i},x_j)\qquad\qquad\qquad\qquad\qquad\qquad
$$
for $i=1, \dots ,m-1$. To see this, let $x_m\in K$ be arbitrary,
let $1\leq k < m$, and suppose  $x_{m-k+1}, \dots ,x_m\in K$ have
been chosen so that the equality~(1) holds for $i=1,\dots ,k-1$
(the case~$k=1$ is vacuously true). By Proposition~\ref{Appendix
Green's fn extension lemma}~(iii), the function
$$
x \mapsto \sum _{j=m-k+1}^m \check G(x,x_j)
$$
is \cont and therefore assumes its minimum value on~$K$ at some
point~$x_{m-k}=x_{m-(k+1)+1}$. The equality~(1) then holds
for~$i=k$. Thus, proceeding inductively, we obtain points
$x_1,\dots ,x_m \in K$ with the required properties.

From~(1) it follows that
$$
iE_i(K)\geq \sum _{j=m-i+1}^m \check G(x_{m-i},x_j)
$$
for $i=1, \dots ,m-1$. Summing over~$i$, we get
\begin{align*}
\sum _{i=1}^{m-1}iE_i(K) &\geq \sum _{i=1}^{m-1}\sum _{j=m-i+1}^m
\check G(x_{m-i},x_j)
=\sum _{1\leq i<j\leq m}\check G(x_i,x_j) \\
&=\binom m2 D_m(x_1,\dots ,x_m) \geq \binom m2 D_m(K).
\end{align*}
Clearly, we may assume that $E(K)<+\infty $. Therefore, for each
integer~$k$ with $0<k<m-1$,
\begin{align*} D_m(K)&\leq \binom m2
\inv \sum _{i=1}^{m-1}iE_i(K)
=E(K)+\binom m2 \inv \sum _{i=1}^{m-1}i(E_i(K)-E(K)) \\
&\leq E(K)+\binom m2 \inv \sum _{i=1}^{k}i|E_i(K)-E(K)| +\binom m2
\inv \sum _{i=k+1}^{m-1}i|E_i(K)-E(K)|.
\end{align*}
Hence
$$
D_m(K)\leq E(K) +\binom m2 \binom {k+1}2 \max _{1\leq i\leq
k}|E_i(K)-E(K)| +\sup _{i>k}|E_i(K)-E(K)|.
$$
By choosing $k$ sufficiently large, we can make the third term on
the right-hand side of the above equation arbitrarily small.
Moreover, for $k$ fixed, the second term approaches~$0$ as~$m \to
\infty $. The claim~(ii) now follows.
\end{pf}

The main step in Nakai's construction of the Evans-Selberg
potential is a proof that the transfinite diameter of the \v Cech
boundary is infinite. For this, the first step is the following
lemma:
\begin{lem}\label{Appendix D of complement is D of boundary lemma}
If $\Omega $ is a~$\cinf $ domain with $M_0\Subset \Omega \Subset
M$, then
$$
D(\check M\setminus \Omega )=D(\partial \Omega ).
$$
\end{lem}
\begin{pf} It suffices to show that, for every integer~$m>2$,
$$
D_m(\check M\setminus \Omega )=D_m(\partial \Omega );
$$
that is, given points $x_1, \dots ,x_m\in \check M\setminus \Omega
$ we have
$$
D_m(x_1, \dots ,x_m) \geq \inf _{y_1, \dots ,y_m \in \partial
\Omega } D_m(y_1,\dots ,y_m).
$$
We will prove by induction on~$k=0,1,2,\dots ,m$ that there exist
points $y_1,\dots ,y_k $ in $\partial \Omega $ \st
$$
D_m(x_1, \dots ,x_m) \geq D_m(y_1,\dots ,y_k , x_{k+1}, \dots
,x_m).
$$
The case~$k=0$ is clear. Suppose $0\leq k\leq m-1$ and there exist
points $y_1,\dots ,y_k$ in $\partial \Omega $ satisfying the above
inequality. By Proposition~\ref{Appendix Green's fn extension
lemma} (parts~(i)~and~(iii)), the function~$\psi $ defined by
\begin{align*}
\psi (y) &= \binom m2 D_m(y_1,\dots ,y_k , y, x_{k+2}, \dots ,x_m)\\
&=\binom {m-1}2 D_{m-1}(y_1,\dots ,y_k , x_{k+2}, \dots ,x_m)
+\sum _{i=1}^k \check G(y_i,y) +\sum _{i=k+2}^m \check G(y,x_i)
\end{align*}
is nonnegative and continuous on~$\check M\sm M_0\supset M
\setminus \Omega $ and super\harm on $M\sm\overline{M_0} \supset
M\setminus \overline \Omega $. Therefore, by the minimum principle
for parabolic manifolds (Sect.~1),
$$
\inf_{\check M\setminus \Omega } \psi =\inf _{M\setminus \Omega }
\psi =\min _{\partial \Omega }\psi =\psi (y_{k+1})
$$
for some point $y_{k+1}\in \partial \Omega $. Therefore
\begin{align*} \binom m2 D_m(x_1, \dots ,x_m) &\geq \binom m2
D_m(y_1,\dots ,y_k , x_{k+1}, \dots ,x_m)
=\psi (x_{k+1}) \\
&\geq \psi (y_{k+1}) =\binom m2 D_m(y_1,\dots ,y_k , y_{k+1},
x_{k+2},\dots ,x_m).
\end{align*}
Thus the claim and the lemma follow by induction.
\end{pf}

\noindent {\it Capacity}. Given a nonempty \cpt subset~$K$
of~$M\sm\overline{M_0} $, we denote by~$W(K)$ the positive number
$\inf \| \mu \| ^2$, where the infimum is over all unit positive
regular Borel measures supported in~$K$ and, for each such
measure~$\mu $,
$$
\| \mu \| ^2 =\int _MG_\mu (y) \, d\mu (y) = \int _M\int _M G(x,y)
\, d\mu (x)d\mu (y)
$$
is the energy of~$\mu $. Observe that
$$
W(K) \geq \min _{K\times K}G >0.
$$
The number $1/W(K)$ is called the {\it capacity} of~$K$.
\begin{lem}\label{Appendix-W of boundary lemma} Let~$\Omega $ be a~$\cinf $ domain
with $M_0\Subset \Omega \Subset M$. Then
$$
W(\partial \Omega )= \frac {(n-2)\tau _{n-1}}{\dir _{\Omega
\setminus \overline M_0}(u)},
$$
where $u$ is the \hm of~$\partial \Omega $ \wrt $\Omega \setminus
\overline M_0$.
\end{lem}
\begin{pf} The nonnegative finitely \cont function $\vphi $
on~$M\sm M_0$ defined by
$$
\vphi (x) =\left\{ \begin{aligned} u(x)
&\qquad \text {if } x\in \overline \Omega \setminus M_0\\
1 &\qquad \text {if } x\in M\setminus \Omega
\end{aligned}
\right.
$$
is a potential on~$M\sm\overline{M_0} $. For~$\vphi $ is
superharmonic and, if~$v$ is a \harm \fn on~$M\sm\overline{M_0} $
with $0\leq v\leq \vphi $, then $v$ vanishes continuously
at~$\partial M_0$. Hence $1-v$ is a nonnegative \cont \fn on~$M\sm
M_0$ which is super\harm on~$M\sm\overline{M_0} $ and equal to~$1$
on~$\partial M_0$. Therefore, by the minimum principle for a
superharmonic function on the complement of a~$\cinf $ domain in a
parabolic manifold, $1\leq 1-v \leq 1$ on~$M\sm M_0$. Therefore
$v\equiv 0$ and $\vphi $ is a potential. Moreover, $\vphi $ is
\harm on ~$(M\sm\overline{M_0} )\setminus \partial \Omega $.
Therefore, by Lemma~\ref{Appendix potential gives borel measure
lemma}, there exists a positive regular Borel measure~$\mu $
supported in~$\partial \Omega $ \st
$$
\vphi (x)=G_\mu (x) =\int _M G(x,y) \, d\mu (y) \quad \forall \,
x\in M\sm\overline{M_0} .
$$
Since $\vphi >0$ on~$M\sm\overline{M_0} $, $\mu (\partial \Omega
)>0$.

The infimum~$W(\partial \Omega )$ is attained by the energy $\|
\hat\mu \| ^2$ of the unit measure
$$
\hat\mu \equiv \mu /\mu (\partial \Omega ).
$$
For
$$
\| \mu \| ^2=\int _{\partial \Omega } G_\mu \, d\mu =\int
_{\partial \Omega } \vphi \, d\mu =\int _{\partial \Omega } \,
d\mu =\mu (\partial \Omega ).
$$
Hence
$$
\| \hat\mu \| ^2=\frac {\| \mu \| ^2}{(\mu (\partial \Omega ))^2}
=\| \mu \| ^{-2}.
$$
If $\lambda $ is a positive regular {\it unit} Borel measure~$\mu
$ supported in~$\partial \Omega $, then, by the energy principle
(Lemma~\ref{Appendix Energy Principle Lemma}),
$$
1=\int_{\partial\Omega}\,d\lambda =\int_{\partial\Omega}\vphi\,
d\lambda=\int _{\partial \Omega}G_\mu\,d\lambda
=\langle\mu,\lambda\rangle\leq\|\mu\|\cdot\|\lambda\|.
$$
Hence $\| \lambda  \| ^2 \geq \| \mu \| ^{-2}=\| \hat\mu \| ^2$.
Therefore
$$
W(\partial \Omega )=\| \hat\mu \| ^2 =1/\| \mu \| ^2=1/\mu
(\partial \Omega )
$$
as claimed.

On the other hand, since $u=0$ on~$\partial M_0$, $u=1$
on~$\partial \Omega $, and $u$ is \harm on~$\Omega \setminus
\overline M_0$, we have
\begin{align*} \dir _{\Omega \setminus
\overline M_0}(u) &=\int _{\Omega \setminus \overline M_0} |\grad
u|^2 =\int _{\partial (\Omega \setminus M_0)} u(x)\ndof u (x)\,
d\sigma (x)
= \int _{\partial \Omega } \ndof u (x)\, d\sigma (x) \\
&= \int _{\partial M_0} \ndof u (x)\, d\sigma (x) =\int _{\partial
M_0} \ndof \vphi (x)\, d\sigma (x)
=\int _{\partial M_0} \ndof {G_\mu }(x)\, d\sigma (x)\\
&=\int _{\partial M_0} \nd \biggl[ \int _{\partial \Omega
}G_y(x)\, d\mu (y) \biggr] \, d\sigma (x) =\int _{\partial M_0}
\int _{\partial \Omega }\ndof {G_y}(x)\, d\mu (y)
\, d\sigma (x) \\
&=\int _{\partial \Omega } \biggl[ \int _{\partial M_0} \ndof
{G_y}(x)\, d\sigma (x) \biggr] \, d\mu (y) =\int _{\partial \Omega
} \biggl[ \int _{\partial M_0} \ndof {\check G_y}(x)\, d\sigma (x)
\biggr]
\, d\mu (y) \\
&=(n-2)\tau _{n-1} \int _{\partial \Omega }\, d\mu (y)
\quad (\text {Proposition~\ref{Appendix Green's fn extension lemma}, part~(v)}) \\
&=(n-2)\tau _{n-1} \mu (\partial \Omega )
\end{align*}
The lemma now follows.
\end{pf}

\begin{lem}\label{Appendix-D bigger than W lemma} If $\Omega $ is a~$\cinf $ domain
with $M_0\Subset \Omega \Subset M$, then
$$
D(\partial \Omega ) \geq W(\partial \Omega ).
$$
\end{lem}
\begin{pf} For each integer $m>1$, there exist points $z_{m1},
\dots ,z_{mm} \in\partial \Omega $ \st
$$
D_m(z_{m1}, \dots ,z_{mm}) \geq D_m(\partial \Omega ) \geq
D_m(z_{m1}, \dots ,z_{mm}) -\frac 1m .
$$
In particular,
$$
D_m(z_{m1}, \dots ,z_{mm}) \to D(\partial \Omega ) \quad \text
{as} \quad m \to \infty .
$$
Let $\mu _m $ be the unit positive regular Borel measure supported
in the subset
$$
\{ z_{m1}, \dots ,z_{mm} \}
$$
of $\partial \Omega $ with
$$
\mu _m (\{ z_{mj} \} )=1/m \quad \text {for} \quad j=1,\dots ,m.
$$
Then, since the sequence $\{ \mu _m (\partial \Omega ) \} $ is
constant (hence bounded), there is a subsequence $\seq \mu {m_k}$
of $\seq \mu m$ converging weakly to a unit positive regular Borel
measure~$\mu $ supported in~$\partial \Omega $. In other words,
$$
\int _{\partial \Omega }fd\mu _{m_k} \to \int _{\partial \Omega
}fd\mu \quad \forall \, f\in C^0(\partial \Omega ).
$$
It follows that, for every \cont function~$f$ on $\partial
\Omega\times  \partial \Omega $,
$$
\int _{\partial \Omega }\int _{\partial \Omega } f(x,y) d\mu
_{m_k}(x) d\mu _{m_k}(y) \to \int _{\partial \Omega }\int
_{\partial \Omega } f(x,y) d\mu _{m}(x) d\mu _{m}(y) .
$$
One can verify this by applying the Stone-Weierstrass theorem to
approximate~$f$ by a linear combination of functions of the form
$a(x)b(y)$ with $a,b \in C^0(\partial \Omega )$.

Given a positive real number~$c$, let $G_c\equiv \min (G,c)$; a
finitely \cont function on $(M\times M) \setminus (\partial
M\times \partial M)$. Then, for each~$k$, we have
\begin{align*}
D_{m_k}(z_{m_k1}, \dots ,z_{m_km_k} ) &\geq \frac
{2}{m_k(m_k-1)}\sum _{1\leq i<j\leq m_k}
G_c(z_{m_ki},z_{m_kj}) \\
&\geq \frac {1}{m_k^2}\sum _{1\leq i,j\leq m_k ; \, i\neq j}
G_c(z_{m_ki},z_{m_kj}) \\
&=\frac {1}{m_k^2}\sum _{i,j=1}^{m_k} G_c(z_{m_ki},z_{m_kj})
-\frac {1}{m_k^2}\sum _{i=1}^{m_k}
G_c(z_{m_ki},z_{m_ki}) \\
&=\frac {1}{m_k^2}\sum _{i,j=1}^{m_k} G_c(z_{m_ki},z_{m_kj})
-\frac {c}{m_k} \\
&=\int _{\partial \Omega }\int _{\partial \Omega } G_c(x,y) \,
d\mu _{m_k}(x) d\mu _{m_k}(y) -\frac {c}{m_k} .
\end{align*}
Thus
$$
D(\partial \Omega ) =\lim _{k \to \infty } D_{m_k}(z_{m_k1},\dots
, z_{m_km_k}) \geq \int _{\partial \Omega }\int _{\partial \Omega
} G_c(x,y) \, d\mu (x) d\mu (y)
$$
for every $c>0$. Therefore, by the monotone convergence theorem,
$$
D(\partial \Omega ) \geq \int _{\partial \Omega }\int _{\partial
\Omega } G(x,y) \, d\mu (x) d\mu (y) =\| \mu \| ^2 \geq W(\partial
\Omega ).
$$
\end{pf}

We now come to the main point:
\begin{lem}\label{Appendix E and D infinite lemma}
$E(\Gamma )=D(\Gamma ) =+\infty $.
\end{lem}
\begin{pf} Let $\seq Mk $ be an \exh of~$M$ by~$\cinf $ \rel
\cpt domains containing~$\overline M_0$. Then, for each positive
integer~$k$,
$$
E(\Gamma )\geq D(\Gamma ) \geq D(\check M \setminus M_k) =
D(\partial M_k)\geq W(\partial M_k) =\frac {(n-2)\tau _{n-1}}{\dir
_{M_k\setminus \overline M_0}(u_k)},
$$
where $u_k$ is the \hm of $\partial M_k$ \wrt $M_k\setminus
\overline M_0$. Since $M$ is parabolic, the sequence~$\seq uk $
converges to~$0$ uniformly on compact sets. Hence
$$
\dir  _{M_k\setminus \overline M_0}(u_k) =\int _{\partial
(M_k\setminus M_0)} u_k\ndof {u_k} =\int _{\partial M_k}\ndof
{u_k} =\int _{\partial M_0}\ndof {u_k} \to 0 \quad \text{as} \quad
k \to \infty .
$$
The lemma now follows.
\end{pf}

\begin{thm}\label{Appendix-Exists harm exh lem} Let $(M,g)$ be a \con non\cpt
oriented parabolic Riemannian manifold of dimension~$n>2$
and~$M_0$ a \rel \cpt domain with $\cinf $ \con boundary. Then
there exists a \cont \exh \fn $\vphi $ on $M\sm M_0$ \st $\vphi $
is \harm on~$M\sm\overline{M_0} $ and equal to~$0$ on~$\partial
M_0$.
\end{thm}
\begin{rmk} Of course, this theorem holds for~$n=2$ as well.
\end{rmk}
\begin{pf*}{Proof of Theorem~\ref{Appendix-Exists harm exh
lem}} In the notation of this section, we have
$$
E_m(\Gamma ) \to E(\Gamma ) =+\infty \quad\text{as}\quad m \to
\infty.
$$
Hence there is an increasing sequence of integers $\seq mk$ \stns
, for each~$k$, $m_k >1$ and $E_{m_k}(\Gamma ) >2^k$. Hence there
exist points $x_{k1}, \dots , x_{km_k}$ in $\Gamma $ \st
$E_{m_k}(\Gamma , x_{k1}, \dots , x_{km_k})>2^k$. Let $\vphi _k$
be the nonnegative \cont function defined by
$$
\vphi _k(x)\equiv \frac {1}{m_k2^k}\sum _{i=1}^{m_k} \check
G(x,x_{ki}) \quad \forall \, x\in \check M\sm M_0.
$$
Then, by Proposition~\ref{Appendix Green's fn extension lemma},
$\vphi _k$ is positive and \harm on~$M\sm\overline{M_0} $ and
equal to~$0$ on~$\partial M_0$. On $\Gamma $, we have
$$
\vphi _k \geq \frac {1}{2^k} E_{m_k}(\Gamma , x_{k1}, \dots ,
x_{km_k}) >1.
$$
Hence there exists an \exh $\seq \Omega k$ of~$M$ by domains \st
$\vphi _k >1$ on~$\check M\setminus \Omega _k$ for each~$k$. If
$K$ is a \cpt subset of~$M\sm M_0$, then $\check G$ is finitely
\contns , hence bounded, on~$K\times \Gamma $. Therefore, for some
positive constant~$a$ (depending on~$K$), we have
$$
0\leq \vphi _k \leq a \sum _{i=1}^{m_k} \frac {1}{m_k2^k} =a2^{-k}
$$
for point $x\in K$ and integer $k>1$. Hence $\sum \vphi _k$
converges uniformly on \cpt subsets of~$M\sm M_0$ to a \cont
function~$\vphi $ which is positive and \harm
on~$M\sm\overline{M_0} $ and equal to~$0$ on~$\partial M_0$.
Moreover, for $x\in M\setminus \Omega _k$, we have
$$
\vphi (x) \geq \sum _{j=1}^k \vphi _j(x) \geq \sum _{j=1}^k 1=k
\to +\infty \quad \text {as} \quad k \to \infty .
$$
Hence $\vphi $ exhausts~$M\sm M_0$.
\end{pf*}

The corresponding result for parabolic ends is applied in Sect. 2
of \cite{NR1}.
\begin{cor}\label{Appendix-Parabolic end cor} Let
$(M,g)$ be a \con non\cpt oriented Riemannian manifold, let
$\Omega $ be a $\cinf $ \rel \cpt domain in~$M$, and let~$E$ be a
\concomp of~$M\setminus \overline \Omega $ which is not \rel \cpt
in~$M$. Assume that~$E$ is a parabolic end. Then there exists a
\cont \exh \fnns~$\psi $ on $\overline E$ which is \harm on~$E$
and equal to~$0$ on~$\partial E$. \rm
\end{cor}
\begin{pf} Clearly, we may assume that~$M$ is the double of~$E$;
a parabolic manifold. Let $M_0$ be a~$\cinf $ \rel \cpt domain in
$\Omega $ \st $M\sm\overline{M_0} $ is \conns , let~$\vphi $ be as
in Theorem~\ref{Appendix-Exists harm exh lem}, and let $\seq Mk$
be an \exh of~$M$ by $\cinf $\rel \cpt domains
containing~$\overline \Omega $, and for each positive integer~$k$,
let $v_k$ be the \cont \fn on $\overline M_k \setminus \Omega $
which is \harm on~$M_k\setminus \overline \Omega $, equal to~$0$
on~$\partial M_k$, and equal to~$\vphi $ on~$\partial \Omega $.
Then $0\leq v_k \leq a \equiv  \max _{\partial \Omega }\vphi $.
Moreover, the sequence $\seq vk $ in nondecreasing and therefore
converges uniformly on \cpt subsets of $M\setminus \Omega $ to a
nonnegative \cont function~$v$ \st $v$ is \harm on~$M\setminus
\overline \Omega $, $v\equiv \varphi $ on~$\partial\Omega$, and
$0\leq v\leq a$ on~$M\setminus\Omega$. Hence the function~$\psi
\equiv \vphi -v$ has the desired properties.
\end{pf}

\begin{rmk} Given a point $x_0\in M_0$, by applying Sario's existence
theorem of principal functions (Theorem~A.1) to the
function~$\vphi $ of Theorem~\ref{Appendix-Exists harm exh lem},
one can construct a unique potential~$p$ on~$M$ such that $p$ is
\harm on~$M\setminus \{ x_0 \} $; $p-G_{x_0}'$ is \harm on~$M_0$,
where $G'$ is the Green's function on~$M_0$; and $p$ exhausts
$M\sm M_0$. The function $p$ is called the {\it Evans-Selberg
potential} of~$M$ with pole at~$x_0$. The existence of the
Evans-Selberg potential is equivalent to parabolicity. For the
details, the reader may refer to~\cite{SaNo}.
\end{rmk}

\bibliographystyle{amsalpha.bst}

\end{document}